\documentclass[final,onefignum]{siamltex}

\usepackage{graphicx}
\usepackage{amsmath}
\usepackage{amssymb}
\usepackage{float}
\usepackage{rotating}

\usepackage[dvipsnames]{xcolor}
\usepackage{multirow}

\newtheorem{remark}{Remark}[section]
\newtheorem{notation}{Notation}
\newtheorem{scheme}{Scheme}

\newtheorem{hypothesis}{Hypothesis}

\title{A stable numerical scheme for stochastic differential equations with multiplicative noise}

\author{
C. M. Mora
\footnotemark[1]
\footnotemark[4]
\and
H. A. Mardones
\footnotemark[1]
\footnotemark[4]
\and
J.C. Jimenez
\footnotemark[2]
\and
M. Selva
\footnotemark[1]
\and
R. Biscay
\footnotemark[3]
}

\begin{document}

\maketitle

\renewcommand{\thefootnote}{\fnsymbol{footnote}}

\footnotetext[1]{Departamento de Ingenier\'{\i}a Matem\'{a}tica, 
Facultad de Ciencias F\'{\i}sicas y Matem\'{a}ticas, Universidad de Concepci\'{o}n,
Casilla 160 C, Concepci\'{o}n, Chile. e-mails: {\tt cmora@ing-mat.udec.cl}, {\tt hmardones@ing-mat.udec.cl}, {\tt mselva@ing-mat.udec.cl}}

\footnotetext[2]{Departamento de Matem\'atica Interdisciplinaria,
Instituto de Cibern\'{e}tica, Matem\'{a}tica y F\'{i}sica. La Habana, Cuba. 
e-mail: {\tt jcarlos@icimaf.cu}
}

\footnotetext[3]{Departamento de Probabilidad y Estad\'istica, Centro de Investigaci\'on en Matem\'atica, Guanajuato, M\'exico.
e-mail: {\tt  rolando.biscay@cimat.mx} 
}

\footnotetext[4]{CI$^2$MA}

\footnotetext[5]{
CMM, HAM and JCJ was partially supported by FONDECYT Grant 1110787.
In addition, CMM and HAM thank the founding of BASAL Grant PFB-03
and  CONICYT Grant  21090691, respectively.}

\renewcommand{\thefootnote}{\arabic{footnote}}

\begin{abstract}
We introduce a new approach for designing 
numerical schemes for stochastic differential equations (SDEs).
The approach, which we have called direction and norm decomposition method, 
proposes to approximate the required solution $X_t$ by integrating the system of coupled SDEs that describes the evolution of the norm of $X_t$ and its projection on the  unit sphere. 
This allows us to develop an explicit scheme for stiff SDEs with multiplicative noise
that shows a solid performance in various numerical experiments.
Under general conditions,
the new integrator preserves the almost sure stability of the solutions for any step-size,
as well as the property of being distant from $0$.
The scheme also has linear rate of  weak convergence 
for a general class of SDEs with locally Lipschitz coefficients,
and one-half  strong order of convergence.
\end{abstract}



\section{Introduction}

This paper deals with the numerical solution of stiff stochastic differential equations (SDEs) with multiplicative noise.
More precisely,
we develop an almost sure stable explicit scheme for the It\^o SDE 
\begin{equation}
\label{eq:1.1}
X_{t}
=
X_{0} 
+ \int_{0}^{t}b\left(  X_{s}\right)  ds+ \sum_{k=1}^{m}\int_{0}^{t}\sigma^{k} \left(X_{s}\right)  dW^{k}_{s} ,
\end{equation}
where 
$W^1,\ldots,W^m$ are independent  Wiener processes on 
a filtered complete probability space 
$\left( \Omega ,\mathfrak{F}, \left(\mathfrak{F}_{t}\right) _{t\geq 0},\mathbb{P}\right) $,
$X_t $ is an adapted $\mathbb{R}^{d}$-valued stochastic process,
and
$b, \sigma^1, \ldots, \sigma^m : \mathbb{R}^{d}\rightarrow\mathbb{R}^{d}$  
have continuous first-order partial derivatives.
In order to do this,
we introduce  a new method for designing numerical schemes for SDEs 
with non-constant diffusion coefficients $\sigma^{k}$. 
We are mainly interested in the computation of
$ \mathbb{E} \varphi \left(X_{t}\right) $,
where $ \varphi \in \mathcal{C}^{\infty} \left( \mathbb{R}^{d},\mathbb{R}\right)$ 
has at most polynomial growth at infinity.

In many cases
the stochastic theta-methods,
like the backward Euler scheme \eqref{scheme:BackwardEuler},
preserve  dynamical  properties of \eqref{eq:1.1}
provided that the step-size $\Delta$  of the time discretization is small enough
(see, e.g., \cite{Higham2007,MaoSINUM2015,Mattingly2010,TalayTubaro1990,Talay2002}).
This does not prevent that the so-called drift-implicit methods (see, e.g., \cite{Milstein1998,Milstein2004})
have poor numerical performance in situations where, for example,  
some partial derivatives of the diffusion coefficients $\sigma^k$  are not small
(see, e.g., \cite{Milstein1998,Milstein2004}).
Though a variety of numerical methods for \eqref{eq:1.1} have been developed in recent times,
the schemes for SDEs with multiplicative noise suffer from step-size restrictions due,
for instance, to stability issues.

Milstein, Platen and Schurz  \cite{Milstein1998} introduced
the general formulation of the balanced implicit methods,
a class of fully implicit schemes for (\ref{eq:1.1})
whose implementation depends on the choice of certain weights
(see, e.g., \cite{Alcock2006,Alcock2012,Milstein1998,Milstein2004,Schurz2005,TretyakovZhang2013}),
which is a complex problem \cite{Schurz2012}.
The reported balanced schemes exhibit low rate of weak convergence,
except incipient progress achieved by \cite{MardonesMora2013PeprintB,Schurz2005}.
Hutzenthaler, Jentzen and Kloeden \cite{Hutzenthaler2012} designed the tamed Euler scheme
for solving \eqref{eq:1.1}  
in case  $b$ satisfies a one-sided Lipschitz condition 
and $\sigma^1,\ldots,\sigma^m$  are globally Lipschitz continuous
(see, e.g., \cite{Hutzenthaler2015,Sabanis2016} for subsequent developments).
Abdulle and Cirilli  \cite{Abdulle2008b} extended Chebyshev's methods 
to solve mean-square stable stiff SDEs (see also, e.g.,  \cite{Abdulle2013}).
Truncating the Brownian motion increments,
Milstein, Repin and Tretyakov \cite{Milstein2002}
constructed a class of fully implicit mean-square schemes 
(see also \cite{Kloeden1992,Milstein2002}).
Multistep, composite and splitting-step methods have been develop, for instance, in 
\cite{AbleidingerBuckwar2016,Anderson2011,BurrageTian2001,MattinglyStuartHigham2002,Milstein2004,Moro2007}.
Using the local linearization method, \cite{Biscay1996} introduced
an exponential scheme for  \eqref{eq:1.1} with $d=1$.
The article \cite{Mora2005} develops an integrator
of  Euler-exponential type for multidimensional SDEs with multiplicative noise
(see also, e.g., \cite{ErdoganLord2016,Komori2014,LordTambue2013,Shoji2011,Stramer99}),
and \cite{JimenezMoraSelva2016} provides a numerical method based on 
the computation of the conditional mean and the square root of the conditional covariance matrix 
of a local linearization approximation to \eqref{eq:1.1}.
Schemes adapted to specific SDEs are given, for instance, in
\cite{Banas2013,Bossy2010,Bossy2015,Mora2005}. 

To the best of our knowledge,
the current numerical methods for \eqref{eq:1.1},
some of them showed up after the first version of this paper \cite{MardonesMora2013PeprintS},
have to use small step-sizes in many cases where the diffusion coefficients 
$\sigma^1,\ldots,\sigma^m$ play an essential role in the dynamics of $X_t$
(see, e.g., Section \ref{Sec:Numerical Experiments}).
This motivates the introduction in Subsection \ref{subsec:EquilibriumAt0} of 
a new technique  for constructing almost sure stable methods for \eqref{eq:1.1}
with equilibrium point at $0$.
In  case $b\left(0\right)=\sigma^1\left(0\right)=\cdots \sigma^m \left(0\right)=0$, 
we divide the numerical approximation of $X_t$ into the computation of
$ Z_t := X_{t} /  \left\Vert X_{t} \right\Vert $ and $ \left\Vert X_{t} \right\Vert $.
Since  $ Z_t$ and $\left\Vert X_{t} \right\Vert $ satisfy 
the system  (\ref{eq:2.3}) and (\ref{eq:2.2}) given below,
we propose to simulate $X_t$
by solving numerically these  two coupled SDEs with smooth coefficients.
We call this approach  Direction and Norm Decomposition method (DND).
We take advantage of the unit norm property of  $ Z_t$
and the one-dimensionality of \eqref{eq:2.3}.
Namely,
this article computes  $ Z_t$ by projecting the Euler-Maruyama scheme applied to  (\ref{eq:2.2}) onto the unit sphere,
and 
the norm of $X_t$ is obtained by applying an exponential scheme to the scalar SDE (\ref{eq:2.3}).
This yields Scheme \ref{scheme:EulerStable} given below. 
Subsection \ref{subsec:GeneralSDEs} provides a way to extend Scheme \ref{scheme:EulerStable} 
to the framework where 
$b\left(0\right), \sigma^1\left(0\right), \cdots$, $\sigma^m\left(0\right)$ 
may be different from $0$.
By adding an auxiliary function,
we transform \eqref{eq:1.1} into a 
$\mathbb{R}^{d+1}$-stochastic differential equation with equilibrium point at $0$.
Then, applying Scheme \ref{scheme:EulerStable},
along with suitable approximations,
we get  Scheme \ref{scheme:EulerStableG}.
It is worth pointing out that Scheme \ref{scheme:EulerStableG} 
becomes Scheme \ref{scheme:EulerStable} 
whenever $b\left(0\right)=\sigma^1\left(0\right)=\cdots \sigma^m \left(0\right)=0$.
Section \ref{Sec:Numerical Experiments} presents various numerical experiments
that illustrate the very good  performance of the new scheme
even for large step-sizes.

Suppose for the moment that $0$ is an equilibrium point of \eqref{eq:1.1}.
In case $b,\sigma^1, \ldots$, $\sigma^m$ have at most linear growth and
\begin{equation}
\label{eq:StabCond}
- \lambda :=
\sup_{ x \neq 0} 
\left(
\frac{
\langle x, b \left( x \right)  \rangle
+
\frac{1}{2} \sum_{k=1}^{m}  \left\Vert \sigma^k \left( x \right)  \right\Vert^2 
 }{
 \left\Vert x \right\Vert^2 
 }
- \frac{ 
\sum_{k=1}^{m}
\left\langle x,  \sigma^k \left( x \right)  \right\rangle^2
}{
\left\Vert  x   \right\Vert^4
}
\right)
< 0 ,
\end{equation}
Higham, Mao and Yuan \cite{Higham2007} prove 
the almost sure exponential stability of the Euler-Maruyama method
for small enough step-sizes (see also \cite{MaoSINUM2015}).
In this situation,  
\eqref{eq:1.1} is almost sure exponential stable. 
If the linear growth condition of $b$ is replaced by the one-sided Lipschitz condition
and \eqref{eq:StabCond} is substituted by a slightly stronger requirement,
then the backward Euler method
\begin{equation}
 \label{scheme:BackwardEuler}
 \bar{E}_{n+1}
=
\bar{E}_{n} +  b\left(   \bar{E}_{n+1} \right)    \Delta
+ \sum_{k=1}^{m} \sigma^{k} \left( \bar{E}_{n} \right) \left( W^k_{T_{n+1}} - W^k_{T_{n}} \right)
\end{equation}
is  almost sure exponential stable for sufficiently small step-sizes $\Delta$ (see \cite{Higham2007}).
Here and subsequently,
$T_k = k \, \Delta$ for any $k \in \mathbb{Z}_+$, where $\Delta > 0$.
The almost sure asymptotic stability of 
the stochastic theta-methods applied to test SDEs
has been studied, for instance, in \cite{Berkolaiko2012,Buckwar2012,Higham2000,MaoSINUM2015,RodkinaSchurz2003,SaitoMitsui1993,Schurz2012}.
Similarly,
the almost sure asymptotic stability of Balanced schemes has been tested, for example, in \cite{Alcock2012,MardonesMora2013PeprintB,Schurz2005,Schurz2012,Szpruch2015}.
Under \eqref{eq:StabCond} and the linear growth condition of $\sigma^1,\ldots,\sigma^m$,
we obtain 
the almost sure exponentially stability of Scheme \ref{scheme:EulerStable} for any step-size $\Delta > 0$
(see Theorem  \ref{th:Stable} of Subsection \ref{subsec:LongTime}). 
To the best of our knowledge,
Scheme \ref{scheme:EulerStable} is the first numerical method that preserves, 
in a large class of SDEs, 
the almost sure asymptotic stability of \eqref{eq:1.1} no matter the value of the step-size.
In case $b, \sigma^1, \ldots,  \sigma^m$ are globally Lipschitz functions,
Mao \cite{MaoSINUM2015} proved that \eqref{eq:1.1} is small-moment exponentially stable 
iff the stochastic theta methods applied to \eqref{eq:1.1} are small-moment exponentially stable 
for sufficiently small step sizes (see also  \cite{Higham2007}).
We recall that  \eqref{eq:1.1} is small-moment exponentially stable 
iff \eqref{eq:1.1} is $p$th moment exponentially stable for a sufficiently small $p \in \left] 0, 1 \right[$.
In Subsection \ref{subsec:LongTime}, 
we show that 
Scheme \ref{scheme:EulerStable}  preserves the small-moment exponential stability of \eqref{eq:1.1}
for any step-size whenever \eqref{eq:StabCond} holds
(see Theorem \ref{th:MomentStable}).

It is important that the numerical solution of \eqref{eq:1.1} captures 
the behavior of $X_t$ when $0$ is a non-stable fixed point of  \eqref{eq:1.1}.
In this direction,
we prove that Scheme \ref{scheme:EulerStable}
preserves the non-stability of the origin for any step-size 
provided that $b$, $\sigma^1,\ldots,\sigma^m$ satisfy 
a general criterion for $0$ being a  non-stable  equilibrium point of \eqref{eq:1.1}
(see Theorem \ref{th:NoStable} below).
Previously,
\cite{Berkolaiko2012} has verified that
the stochastic theta methods applied to a test SDE
reproduce the almost sure instability of $0$ when the step-sizes are sufficiently small.
On the other hand, 
Scheme \ref{scheme:EulerStable} also keeps intact  the sign of $X_0$ 
in case $d=1$,
which is an interesting property (see, e.g., \cite{Bossy2015,Moro2007}).

Many applications deal with the computation of $ \mathbb{E} \varphi \left(X_{t}\right) $,
with $\varphi : \mathbb{R}^d \rightarrow \mathbb{R}$.
This motivates the study of the weak errors, i.e., 
the difference between $ \mathbb{E} \varphi \left(X_{t}\right) $
and the expectation of the approximate value of $\varphi \left(X_{t}\right) $.
Using the Kolmorogov equation,
Talay \cite{Talay1984,Talay1986} and Milstein \cite{Milstein1985}
developed a methodology for obtaining the rate of weak convergence 
of  the numerical schemes for  \eqref{eq:1.1}
(see also, e.g., \cite{GrahamTalay2013,Milstein2004,Talay1990} and \cite{ClementKohatsuLamberton2006}).
Thus, 
\cite{Milstein1985,Talay1984,Talay1986} got the linear weak convergence rate
of the Euler Maruyama scheme under the global Lipschitz condition.
Few weak convergence results are available for SDEs with non-globally Lipschitz continuous coefficients,
class of SDEs that appears in important applications.
In case $b,\sigma^1,\ldots,\sigma^m$ are locally Lipschitz functions,
Milstein and Tretyakov \cite{MilsteinTretyakov2005}
proposed to discard the numerical trajectories leaving  a sufficiently large sphere
and studied the weak error involved in this  procedure.
Hairer, Hutzenthaler and Jentzen \cite{Hairer2015} showed 
the existence of a locally Lipschitz SDE with smooth coefficients 
for which the Euler-Maruyama converges in the weak sense
(also in the strong one) without any arbitrarily small polynomial rate of convergence
(see also \cite{Hutzenthaler2015,Hutzenthaler2011}).
Bossy and Diop \cite{Bossy2015} obtained that the symmetrized {E}uler scheme 
attains the order $1$ of weak convergence 
when it is applied to \eqref{eq:1.1}
with  $d=m=1$, $\sigma^1 \left( x \right) = \left\vert x\right\vert^{\alpha}$
and $b$ globally Lipschitz continuous, where $\alpha\in\left[1/2,1\right[$.
In the ergodic case,
Talay \cite{Talay2002} addressed the computation of integrals with respect to the invariant probability law of 
a kind of stochastic Hamiltonian system, with additive noise and $b$ locally Lipschitz continuous;
Talay \cite{Talay2002} showed that the discretization error of the backward Euler scheme
has the same expansion as in the globally Lipschitz case \cite{TalayTubaro1990}.
We prove that Scheme \ref{scheme:EulerStableG} converges weakly with order $1$
under a global coercivity condition 
and the smoothness of the solution of the Kolmogorov equation associated to \eqref{eq:1.1}
(see Theorem \ref{th:RateConvergence} below),
which is a general class of SDEs with multiplicative noise.
To this end,
we derive a fundamental weak convergence theorem (see Theorem \ref{th:GeneralWeakConvergence}).

Tretyakov and Zhang \cite{TretyakovZhang2013}
gave a fundamental mean-square convergence theorem 
for SDEs with non-global Lipschitz coefficients that satisfy  a global monotonicity condition.
Moreover, 
\cite{TretyakovZhang2013} introduced a particular balanced scheme 
which has rate $1/2$ of strong convergence in a non-global Lipschitz setting
(see, e.g., \cite{Hairer2015,Hutzenthaler2015,MaoSzpruch2013,TretyakovZhang2013} 
for a recent account of  strong convergence results for SDEs with non-global Lipschitz coefficients).
We prove that Scheme \ref{scheme:EulerStableG} converge in 
$L^p \left( \mathbb{P} \right)$ with rate $1/2$, where $p \in \mathbb{N}$,
under the assumptions of the fundamental mean-square convergence theorem proved in \cite{TretyakovZhang2013}
(see Theorem \ref{th:StrongRateConvergence} below),

The paper is organized as follows.
In Section \ref{sec:2} we introduce the direction and norm decomposition method (DND).
Section \ref{sec:Convergence} is devoted to the stability and convergence properties of the DND scheme.
Section \ref{Sec:Numerical Experiments} provides numerical experiments.
All proofs are deferred to Section \ref{Sec:Proofs},
and Section \ref{sec:Conclusions} presents our conclusions.

\subsection{Notation}
For simplicity, we consider the equidistant  time discretization
$T_n = n \Delta$, where $\Delta >0$ and $n=0,1,\ldots$
We will use the same symbol $K$ (resp. $q$ and $K\left( \cdot\right) $)
for different non-negative real numbers 
(resp. natural numbers and non-negative increasing functions)
that have the common property to be independent of  $\Delta$.
We set $\mathcal{P}_{0}=\left\{ 0\right\}$,
and $\partial^{0}$ denotes the identity operator.
For any $\ell\in\mathbb{N}$
we define  
$\mathcal{P}_{\ell} = \left\{1,\ldots,d\right\} ^{\ell}$,
and we write 
$
\partial_{x}^{\vec{p}}
= 
\frac{\partial}{\partial x_{p_{1}}} \cdots\frac{\partial}{\partial x_{p_{\ell}}}
$
whenever $\vec{p} = \left( p_1, \ldots, p_{\ell} \right) \in\mathcal{P}_{\ell}$.
Let 
$\mathcal{C}_{P}^{L}\left( \left[ 0, T \right] \times \mathbb{R}^{d},\mathbb{R}\right) $
be the set of all 
$
f  :
\left[ 0, T \right] \times \mathbb{R}^{d}\rightarrow \mathbb{R}
$ 
such that for any  $\vec{p}\in\mathcal{P}_{\ell}$, with $\ell \leq L$,
(i)  $\partial_{x}^{\vec{p}} f $ is continuous;
and 
(ii) 
$\left\vert \partial_{x}^{\vec{p}} f \left( t,x\right)
\right\vert \leq K  \left( 1+\left\Vert x\right\Vert ^{q}\right) $ for
all $t \in \left[ 0, T \right]$ and $x\in\mathbb{R}^{d}$.
Here and below,
$ 
\left\Vert \cdot \right\Vert 
$
and 
$ \langle \cdot ,  \cdot \rangle $ 
stand for the norm and the dot product 
(the usual Euclidean scalar product)
on $\mathbb{R}^{d}$, respectively.
We say that
$
f : \mathbb{R}^{d} \rightarrow \mathbb{R}
$
is in the class 
$
\mathcal{C}_{P}^{L}\left(\mathbb{R}^{d},\mathbb{R}\right)
$
if
$
\left( t, x \right) \mapsto f \left( x \right)
$
belongs to
$
\mathcal{C}_{P}^{L}\left( \left[ 0, T \right] \times \mathbb{R}^{d},\mathbb{R}\right)
$.
The Jacobian matrix of 
$ g: \mathbb{R}^d \rightarrow \mathbb{R}^d$
is denoted by $Jg$.

%
%


\section{Direction and norm decomposition method}
\label{sec:2}

Since $b$ and $\sigma^1, \ldots, \sigma^m$ are locally Lipschitz,
(\ref{eq:1.1})  has a unique continuous strong solution up to an explosion time 
(see, e.g., \cite{Protter2005}),
which we assume to be $+ \infty$ a.s.
This happens, for instance, when 
$
\langle x, b \left( x \right)  \rangle
+
\frac{1}{2} \sum_{k=1}^{m}  \left\Vert \sigma^k \left( x \right)  \right\Vert^2
\leq 
K \left( 1 +   \left\Vert  x \right\Vert^2 \right) 
$
for all  $x \in \mathbb{R}^d$
(see, e.g., \cite{Gihman1972,Mao2007}).

\subsection{SDEs with equilibrium at $0$}
\label{subsec:EquilibriumAt0}

Suppose that $b\left(0\right)=\sigma^1\left(0\right)=\cdots=\sigma^m\left(0\right)=0$.
Then, there is no loss of generality in assuming $X_0 \neq 0$ a.s.,
and so,  almost surely, $X_t \neq 0$ for all $t \geq 0$.
In this case,
we propose to  divide the computation of $X_t$ into the numerical approximations of 
$ Z_t := X_{t} /  \left\Vert X_{t} \right\Vert $ and the norm of 
$X_{t}$.

We start by obtaining the SDEs describing the evolution of 
$ \left\Vert X_{ t } \right\Vert $ and $X_{t} /  \left\Vert X_{t} \right\Vert $.
Applying It\^o's formula to $\sqrt{\left\| X_{t \wedge \widetilde{\tau}_j} \right\|^2}$ we obtain
\begin{align*}
 \left\Vert X_{t \wedge \widetilde{\tau}_j} \right\Vert
 & =
\left\Vert X_{0} \right\Vert
+ \sum_{k=1}^m \int_{0}^{t \wedge \widetilde{\tau}_j} \frac{\langle X_s,\sigma^k \left(X_s\right)\rangle}{\left\Vert X_s \right\Vert}  dW^k_s 
\\
& \quad
+\int_{0}^{t \wedge \widetilde{\tau}_j} \left(\frac{\langle X_s,b\left(X_s\right)\rangle 
+\frac{1}{2} \sum_{k=1}^m \left\Vert\sigma^k \left(X_s\right)\right\Vert^2}{\left\Vert X_s \right\Vert}
- \frac{1}{2} \sum_{k=1}^m \frac{ \langle X_s,\sigma^k \left(X_s\right)\rangle^2}{\left\Vert X_s \right\Vert^3} \right) ds ,
\end{align*}
where 
$ \widetilde{\tau}_j := \inf \left\{ t > 0 : \left\| X_t \right\| < 1/j  \right\}   $. 
Since $X_t \neq 0$  for all $t \geq 0$,
$ \widetilde{\tau}_j \longrightarrow_{j \rightarrow \infty} \infty $,
and so taking limit as  $j \rightarrow \infty$ gives 
\begin{align}
\label{eq:2.1}
\left\Vert X_{t} \right\Vert
&=
\left\Vert X_{0} \right\Vert
+ \sum_{k=1}^m \int_{0}^{t} \frac{\langle X_s,\sigma^k \left(X_s\right)\rangle}{\left\Vert X_s \right\Vert}  dW^k_s 
\\
\nonumber
& \quad
+\int_{0}^{t} \left(\frac{\langle X_s,b\left(X_s\right)\rangle 
+\frac{1}{2} \sum_{k=1}^m \left\Vert\sigma^k \left(X_s\right)\right\Vert^2}{\left\Vert X_s \right\Vert}
- \frac{1}{2} \sum_{k=1}^m \frac{ \langle X_s,\sigma^k \left(X_s\right)\rangle^2}{\left\Vert X_s \right\Vert^3} \right) ds  .
\end{align}
We rewrite (\ref{eq:2.1}) as 
\begin{align}
\label{eq:2.3}
 \left\Vert X_{t} \right\Vert
  =
\left\Vert X_{0} \right\Vert
+ \int_{0}^{t} \beta \left( \left\Vert X_s \right\Vert,  Z_s \right)
\left\Vert X_s \right\Vert ds
+ \sum_{k=1}^m \int_{0}^{t} \left\langle Z_s,\bar{\sigma}^k \left( \left\Vert X_s \right\Vert,  Z_s \right) \right\rangle 
\left\Vert X_s \right\Vert dW^k_s ,
\end{align}
where for any $\eta \in \mathbb{R}$ and $z \in \mathbb{R}^d$ we define
$$
\beta \left( \eta , z \right)
 =
\left\langle z , \bar{b} \left( \eta , z \right) \right\rangle 
+\frac{1}{2} \sum_{k=1}^m \left\Vert\  \bar{\sigma}^k \left(\eta ,  z \right) \right\Vert^2
- \frac{1}{2} \sum_{k=1}^m \left\langle z , \bar{\sigma}^k \left(\eta ,  z \right) \right\rangle^2 ,
$$
\begin{equation}
\label{eq:3.5}
 \bar{b} \left( \eta, z \right)
= 
\begin{cases}
b \left( \eta \, z \right) / \eta
&
\text{if } \eta \neq 0 ,
\\
Jb  \left( 0 \right) z
&
\text{if } 
\eta = 0
\end{cases}
\text{ and }
\bar{ \sigma }^k \left( \eta, z \right)
= 
\begin{cases}
\sigma^k \left( \eta \, z \right) / \eta
&
\text{if } \eta \neq 0 ,
\\
J \sigma^k  \left( 0 \right) z
&
\text{if } 
\eta = 0 .
\end{cases}
\end{equation}

Applying It\^o's formula to 
$ X_{t \wedge \widetilde{\tau}_j} / \sqrt{\left\| X_{t \wedge \widetilde{\tau}_j} \right\|^2}$
we get after a long calculation that 
\begin{align*}
\frac{X_{t\wedge \widetilde{\tau}_j}}{\left\Vert X_{t\wedge \widetilde{\tau}_j} \right\Vert}
& = 
\frac{X_{0}}{\left\Vert X_{0} \right\Vert}
+\int_{0}^{t\wedge \widetilde{\tau}_j} \left(\frac{b\left(X_s\right)}{\left\Vert X_s \right\Vert}
- \left\langle \frac{X_s}{\left\Vert X_s \right\Vert}, \frac{b\left(X_s\right)}{\left\Vert X_s \right\Vert}\right \rangle \frac{X_s}{\left\Vert X_s \right\Vert} \right) ds 
\\
\nonumber
& \quad
+\frac{1}{2} \sum_{k=1}^m \int_{0}^{t\wedge \widetilde{\tau}_j}  \left(3 \left\langle \frac{X_s}{\left\Vert X_s \right\Vert},\frac{\sigma^k\left(X_s\right)}{\left\Vert X_s \right\Vert}\right\rangle^2
-\left\langle \frac{\sigma^k\left(X_s\right)}{\left\Vert X_s \right\Vert},\frac{\sigma^k\left(X_s\right)}{\left\Vert X_s \right\Vert}\right\rangle \right)  \frac{X_s}{\left\Vert X_s \right\Vert} ds 
\\
\nonumber
& \quad
- \sum_{k=1}^m \int_{0}^{t\wedge \widetilde{\tau}_j} 
 \left\langle \frac{X_s}{\left\Vert X_s \right\Vert},\frac{\sigma^k\left(X_s\right)}{\left\Vert X_s \right\Vert}\right\rangle \frac{\sigma^k\left(X_s\right)}{\left\Vert X_s \right\Vert} 
ds 
\\
\nonumber
& \quad
+\sum_{k=1}^m\int_{0}^{t\wedge \widetilde{\tau}_j} 
\left(\frac{\sigma^k\left(X_s\right)}{\left\Vert X_s \right\Vert} - \left\langle \frac{X_s}{\left\Vert X_s \right\Vert},\frac{\sigma^k\left(X_s\right)}{\left\Vert X_s \right\Vert}\right\rangle \frac{X_s}{\left\Vert X_s \right\Vert} \right)dW_s^k .
\end{align*}
This implies
\begin{align}
\nonumber
Z_t
& = 
Z_{0}
+\int_{0}^t \left( \bar{ b} \left( \left\Vert X_s \right\Vert,  Z_s \right) 
- \left\langle Z_s , \bar{ b } \left( \left\Vert X_s \right\Vert,  Z_s \right)  \right\rangle  Z_s 
+ \Psi  \left( \left\Vert X_s \right\Vert,  Z_s \right)
\right) ds 
\\
\label{eq:2.2} 
& \quad
+\sum_{k=1}^m\int_{0}^t \left( \bar{ \sigma }^k \left( \left\Vert X_s \right\Vert,  Z_s \right) 
- \left\langle Z_s ,\bar{ \sigma }^k \left( \left\Vert X_s \right\Vert,  Z_s \right)\right\rangle Z_s \right)dW_s^k ,
\end{align}
where $Z_t$ was defined by $ Z_t = X_{t} /  \left\Vert X_{t} \right\Vert $ and 
\begin{equation}
\label{eq:2.21} 
\Psi  \left( \eta, z \right)
\hspace{-1pt} = \hspace{-1pt}
\sum_{k=1}^m  \hspace{-2pt} \left( 
 \left( \frac{3}{2} \left\langle z , \bar{ \sigma }^k \left( \eta ,  z \right)
\right\rangle^2
-
\frac{1}{2} \left\Vert  \bar{ \sigma}^k \left( \eta ,  z \right) \right\Vert ^2 \hspace{-1pt}  \right)  z
-
\left\langle z , \bar{ \sigma }^k \left(  \eta , z \right) \right\rangle \bar{ \sigma }^k \left(  \eta ,  z \right)
 \right) \hspace{-1pt} .
\end{equation}

Since 
$
X_{t}
=
 \left\Vert X_t \right\Vert Z_t
$,
the main idea of this paper is to compute $ X_{t} $
by solving the system of SDEs formed by (\ref{eq:2.3}) and (\ref{eq:2.2}),
which defines the Direction and Norm Decomposition method (DND).
Taking advantage of 
the unit norm property of  $ Z_t = X_{t} /  \left\Vert X_{t} \right\Vert $
and the one-dimensionality of (\ref{eq:2.3}),
we next design a simple numerical scheme based on our DND method.

Suppose that the $\mathfrak{F}_{0}$-measurable random variable $\bar{X}_0 \neq 0$
simulates the initial condition $X_0$.
Then, we set 
$ \bar{\eta}_{0} = \left\Vert \bar{X}_0 \right\Vert$ and $\hat{Z}_{0} = \bar{X}_0  / \left\Vert \bar{X}_0 \right\Vert$.
In order to compute  
$\left( \left\Vert X_{T_n} \right\Vert \right)_{n \in \mathbb{Z}_+}$ and
$ \left( Z_{T_n}  \right)_{n \in \mathbb{Z}_+}$,
we next generate recursively the pairs of $\mathfrak{F}_{T_n}$-measurable random variables
$\left( \bar{\eta}_{n}, \hat{Z}_n \right)$ 
such that  
$\bar{\eta}_{n}$ takes values in $ \left[ 0 , \infty \right[$
and 
$ \hat{Z}_n$ lies on the unit sphere of $\mathbb{R}^d$
for any $n \geq 0$.
Here,
$\bar{\eta}_{n}$ and $\hat{Z}_n$ will approximate $\left\Vert X_{T_n} \right\Vert $ and $  Z_{T_n}$, respectively.
Fix $\bar{\eta}_{n} $ and $\hat{Z}_n$,
which satisfy $\bar{\eta}_{n} \geq 0$ and $\left\Vert \hat{Z}_n \right\Vert = 1$.
From \eqref{eq:2.1} it follows that for all $t \in \left[ T_n , T_{n+1} \right]$, 
$$
 \left\Vert X_{t} \right\Vert
  =
\left\Vert X_{T_n} \right\Vert
+ \int_{T_n}^{t} \beta \left( \left\Vert X_s \right\Vert,  Z_s \right)
\left\Vert X_s \right\Vert ds
+ 
\sum_{k=1}^m \int_{T_n}^{t} \left\langle Z_s,\bar{\sigma}^k \left( \left\Vert X_s \right\Vert,  Z_s \right) \right\rangle 
\left\Vert X_s \right\Vert dW^k_s .
$$
Freezing $\bar{ b } \left( \left\Vert X_s \right\Vert,  Z_s \right) $,
$\bar{  \sigma }^1 \left( \left\Vert X_s \right\Vert,  Z_s \right) , \ldots, 
\bar{  \sigma }^m \left( \left\Vert X_s \right\Vert,  Z_s \right)$
over $s \in \left[ T_n , T_{n+1} \right]$ at  the values 
$\bar{ b } \left( \left\Vert X_{T_n} \right\Vert,  Z_{T_n} \right) $,
$\bar{  \sigma }^1 \left( \left\Vert X_{T_n} \right\Vert,  Z_{T_n} \right), \ldots,  
\bar{  \sigma }^m \left( \left\Vert X_{T_n} \right\Vert,  Z_{T_n} \right)$,
and replacing the pair
$\left\Vert X_{T_n} \right\Vert$, $Z_{T_n}$
by $\bar{\eta}_n$,  $\hat{Z}_n$, 
we obtain the linear scalar SDE
\begin{equation*}
 \eta_t
 =
\bar{\eta}_n 
+ \int_{T_n}^{t} \beta \left( \bar{\eta}_n , \hat{Z}_n \right) \eta_s ds
+ \sum_{k=1}^m \int_{T_n}^{t} \left\langle \hat{Z}_n  ,\bar{  \sigma }^k \left(\bar{\eta}_n  , \hat{Z}_n \right) \right\rangle \eta_s dW^k_s ,
\end{equation*}
whose solution at time $T_{n+1}$ is
\begin{equation}
\label{eq:2.22}
 \begin{aligned}
\eta_{T_{n+1}} = \bar{\eta}_n  \exp  \left(
 \left( \beta  \left( \bar{\eta}_n , \hat{Z}_n \right) 
 - \frac{1}{2} \sum_{k=1}^m \left\langle \hat{Z}_n  ,\bar{ \sigma}^k \left(\bar{\eta}_n  , \hat{Z}_n \right) \right\rangle^2  \right) \Delta
 \right. \hspace{1cm} &
 \\
 \left. 
 + \sum_{k=1}^m \left\langle \hat{Z}_n  ,\bar{ \sigma}^k \left(\bar{\eta}_n  , \hat{Z}_n \right) \right\rangle 
 \left(W^k_{T_{n+1}} - W^k_{T_n} \right)
 \right) . &
\end{aligned}
\end{equation}

We can solve (\ref{eq:2.2}) using various iterative schemes.
For simplicity,
we apply the Euler approximation to
\begin{align*}
Z_t
& = 
Z_{T_n}
+\int_{T_n}^t \left( \bar{ b} \left( \left\Vert X_s \right\Vert,  Z_s \right) 
- \left\langle Z_s , \bar{ b } \left( \left\Vert X_s \right\Vert,  Z_s \right)  \right\rangle  Z_s 
+ \Psi  \left( \left\Vert X_s \right\Vert,  Z_s \right)
\right) ds 
\\
& \quad
+\sum_{k=1}^m\int_{T_n}^t \left( \bar{ \sigma }^k \left( \left\Vert X_s \right\Vert,  Z_s \right) 
- \left\langle Z_s ,\bar{ \sigma }^k \left( \left\Vert X_s \right\Vert,  Z_s \right)\right\rangle Z_s \right)dW_s^k ,
\end{align*}
which yields 
\begin{align*}
Z_{T_{n+1}} 
& \approx 
Z_{T_n} + \left( 
 \bar{ b } \left(  \left\Vert X_{T_n} \right\Vert , Z_{T_n} \right) 
- \left\langle  Z_{T_n} , \bar{ b } \left(  \left\Vert X_{T_n} \right\Vert , Z_{T_n} \right) \right\rangle  Z_{T_n}
+ \Psi  \left(  \left\Vert X_{T_n} \right\Vert , Z_{T_n} \right) 
\right) \Delta
\\
& \quad
+ \sum_{k=1}^m  \left( \bar{ \sigma }^k \left(  \left\Vert X_{T_n} \right\Vert  , Z_{T_n} \right)  
- \left\langle  Z_{T_n} ,  \bar{ \sigma }^k \left(  \left\Vert X_{T_n} \right\Vert  , Z_{T_n} \right) \right\rangle Z_{T_n} \right)
\left( W_{T_{n+1}}^k-W_{T_{n}}^k \right) .
\end{align*}
Hence,
substituting  
$\left\Vert X_{T_n} \right\Vert$ and $Z_{T_n}$
by $\bar{\eta}_n$ and  $\hat{Z}_n$, respectively, we obtain
\begin{equation}
\label{eq:2.23}
 \begin{aligned}
\widetilde{Z}_{n+1} 
 & : =
 \hat{Z}_n + \left( 
 \bar{ b } \left( \bar{\eta}_n , \hat{Z}_n \right) 
- \left\langle  \hat{Z}_n , \bar{ b } \left( \bar{\eta}_n , \hat{Z}_n \right) \right\rangle  \hat{Z}_n
+ \Psi  \left( \bar{\eta}_n , \hat{Z}_n \right) 
\right) \Delta
\\
& \quad
+ \sum_{k=1}^m  \left( \bar{ \sigma }^k \left(\bar{\eta}_n  , \hat{Z}_n \right)  
- \left\langle  \hat{Z}_n ,  \bar{ \sigma }^k \left(\bar{\eta}_n  , \hat{Z}_n \right) \right\rangle  \hat{Z}_n \right)
\left( W_{T_{n+1}}^k-W_{T_{n}}^k \right) . 
\end{aligned}
\end{equation}

As we are interested in simulating the distribution of $\left(  X_{T_{n+1}} , Z_{T_{n+1}}\right)$,
in \eqref{eq:2.22} and \eqref{eq:2.23} we replace  $W_{T_{n+1}}^k-W_{T_{n}}^k$ by $ \sqrt{\Delta} \hat{W}^k_{n+1}$,
where $\hat{W}^1_{n+1}, \ldots, \hat{W}^m_{n+1}$ are 
independent and identically distributed (i.i.d.) $\mathfrak{F}_{T_{n+1}}$-measurable random variables with symmetric law and variance $1$,
which are also independent of $\mathfrak{F}_{T_{n}}$.
Then $\eta_{T_{n+1}}$ and $\widetilde{Z}_{n+1} $ become
\begin{align*}
 \bar{\eta}_{n+1} := \bar{\eta}_n  \exp  \left( 
 \left\langle \hat{Z}_n , \bar{ b } \left( \bar{\eta}_n , \hat{Z}_n \right) \right\rangle  \Delta
+
 \sum_{k=1}^m \left\langle \hat{Z}_n  ,\bar{ \sigma }^k \left(\bar{\eta}_n  , \hat{Z}_n \right) \right\rangle \sqrt{\Delta} \hat{W}^k_{n+1} 
 \right.
 \hspace{1cm}
 &
 \\
 \left.
+ \left( 
 \frac{1}{2} \sum_{k=1}^m \left\Vert\  \bar{  \sigma }^k \left( \bar{\eta}_n ,  \hat{Z}_n \right) \right\Vert^2
-  \sum_{k=1}^m \left\langle \hat{Z}_n , \bar{ \sigma }^k \left( \bar{\eta}_n ,  \hat{Z}_n \right) \right\rangle^2
  \right) \Delta
 \right)
 &
\end{align*}
and
\begin{align*}
 \bar{Z}_{n+1}
& :=
 \hat{Z}_n + \left( 
 \bar{ b } \left( \bar{\eta}_n , \hat{Z}_n \right) 
- \left\langle  \hat{Z}_n , \bar{ b } \left( \bar{\eta}_n , \hat{Z}_n \right) \right\rangle  \hat{Z}_n
+ \Psi  \left( \bar{\eta}_n , \hat{Z}_n \right) 
\right) \Delta
\\
\nonumber
& \quad
+ \sum_{k=1}^m  \left( \bar{ \sigma }^k \left(\bar{\eta}_n  , \hat{Z}_n \right)  
- \left\langle  \hat{Z}_n ,  \bar{ \sigma }^k \left(\bar{\eta}_n  , \hat{Z}_n \right) \right\rangle  \hat{Z}_n \right)
 \sqrt{\Delta} \hat{W}^k_{n+1} .
\end{align*}

Finally, we take  $\hat{Z}_{n+1}:= \bar{Z}_{n+1} /  \left\Vert  \bar{Z}_{n+1} \right\Vert$,
and so we have defined the new pair $\left( \bar{\eta}_{n+1}, \hat{Z}_{n+1} \right)$.
The fact that  
$\left\Vert Z_{T_{n+1}} \right\Vert = 1$
allows us to improve the accuracy of $\bar{Z}_{n+1}$ 
by projecting $\bar{Z}_{n+1}$ on the unit sphere,
a normalization procedure that
has been used with success
in the numerical solution of the non-linear Sch\"odinger equations 
(see, e.g., \cite{Mora2005,Percival1998})
and the computation of Lyapunov exponents (see, e.g., \cite{Carbonell2010,Talay1991}).
Since $\left\Vert \bar{Z}_{n+1} \right\Vert$ approximates the law of $\left\Vert Z_{T_{n+1}} \right\Vert = 1$,
we can expect that $\left\Vert \bar{Z}_{n+1} \right\Vert $ is not close to $0$,
and so the weak approximation $\hat{Z}_{n+1}= \bar{Z}_{n+1} /  \left\Vert  \bar{Z}_{n+1} \right\Vert$
of $ Z_{T_{n+1}}$
should reproduce efficiently  the unit-norm property of  $Z_{T_{n+1}}$
without incurring round-off errors.
In summary,
we have designed the following numerical scheme.


\begin{scheme}
 \label{scheme:EulerStable}
 Let $\bar{\eta}_{0}  \in \mathbb{R}$  and $\hat{Z}_0 \in \mathbb{R}^d$
be random variables  
satisfying $\bar{\eta}_{0} \geq 0$ and $\left\Vert \hat{Z}_0 \right\Vert =1$.
Consider the i.i.d. symmetric random variables
$\hat{W}^1_1, \hat{W}^2_1, \ldots, \hat{W}^m_1, \hat{W}^1_2, \ldots$ 
with variance $1$
that are independent of  $\bar{\eta}_{0}, \hat{Z}_0 $.
For any $\Delta > 0$, 
we define recursively the pair $\left( \bar{\eta}_{n+1} ,  \hat{Z}_{n+1} \right)$ by 
\begin{equation}
\label{eq:2.5}
\bar{\eta}_{n+1} := \bar{\eta}_n  \exp  \left( 
\mu_n \Delta
+
 \sum_{k=1}^m \left\langle \hat{Z}_n  ,\bar{ \sigma }^k \left(\bar{\eta}_n  , \hat{Z}_n \right) \right\rangle \sqrt{\Delta} \hat{W}^k_{n+1} 
 \right)
\end{equation}
and
$
\hat{Z}_{n+1} =
\begin{cases}
 \bar{Z}_{n+1} /  \left\Vert  \bar{Z}_{n+1} \right\Vert 
 & \text{ if }   \bar{Z}_{n+1} \neq 0 ,
 \\ 
\hat{Z}_n  & \text{ if }  \bar{Z}_{n+1} = 0 ,
\end{cases}
$
where
\begin{align}
\label{eq:2.4}
 \bar{Z}_{n+1}
& :=
 \hat{Z}_n + \left( 
 \bar{ b } \left( \bar{\eta}_n , \hat{Z}_n \right) 
- \left\langle  \hat{Z}_n , \bar{ b } \left( \bar{\eta}_n , \hat{Z}_n \right) \right\rangle  \hat{Z}_n
+ \Psi  \left( \bar{\eta}_n , \hat{Z}_n \right) 
\right) \Delta
\\
\nonumber
& \quad
+ \sum_{k=1}^m  \left( \bar{ \sigma }^k \left(\bar{\eta}_n  , \hat{Z}_n \right)  
- \left\langle  \hat{Z}_n ,  \bar{ \sigma }^k \left(\bar{\eta}_n  , \hat{Z}_n \right) \right\rangle  \hat{Z}_n \right)
 \sqrt{\Delta} \hat{W}^k_{n+1} ,
\end{align}
$\Psi$ is given by \eqref{eq:2.21}, 
the functions
$\bar{ b }$, $\bar{ \sigma }^k$ are described by (\ref{eq:3.5}),
and 
$$
\mu_n 
=
 \left\langle \hat{Z}_n , \bar{ b } \left( \bar{\eta}_n , \hat{Z}_n \right) \right\rangle 
+\frac{1}{2} \sum_{k=1}^m \left\Vert\  \bar{  \sigma }^k \left( \bar{\eta}_n ,  \hat{Z}_n \right) \right\Vert^2
-  \sum_{k=1}^m \left\langle \hat{Z}_n , \bar{ \sigma }^k \left( \bar{\eta}_n ,  \hat{Z}_n \right) \right\rangle^2 
.
$$ 
Here, $\bar{X}_n := \bar{\eta}_{n} \hat{Z}_{n} $ approximates the solution $ X_{T_n}$
of the SDE \eqref{eq:1.1} with $b\left(0\right)=\sigma^1\left(0\right)=\cdots=\sigma^m\left(0\right)=0$ 
for all $n \in \mathbb{Z}_+$.
\end{scheme}

\begin{remark}
Since $\bar{X}_n $ could be approximately $0$,
we implement Scheme \ref{scheme:EulerStable}
by computing  $\bar{\eta}_{n}$ and $\hat{Z}_{n}$
rather than  $\bar{X}_n $,
which avoids possible round-off errors (see Subsection \ref{subsec:R_Errors}).
\end{remark}

If (\ref{eq:1.1}) reduces to the bilinear SDE
\begin{equation}
\label{eq:BilinealSDE}
X_{t}
=
X_{0} 
+ \int_{0}^{t} B X_{s}  \ ds + \sum_{k=1}^{m} \int_{0}^{t} \sigma^{k} X_{s} \ dW^{k}_{s} 
\end{equation}
with 
$B,  \sigma^{1}, \ldots,  \sigma^{m} \in \mathbb{R}^{d \times d}$,
then Scheme \ref{scheme:EulerStable} becomes
\begin{scheme}
\label{scheme:EulerStableBilinear}
Define recursively 
 $
 \hat{Z}_{n+1} 
 =
 \begin{cases}
 \bar{Z}_{n+1} /  \left\Vert  \bar{Z}_{n+1} \right\Vert 
&
\text{if } \bar{Z}_{n+1} \neq 0 ,
\\
\hat{Z}_{n} 
&
\text{if } 
\bar{Z}_{n+1} = 0 ,
\end{cases}
$
where
\begin{equation}
\label{eq:2.12}
 \bar{Z}_{n+1}
  = 
\hat{Z}_{n}
+  B_n \hat{Z}_n \Delta
+ \sum_{k=1}^m \left(\sigma^k-\left\langle \hat{Z}_n ,\sigma^k \hat{Z}_n \right\rangle \right) \hat{Z}_{n} \sqrt{\Delta}  \hat{W}^{k}_{n+1}
\end{equation}
with $\bar{\eta}_{0}, \hat{Z}_0, \hat{W}^k_n$ as in Scheme \ref{scheme:EulerStable},
and 
$$
B_n
=
B - \left\langle \hat{Z}_n ,B \hat{Z}_n \right\rangle
+
\sum_{k=1}^m 
\left(
\frac{3}{2} \left\langle \hat{Z}_n ,\sigma^k \hat{Z}_n \right\rangle^2 
- \left\langle \hat{Z}_n ,\sigma^k \hat{Z}_n \right\rangle\sigma^k
- \frac{1}{2} \left\Vert\sigma^k \hat{Z}_n \right\Vert^2  
\right) .
$$
The stochastic process $\bar{\eta}_{n+1}$ is given by the iterative formula 
\begin{multline*}
 \bar{\eta}_{n+1} 
 = 
  \bar{\eta}_n \exp  \left( 
 \left( 
 \langle \hat{Z}_{n} , B \hat{Z}_{n} \rangle 
 + \frac{1}{2}\sum_{k=1}^m \left\Vert \sigma^k \hat{Z}_{n}  \right\Vert^2
- \sum_{k=1}^m 
\langle \hat{Z}_{n} , \sigma^k \hat{Z}_{n}  \rangle^2
\right)  \Delta
\right.
\\
\left.
 +
 \sum_{k=1}^m \langle \hat{Z}_{n}, \sigma^k \hat{Z}_{n} \rangle
 \sqrt{\Delta} \hat{W}^k_{n+1} 
 \right) .
\end{multline*}
Thus, $\bar{X}_n := \bar{\eta}_{n} \hat{Z}_n $
approximates the solution $X_{T_n}$ of \eqref{eq:BilinealSDE} for all $n \in \mathbb{Z}_+$.
\end{scheme}

In case $d=1$,
from \eqref{eq:2.4} we have
$\bar{Z}_n = \hat{Z}_0 = \pm 1$ for all $n \geq 0$.
Therefore, Scheme \ref{scheme:EulerStable} reduces to the follow
\begin{scheme}
\label{scheme:EulerStableScalar}
Given $\bar{X}_{0}$,
 \begin{equation*}
\bar{X}_{n+1}
=
 \bar{X}_{n} \exp\left(
\left(
\bar{b} \left(  \bar{X}_n\right)
-
\frac{1}{2} \sum_{k=1}^m   \bar{\sigma}^k \left(\bar{X}_n \right)^2  \right)\Delta
+
\sum_{k=1}^m  \bar{\sigma}^k \left(\bar{X}_n \right)  \sqrt{\Delta} \hat{W}^k_{n+1} \right) 
\end{equation*}
approximates the solution $X_{T_{n+1}}$ of \eqref{eq:1.1} with $d=1$,
where
$
\bar{b} \left( x \right)
= 
\begin{cases}
b \left( x \right) / x 
&
\text{if } x \neq 0 ,
\\
b^{\prime}  \left( 0 \right) 
&
\text{if } 
x = 0 ,
\end{cases}
$
$
\bar{\sigma}^k  \left( x \right)
= 
\begin{cases}
 \sigma^{k} \left( x \right) / x 
&
\text{if } x \neq 0 ,
\\
 \left(\sigma^{k}\right) ^{\prime}  \left( 0 \right)
&
\text{if } 
x = 0 
\end{cases}
$
and
the $\hat{W}^k_n$'s 
are i.i.d.  symmetric  random variables with variance $1$ that are independent of  $\bar{X}_{0}$.
\end{scheme}

\begin{remark}
A projection technique on the sphere was introduced in \cite{Talay1991} to approximate 
the upper Lyapunov exponent of \eqref{eq:1.1} with $b$, $\sigma^{1}, \ldots, \sigma^{m}$ linear.
As a difference with the DND approach considered here,
the method of \cite{Talay1991} does not involve the numerical solution of the coupled system of the SDEs (\ref{eq:2.3}) and (\ref{eq:2.2}),
which describes the evolution of the norm of $X_t$ and the projection of $X_t$  on the unit sphere, respectively.
\end{remark}

\subsection{General SDEs}
\label{subsec:GeneralSDEs}


Suppose that at least one of the vectors 
$b\left(0\right), \sigma^1\left(0\right), \ldots$, $\sigma^m\left(0\right)$ 
is different from $0$.
Consider the constant function $V_t = \alpha$, 
where $\alpha \neq 0$.
From \eqref{eq:1.1} it follows that
\begin{equation}
\label{eq:3.4}
 \begin{pmatrix}
 X_t 
 \\ V_t
\end{pmatrix}  
= 
\begin{pmatrix}
 X_{0} \\ \alpha
\end{pmatrix} 
+
\int_{0}^t 
f \left( \begin{pmatrix} X_s  \\ V_s \end{pmatrix}   \right)
ds
+ \sum_{k=1}^{m}\int_{0}^{t}  
g^k \left( \begin{pmatrix} X_s  \\ V_s \end{pmatrix} \right)
dW^{k}_{s} ,
\end{equation}
with 
$
f \left( 
\begin{pmatrix}
 x
 \\ v
\end{pmatrix}  \right) 
=
\begin{pmatrix}
 b\left(  x \right) - b\left( 0 \right) + \frac{ b \left( 0 \right) }{ \alpha } v 
 \\ 0
\end{pmatrix} 
$
and
$
g^k \left( 
\begin{pmatrix}
 x
 \\ v
\end{pmatrix}  \right) 
=
\begin{pmatrix}
 \sigma^{k} \left( x \right) -  \sigma^{k} \left( 0 \right) +   \frac{ \sigma^k \left( 0 \right)}{\alpha}  v 
 \\ 0
\end{pmatrix} 
$
for all  $x \in \mathbb{R}^d$ and $v \in \mathbb{R}$.
As $f \left( 0 \right) = g^1 \left( 0 \right)= \cdots  g^m \left( 0 \right) = 0$,
we can compute $X_{t}$ by applying Scheme \ref{scheme:EulerStable} to \eqref{eq:3.4}.
Since $V_{T_n} = \alpha$, 
a better alternative is to approximate 
 $X_t$ and $V_t$ in  $\left[ T_n , T_{n+1} \right]$ by the solution of 
\begin{equation}
\label{eq:3.1}
 \begin{pmatrix}
 \widetilde{X}_t 
 \\  \widetilde{V}_t
\end{pmatrix}  
= 
\begin{pmatrix}
 \bar{X}_n \\ \alpha
\end{pmatrix} 
+
\int_{T_n}^t 
f \left( \begin{pmatrix}  \widetilde{X}_s  \\  \widetilde{V}_s \end{pmatrix}   \right)
ds
+ \sum_{k=1}^{m}\int_{T_n}^{t}  
g^k \left( \begin{pmatrix}  \widetilde{X}_s  \\  \widetilde{V}_s \end{pmatrix} \right)
dW^{k}_{s} ,
\end{equation}
with $\bar{X}_n \approx X_{T_n}$.
Then
$ 
\left( \widetilde{X}_{T_{n+1}} , \widetilde{V}_{T_{n+1}} \right)^{\top}
\approx 
\bar{\rho}_{n+1} \bar{Z}_{n+1}
$,
where 
$\bar{\rho}_{n+1}$ and $\bar{Z}_{n+1}$ are given by one iteration of  
Scheme \ref{scheme:EulerStable} applied to  (\ref{eq:3.1}).
Hence,
we compute  $X_{T_{n+1}}$ by projecting $\bar{\rho}_{n+1} \bar{Z}_{n+1}$ 
onto its first $d$ coordinates.
After a short algebraic manipulation we get Scheme \ref{scheme:EulerStableG} 
with $\alpha \neq 0$
(see Remark \ref{rem:ObtencionEulerStableG} for details).

\begin{scheme}
\label{scheme:EulerStableG}
Let $\bar{X}_0$ be a random variable with values in $\mathbb{R}^d$
such that $ \bar{X}_{0} \neq 0$
in case $ b\left(0\right) = \sigma^1\left(0\right) = \cdots = \sigma^m\left(0\right) = 0 $.
Suppose that $\hat{W}^1_1, \hat{W}^2_1, \ldots, \hat{W}^m_1, \hat{W}^1_2, \ldots$ 
are i.i.d. symmetric random variables with variance $1$ that are independent of  $\bar{X}_{0}$.
Choose $\alpha \in \mathbb{R}$ such that 
$ \alpha = 0 $  if $ b\left(0\right) = \sigma^1\left(0\right) = \cdots = \sigma^m\left(0\right) = 0 $, 
and $ \alpha \neq 0 $ otherwise.
Then,  for all $n \in \mathbb{Z}_+$ the solution $X_{T_{n+1}}$ of \eqref{eq:1.1}
is recursively approximated by
\begin{equation}
 \label{eq:3.3}
 \bar{X}_{n+1} 
 = 
\begin{cases}
\bar{\rho}_{n+1} \bar{U}_{n+1} / \sqrt{ \left\Vert \bar{U}_{n+1} \right\Vert^2 +  \left( \bar{V}_{n+1} \right)^2}
&
\text{if } \left( \bar{U}_{n+1}, \bar{V}_{n+1}  \right) \neq 0 ,
\\
\bar{\rho}_{n+1} \hat{U}_{n} 
&
\text{if } 
\left( \bar{U}_{n+1}, \bar{V}_{n+1}  \right) = 0 ,
\end{cases} 
\end{equation}
where 
$
\bar{\rho}_{n+1} = \bar{\eta}_n  \exp  \left( 
\mu_n \Delta + 
\sum_{k=1}^m \left\langle \hat{U}_n  ,\bar{ \sigma }^k \left(\bar{\eta}_n  , \hat{U}_n \right) \right\rangle \sqrt{\Delta} \hat{W}^k_{n+1}
\right) 
$,
\begin{equation}
\label{eq:7.13}
 \begin{aligned}
 \bar{U}_{n+1}
& =
\hat{U}_n 
+ \left( 
 \bar{ b } \left( \bar{\eta}_n , \hat{U}_n \right) 
  - \left\langle  \hat{U}_n , \bar{ b } \left( \bar{\eta}_n , \hat{U}_n \right) \right\rangle  \hat{U}_n
  + \Psi  \left( \bar{\eta}_n , \hat{U}_n \right)
 \right) \Delta
 \\
 & \quad
 +   \sum_{k=1}^m  \left( \bar{ \sigma}^k \left(\bar{\eta}_n  , \hat{U}_n \right)  
- \left\langle  \hat{U}_n ,  \bar{ \sigma}^k \left(\bar{\eta}_n  , \hat{U}_n \right) \right\rangle  \hat{U}_n \right)
\sqrt{\Delta} \hat{W}^k_{n+1} ,
\end{aligned}
\end{equation}
and
\begin{align*}
 \bar{V}_{n+1}
 & =
\frac{ \alpha }{ \bar{\eta}_n }
 -  \frac{ \alpha }{ \bar{\eta}_n } \left\langle  \hat{U}_n , \bar{ b } \left( \bar{\eta}_n , \hat{U}_n \right) \right\rangle  \Delta
 \\
 & \quad +
\frac{ \alpha }{ \bar{\eta}_n }  \Delta \sum_{k=1}^m   \left( \frac{3}{2} \left\langle  \hat{U}_n ,\bar{  \sigma }^k \left(\bar{\eta}_n  , \hat{U}_n \right) \right\rangle^2
- \frac{1}{2} \left\Vert \bar{ \sigma }^k \left(\bar{\eta}_n  , \hat{U}_n \right) \right\Vert ^2 
 \right)  
 \\
 & \quad
 - \frac{ \alpha }{ \bar{\eta}_n } \sum_{k=1}^m   
  \left\langle  \hat{U}_n ,  \bar{ \sigma}^k \left(\bar{\eta}_n  , \hat{U}_n \right) \right\rangle   
 \sqrt{\Delta} \hat{W}^k_{n+1} .
\end{align*}
Here, $\bar{b}, \bar{ \sigma }^k$ are given by \eqref{eq:3.5}, 
$\Psi$ is defined by \eqref{eq:2.21}, 
$
\bar{\eta}_{n}  = \sqrt{\left\Vert  \bar{X}_n  \right\Vert^2 + \alpha^2}
$,
$
\hat{U}_{n} 
:=
\bar{X}_{n} /  \bar{\eta}_{n} 
$
and 
$$
\mu_n =
 \left\langle \hat{U}_n , \bar{ b } \left( \bar{\eta}_n , \hat{U}_n \right) \right\rangle 
+\frac{1}{2} \sum_{k=1}^m \left\Vert\  \bar{  \sigma }^k \left( \bar{\eta}_n ,  \hat{U}_n \right) \right\Vert^2
-  \sum_{k=1}^m \left\langle \hat{U}_n , \bar{ \sigma }^k \left( \bar{\eta}_n ,  \hat{U}_n \right) \right\rangle^2 .
$$
\end{scheme}

\begin{remark}
\label{rem:ObtencionEulerStableG}
Consider Scheme \ref{scheme:EulerStableG} with $ \alpha \neq 0 $.
Then
$\left( \bar{U}_{n+1}, \bar{V}_{n+1}  \right) = \bar{Z}_{n+1}$,
where $\bar{Z}_{n+1}$ is given by (\ref{eq:2.4}) with 
$b$, $\sigma^k$, $\bar{\eta}_{n}$ and $\hat{Z}_n$ replaced by 
$f$, $g^k$, $\sqrt{\left\Vert  \bar{X}_n  \right\Vert^2 + \alpha^2}$  and 
$ \left(  \bar{X}_n , \alpha \right) / \sqrt{\left\Vert  \bar{X}_n  \right\Vert^2 + \alpha^2} $, respectively.
Similarly, 
$\bar{\rho}_{n+1}$ arises from evaluating (\ref{eq:2.5}) in the context of  \eqref{eq:3.1}.
Therefore, 
one iteration of  Scheme \ref{scheme:EulerStable} applied to  (\ref{eq:3.1})
maps 
$
\left( 
\sqrt{\left\Vert  \bar{X}_n  \right\Vert^2 + \alpha^2}, 
\frac{
\left(  \bar{X}_n , \alpha \right) 
}{ 
\sqrt{\left\Vert  \bar{X}_n  \right\Vert^2 + \alpha^2 }  } 
\right)
$
to 
$
\left(
 \bar{\rho}_{n+1} , 
\frac{
\left( \bar{U}_{n+1} , \bar{V}_{n+1} \right) }{
 \sqrt{ \left\Vert \bar{U}_{n+1} \right\Vert^2 +  \left( \bar{V}_{n+1} \right)^2}
 } 
\right) 
$
if $\left( \bar{U}_{n+1}, \bar{V}_{n+1}  \right) \neq 0$,
a computation which is not affected by round-off errors
as we discuss in Subsection \ref{subsec:R_Errors}.
Hence,
$\bar{X}_{n+1}$ coincides with the first $d$ coordinates of 
the $\mathbb{R}^{d+1}$-vector obtained by applying 
one-step of Scheme \ref{scheme:EulerStable} to  (\ref{eq:3.1}).
\end{remark}

\begin{remark}
In this paper, using heuristic arguments we select 
$$
 \alpha = \max \left\{ \left\Vert b \left( 0 \right) \right\Vert_{\infty} ,    \left\Vert \sigma^1 \left( 0 \right)   \right\Vert_{\infty} , 
\ldots , \left\Vert  \sigma^m \left( 0 \right)  \right\Vert_{\infty} \right\} /2 .
$$
We will address
the problem of finding optimal values for $\alpha$
in a further work.
\end{remark}

\begin{remark}
Let $ b\left(0\right) = \sigma^1\left(0\right) = \cdots = \sigma^m\left(0\right) = 0 $.
Then,
Scheme \ref{scheme:EulerStableG} coincides with Scheme \ref{scheme:EulerStable} with
$\hat{Z}_0 =  \bar{X}_{0}  /  \left\Vert  \bar{X}_{0}  \right\Vert$ and $\bar{\eta}_{0} =  \left\Vert  \bar{X}_{0}  \right\Vert$.
In this case,
the rewriting of $\bar{X}_{n+1}$ as the pair  $\left( \hspace{-1.5pt} \bar{\eta}_{n+1} ,  \hat{Z}_{n+1} \hspace{-1.5pt} \right)$,
given by Scheme \ref{scheme:EulerStable}, 
is a key factor in the implementation of Scheme \ref{scheme:EulerStableG}.
\end{remark}


\subsection{Round-off errors}
\label{subsec:R_Errors}

The computer implementation of Scheme \ref{scheme:EulerStable} 
is seldom influenced by the effect of round-off errors.
First, Scheme \ref{scheme:EulerStable} involves the calculation of 
$\bar{ \sigma }^{k} \left( \eta, z \right)$
for  $\left\Vert z \right\Vert = 1$ and $\eta \geq 0$,
where $k=0,\ldots,m$ and $\bar{ \sigma }^0 := \bar{ b }$.
Sometimes we can evaluate efficiently closed analytical expressions for 
$\bar{ \sigma }^{k} $ like in Subsections \ref{sec:NumExpBilinear}
and  \ref{sec:NumExpScalar}.
In general,
using the smoothness of  $\bar{ \sigma }^{k}$
we avoid the effect of
round-off errors in the implementation of 
$  \bar{ \sigma }^k \left( \eta, z \right)$
when $\eta$ is near $0$.
Indeed, if
$\sigma^k \in \mathcal{C}^{ \ell +1 } \left( \mathbb{R}^d,\mathbb{R}^d\right)  $,
then applying Lemma \ref{lem:Derivadas} we deduce that 
$  \eta \mapsto \bar{ \sigma }^{k,j} \left( \eta, z \right)$
is $\ell$-times  continuously differentiable,
and 
$$
\frac{d^ \ell}{d \eta ^ \ell} \bar{ \sigma }^{k,j} \left( 0, z \right)
= \frac{1}{\ell+1} \left(
\sum_{j_1, \ldots, j_{\ell+1} =1}^d 
\frac{\partial^{\ell+1} \sigma^{k,j}}{\partial x^{j_1} \cdots \partial x^{j_{\ell+1} }} 
 \left( 0\right) z^{j_1} \cdots  z^{j_{\ell+1}} 
\right) ,
$$
where 
$ \bar{ \sigma }^{k,j} \left( \eta, z \right) $, $\sigma^{k,j} \left( \eta, z \right) $ 
denote the $j$-th coordinate of $\bar{ \sigma }^{k} \left( \eta, z \right)$, $\sigma^{k} \left( \eta, z \right) $, respectively.
This allows us to approximate successfully 
$\bar{ \sigma }^{k,j} \left( \eta, z \right)$, with $\eta \approx 0$,
by means of the truncated Taylor expansions of 
$  \eta \mapsto \bar{ \sigma }^{k,j} \left( \eta, z \right)$
around $0$.
Alternatively, 
we can interpolate 
$ \bar{ \sigma }^{k,j} \left( \cdot, z \right)$ in a neighborhood of $0$.

\begin{lemma}
 \label{lem:Derivadas}
Consider
$
g \left( x \right)
= 
\begin{cases}
\left( f \left( x \right) - f \left( 0 \right) \right) / x
&
\text{if } x \neq 0 ,
\\
f^{\prime}  \left( 0 \right) 
&
\text{if } 
x = 0 ,
\end{cases}
$
where $x \in  \mathbb{R}$ and $f \in \mathcal{C}^{ \ell +1 } \left( \mathbb{R},\mathbb{R}\right)  $
with $\ell \in \mathbb{N}$.
Then $g$ is $\ell$-times continuously differentiable and
$$
\frac{d^ \ell}{dx^\ell} 
g  \left( 0 \right) = \frac{1}{\ell+1} \frac{d^{\ell+1}}{dx^{\ell+1}} f \left( 0 \right) .
$$ 
\end{lemma}

\begin{proof}
 Combining Leibniz's rule with Taylor's theorem we obtain the assertion of the lemma.
\end{proof}

Second,
we project $\bar{Z}_{n+1}$ on the unit sphere.
Since
$\left\Vert \bar{Z}_{n+1} \right\Vert \approx \left\Vert Z_{T_{n+1}} \right\Vert = 1$,
$\left\Vert \bar{Z}_{n+1} \right\Vert $ usually keeps away from $0$.
In fact,
$\left\Vert  \bar{Z}_{n+1} \right\Vert$ may take  small values 
only for certain special combinations  of $\bar{\eta}_n$, $\hat{Z}_n$ and $\Delta$.
For example, 
$\bar{Z}_{n+1} \approx 0$
implies 
$
  \left\langle \hat{Z}_n , \bar{Z}_{n+1}  \right\rangle \approx 0
$,
and hence 
$$
\sum_{k=1}^m 
\left(
 \left\Vert  \bar{ \sigma}^k \left(\bar{\eta}_n  , \hat{Z}_n \right)  \right\Vert^2
 -
 \left\langle  \hat{Z}_n   ,  \bar{ \sigma}^k \left(\bar{\eta}_n  , \hat{Z}_n \right)  \right\rangle^2 
\right)
\approx
2 / \Delta .
$$ 
Let $\hat{W}_n^{k}$ be distributed uniformly on $\left[ - \sqrt{3}, \sqrt{3} \right]$. 
According to Lemma \ref{lem:ZnoNulo} below we have that
$\bar{Z}_{n+1} \neq 0$ a.s. for all $\Delta > 0$,
and the proof of Lemma \ref{lem:ZnoNulo} suggests us that
the special cases where 
$\left\Vert  \bar{Z}_{n+1} \right\Vert \approx 0$
could happen only with quite small probability.
In the latter situation,
we can use a preconditioner like
$$
 \hat{Z}_{n+1} 
 = 
 \left( \bar{Z}_{n+1} /  \left\Vert  \bar{Z}_{n+1} \right\Vert_\infty \right)
 /
 \left\Vert \bar{Z}_{n+1} /  \left\Vert  \bar{Z}_{n+1} \right\Vert_\infty \right\Vert ,
$$
or set $\hat{Z}_{n+1}=\hat{Z}_{n}$ in the worst case.
Furthermore,
if $b$ and $\sigma^k$  have at most linear growth,
then 
combining $\left\Vert  \hat{Z}_{n} \right\Vert =1$ with 
the compact support property of $\hat{W}_n^{k}$
we deduce that
$\left\Vert  \bar{Z}_{n+1} \right\Vert $ is uniformly bounded from below by a positive constant
whenever $\Delta$ is small enough.

\begin{lemma}
\label{lem:ZnoNulo}
 Adopt the framework of Scheme \ref{scheme:EulerStable}.
 Let the distribution of $\hat{W}^{k}_{n}$ be absolutely continuous with respect to the Lebesgue measure.
 Then, for all $n \geq 0$, $\bar{Z}_{n} \neq 0$ a.s.
\end{lemma}

\section{Convergence properties}
\label{sec:Convergence}

\subsection{Long time behavior}
\label{subsec:LongTime}

In our numerical experiments, 
Scheme \ref{scheme:EulerStable} reproduces very well the long time behavior of $ X_{t} $.
Next,  we assert that  for all $\Delta > 0$,
Scheme \ref{scheme:EulerStable}  converges almost sure exponentially fast to $0$
under  a classical condition for 
$
\limsup_{t \rightarrow \infty}  \left( \log \left\Vert X_t \right\Vert \right) / t
< 0 
$ 
(see, e.g., \cite{Higham2007,Mao2007}).

\begin{hypothesis}
\label{hyp:Basica}
The functions  
$b, \sigma^1, \ldots, \sigma^m$
have continuous first-order partial derivatives.
Moreover,
$ \bar{\eta}_{0} >0$,
$b \left( 0 \right) = 0$
and  
$
 \left\Vert  \sigma^{k} \left(x\right) \right\Vert   \leq  K  \left\Vert x \right\Vert
$
for all $x \in \mathbb{R}^d$.
\end{hypothesis}

\begin{theorem}
\label{th:Stable}
Assume Hypothesis \ref{hyp:Basica},
together with the condition \eqref{eq:StabCond}.
Let $ \bar{X}_{n}$ be given by Scheme \ref{scheme:EulerStable}.
Then  
\begin{equation}
\label{eq:7.31}
 \limsup_{n  \rightarrow\infty}\frac{1}{n \Delta }\log \left( \left\Vert \bar{X}_{n} \right\Vert  \right) \leq -\lambda 
\hspace{1cm}
a.s.
\end{equation}
\end{theorem}

Under the assumptions of Theorem \ref{th:Stable},
the $p$-th moment of $\left\Vert X_t \right\Vert$ converges exponentially to $0$ as $t \rightarrow + \infty$
provided that $p > 0$ is sufficiently small (see, e.g., \cite{Higham2007}).
Theorem \ref{th:MomentStable}  shows that  Scheme \ref{scheme:EulerStable} preserves this behavior
for any step-size $\Delta > 0$.

\begin{theorem}
\label{th:MomentStable}
In addition to the assumptions of Theorem \ref{th:Stable},
we suppose that there exists $\gamma \geq 0 $ satisfying 
\begin{equation}
\label{eq:7.2}
\mathbb{E} \exp \left(  s \, \hat{W}^1_1 \right)
\leq 
\exp \left( \gamma \, s^2  \right)
\qquad \quad
\forall s \in \mathbb{R} .
\end{equation}
Then, 
for any $\epsilon \in \left] 0 , \lambda \right[$ there exists $q \in \left] 0 , 1\right[$ such that for all $p \in \left] 0 , q \right[$ 
and $\Delta > 0$,
\begin{equation}
\label{eq:5.3}
\mathbb{E} \left(  \left\Vert \bar{X}_{n} \right\Vert ^p \right)
\leq 
\exp \left( - \left( \lambda - \epsilon \right)  p \, n \, \Delta \right)  \mathbb{E} \left(  \left\Vert \bar{X}_{0} \right\Vert ^p \right)
\qquad \quad
\forall n \in \mathbb{N}. 
\end{equation}
\end{theorem}

\begin{remark}
If 
$\mathbb{P} \left( \hat{W}^{1}_{1} \in \left[ a , b \right] \right) = 1$,
then according to Hoeffding's lemma we have that
$
\mathbb{E} \exp \left(  s \, \hat{W}^1_1 \right)
\leq
\exp \left( s^2  \left( b - a \right)^2 / 8 \right) 
$
for any $s \in \mathbb{R}$.
On the other hand,
$
\mathbb{E} \exp \left(  s \, \hat{W}^1_1 \right)
=
\exp \left( s^2  / 2 \right)
$
whenever
$ \hat{W}^1_1 $ 
is normally distributed.
\end{remark}

Now, we establish that 
 the Scheme \ref{scheme:EulerStable} approximations are 
away from $0$ for any step-size 
in a case where $0$ is an unstable equilibrium point of \eqref{eq:1.1}
(see, e.g., \cite{Appleby2008}).

\begin{theorem}
\label{th:NoStable}
Let Hypothesis \ref{hyp:Basica} hold.
Assume  the existence of $\theta > 0$ such that 
\begin{equation}
\label{eq:2.7}
\frac{
\langle x, b \left( x \right)  \rangle
+
\frac{1}{2} \sum_{k=1}^{m}  \left\Vert \sigma^k \left( x \right)  \right\Vert^2 
 }{
 \left\Vert x \right\Vert^2 
 }
- \left( 1 + \theta \right) 
\frac{ 
\sum_{k=1}^{m}
\left\langle x,  \sigma^k \left( x \right)  \right\rangle^2
}{
\left\Vert x \right\Vert^4
}
\geq 
0
\end{equation}
for all $x \neq 0$.
Then 
$
\liminf_{n \rightarrow \infty}  \left\Vert \bar{X}_{n} \right\Vert  >  0
$
a.s.,
where $\bar{X}_{n}$ is given by Scheme \ref{scheme:EulerStable}.
\end{theorem}

\subsection{Rates of weak and strong convergence}
\label{subsec:RateConvergence}

We will estimate the errors arising from the computation of 
$
\mathbb{E} \varphi \left( X_t \right)
$
with $\varphi \in \mathcal{C}_{P}^{4}\left( \mathbb{R}^{d},\mathbb{R}\right) $.
To this end,
using classical arguments (see, e.g., \cite{GrahamTalay2013,Kloeden1992,Milstein2004}) 
we obtain the following general criterion for  ensuring the linear weak convergence of  
the numerical methods. 

\begin{hypothesis}
\label{hyp:EDP}
Fix $\varphi \in \mathcal{C}_{P}^{4}\left( \mathbb{R}^{d},\mathbb{R}\right) $ and $T>0$.
Let  $b, \sigma^1, \ldots, \sigma^m \in \mathcal{C}_{P}^{4}\left( \mathbb{R}^{d},\mathbb{R}\right) $.
Suppose that there exists 
$u \in \mathcal{C}^{1,4} \left(\left[ 0, T \right] \times \mathbb{R}^{d},\mathbb{R}\right) $
such that
\begin{equation}
\label{eq:Kolmogorov}
\begin{cases}
\frac{\partial }{\partial t }u \left( t ,x\right) 
 =
-\mathcal{L}\left( u \right) \left( t ,x\right)  
&
\text{if } t \in \left[ 0,T \right] \text{ and } x\in \mathbb{R}^{d},
\\
u \left( T,x\right) = \varphi \left( x\right)  
& \text{ if } x \in \mathbb{R}^{d} ,
\end{cases}
\end{equation}
where 
$
\mathcal{L}=\sum_{k=1}^{d}b_{k}\partial _{x}^{k}+\frac{1}{2}%
\sum_{k, \ell =1}^{d}
\left( \sum_{j =1}^{m} \sigma_k^j  \sigma_{\ell}^j 
\right)
\partial_{x}^{k, \ell}
$.
Moreover, we assume that 
$
u \in  \mathcal{C}_{P}^{ 4 } \left( \left[ 0, T \right] \times \mathbb{R}^d,\mathbb{R}\right)
$
and 
\begin{equation}
 \label{eq:KolmogorovG}
 \frac{\partial }{\partial t}\partial _{x}^{\vec{p}}u
=
-\partial_{x}^{\vec{p}}\mathcal{L}\left( u \right) 
\hspace{1cm}
\forall \vec{p}\in \mathcal{P}_{1} \bigcup \mathcal{P}_{2}.
\end{equation}
\end{hypothesis}

\begin{theorem}
 \label{th:GeneralWeakConvergence}
 Let Hypothesis \ref{hyp:EDP} hold,
 and let  $\mathbb{E} \left\Vert X_{0} \right\Vert ^{q} < + \infty$ for all $q \in \mathbb{N}$.
 Consider the sequence of the random variables 
 $\left( \bar{Y}^N_0, \ldots, \bar{Y}^N_{\ell \left( N \right) } \right)_{N \in \mathbb{N}}$,
 real numbers
$ 0 = \tau^N_0 < \tau^N_1 < \cdots < \tau^N_{\ell \left( N \right) } = T$ 
 and filtrations  
 $ \left( \mathfrak{G}^N_{0}, \ldots, \mathfrak{G}^N_{\ell \left( N \right)} \right)_{N \in \mathbb{N}}$
 such that
 $\tau^N_{n+1} - \tau^N_n \leq T/N$
 and 
 $W^k_{\tau^N_{n+1}} - W^k_{\tau^N_n}$ is independent of $\mathfrak{G}^N_{n}$
 for all $n=0,\ldots, \ell \left( N \right) -1$ and $N \in \mathbb{N}$,
 as well as 
 $\bar{Y}^N_n$, $W^k_{\tau^N_n}$ are $\mathfrak{G}^N_{n}$-measurable
 for any  $n=0,\ldots, \ell \left( N \right)$ and $N \in \mathbb{N}$. 
 Assume: 
\begin{description}

\item[(i)] 
For every $ \phi \in C_{p}^{4} \left(  \mathbb{R}^{d},\mathbb{R}\right)  $
there exist $q \in \mathbb{N}$ and $K \geq 0$ such that for all $N \in \mathbb{N}$,
$
 \left\vert \mathbb{E} \phi \left(X_{0}\right) -  \mathbb{E} \phi \left(  \bar{Y}^N_{0}  \right)   \right\vert 
 \leq K\left(  1+\mathbb{E}\left\Vert X_{0} \right\Vert ^{q}\right) T/N  
$.

\item[(ii)] 
For any $p \in \mathbb{N}$ there exist $q \in \mathbb{N}$ and $K \geq 0$ such that 
for all $ n=0, \ldots \ell \left( N \right)$ and $N \in \mathbb{N}$,
$
\mathbb{E} \left(  \left\Vert  \bar{Y}^N_{n} \right\Vert^{p} \right)
 \leq
 K  \left( 1 + \mathbb{E} \left(  \left\Vert  \bar{Y}^N_{0} \right\Vert^{q } \right) \right) 
$.

\item[(iii)]  
For any $p \in \mathbb{N}$ there exist $q \in \mathbb{N}$ and $K \geq 0$ such that 
for all $N \in \mathbb{N}$ and $ n=0, \ldots \ell \left( N \right)-1$,
$
\mathbb{E} \left(  \left\Vert  \bar{Y}^N_{n+1}  - \bar{Y}^N_{n} \right\Vert^{2p} \right)
 \leq
 K   \left( 1 + \mathbb{E} \left(  \left\Vert  \bar{Y}^N_{n} \right\Vert^{q } \right) \right) 
 \left(\tau^N_{n+1} - \tau^N_n \right)^p
$.

\item[(iv)]   
For all $\ell=1,2,3$ and $\vec{p}\in P_{\ell}$, 
there exist $q \in \mathbb{N}$ and $K \geq 0$ satisfying
\begin{equation*}
 \begin{aligned}
&
\left\vert
\mathbb{E} \left(  
F_{\vec{p}}\left(  \bar{Y}^N_{n+1} -  \bar{Y}^N_{n} \right) 
-
F_{\vec{p}}\left( E^n_{\tau^N_{n+1}}  -  \bar{Y}^N_{n} \right) 
 \diagup  \mathfrak{G}^N_{n} \right)
\right\vert
\\
& \leq
K \left( 1+\left\Vert \bar{Y}^N_n \right\Vert ^{q}\right) 
\left(\tau^N_{n+1} - \tau^N_n \right) T/N
\end{aligned}
\hspace{1cm} \forall n=0, \ldots \ell \left( N \right) -1
\end{equation*}
for all natural number $N$.
Here, 
$$
E^n_t : = \bar{Y}^N_{n}
+ b\left( \bar{Y}^N_{n} \right) \left( t - \tau^N_{n} \right)
+ \sum_{k=1}^{m} \sigma^{k} \left(\bar{Y}^N_{n} \right) \left(W^{k}_{t} - W^{k}_{\tau^N_{n}} \right) 
\quad \forall t \in \left[ \tau^N_{n} , \tau^N_{n+1} \right] .
$$
\end{description}
Then, 
for any $ \varphi \in C_{p}^{4} \left(  \mathbb{R}^{d},\mathbb{R}\right)  $
there exist $q \in \mathbb{N}$ and $K \geq 0$ such that  
\begin{equation}
\label{eq:7.16}
\left\vert \mathbb{E} \varphi \left(X_{T}\right) -  \mathbb{E} \varphi \left(  \bar{Y}^N_{\ell \left( N \right)}  \right)   \right\vert 
\leq K \left(  1+\mathbb{E}\left\Vert X_{0} \right\Vert ^{q}\right)  T/N
\hspace{1cm} \forall N \in \mathbb{N}. 
\end{equation}
\end{theorem}

\begin{remark}
 Let $u \in \mathcal{C}^{1,4} \left(\left[ 0, T \right] \times \mathbb{R}^{d},\mathbb{R}\right) $
 satisfy \eqref{eq:Kolmogorov}.
 If 
 $\frac{\partial }{\partial t}\partial _{x}^{\vec{p}}u$
 and 
  $\partial _{x}^{\vec{p}} \frac{\partial }{\partial t}u$
 are continuous function for all $\vec{p}\in \mathcal{P}_{1} \bigcup \mathcal{P}_{2}$,
 then
 $
 \frac{\partial }{\partial t}\partial _{x}^{\vec{p}}u
=
\partial _{x}^{\vec{p}} \frac{\partial }{\partial t}u
=
-\partial_{x}^{\vec{p}}\mathcal{L}\left( u \right)
 $
provided that $\vec{p}\in \mathcal{P}_{1} \bigcup \mathcal{P}_{2}$,
and so \eqref{eq:KolmogorovG} holds. 
\end{remark}

\begin{remark}
Hypothesis \ref{hyp:EDP} holds in case 
$b, \sigma^1, \ldots, \sigma^m  \in \mathcal{C}_{P}^{4}\left( \mathbb{R}^{d},\mathbb{R}\right)$
have  uniformly bounded derivatives (see, e.g., \cite{GrahamTalay2013,Kloeden1992,Krylov1999,Milstein2004}).
A class of SDEs with non-globally Lipschitz coefficients that satisfy Hypothesis \ref{hyp:EDP}
is studied in \cite{Talay2002}.
In general,
the solution of the Kolmogorov equation \eqref{eq:Kolmogorov} 
can lose the regularity of the data $\varphi$
(see \cite{Hairer2015}).
\end{remark}

Combining Theorem \ref{th:GeneralWeakConvergence}
with the study of the local asymptotic behavior of $\bar{X}_n$,
we next deduce that Schemes \ref{scheme:EulerStable}  and  \ref{scheme:EulerStableG}
converge with weak rate $1$ as $\Delta \rightarrow 0+$.

\begin{hypothesis}
\label{hyp:Acotacion}
Let $b, \sigma^1, \ldots, \sigma^m$
have continuous first-order partial derivatives.
Suppose that for any $\beta \geq 1$,
\begin{equation}
 \label{eq:7.1}
 \langle x , b \left( x \right) \rangle
 +
\beta \sum_{k=1}^{m}  \left\Vert \sigma^k \left( x \right)  \right\Vert^2
\leq 
K_{\beta} \left( 1 + \left\Vert x \right\Vert^2 \right)
\qquad \quad
\forall x \in \mathbb{R}^d ,
\end{equation}
with $K_{\beta} \geq 0$.
\end{hypothesis}

\begin{theorem}
\label{th:RateConvergence}
Let Hypotheses \ref{hyp:EDP} and  \ref{hyp:Acotacion}  hold,
together with \eqref{eq:7.2} and 
$\mathbb{E} \left\Vert X_{0} \right\Vert ^{q} < + \infty$ for all $q \in \mathbb{N}$.
Consider Scheme \ref{scheme:EulerStableG} with 
$\Delta = T/N$, where $N \in \mathbb{N}$.
Suppose that 
the distribution of $\hat{W}^{k}_{n}$ 
is absolutely continuous with respect to the Lebesgue measure,
and that 
for every $ \phi \in C_{p}^{4} \left(  \mathbb{R}^{d},\mathbb{R}\right)  $
there exist  $q \in \mathbb{N}$ and $K \geq 0$ such that for all $N \in \mathbb{N}$,
$
\left\vert \mathbb{E} \phi \left(X_{0}\right) -  \mathbb{E} \phi \left(  \bar{X}_{0}  \right)   \right\vert 
 \leq K\left(  1+\mathbb{E}\left\Vert X_{0} \right\Vert ^{q}\right)  T/N
$.
Then
\begin{equation}
\label{eq:7.4}
\left\vert \mathbb{E} \varphi \left(X_{T}\right) -  \mathbb{E} \varphi \left(  \bar{X}_{N}  \right)   \right\vert 
\leq K \left(  1+\mathbb{E}\left\Vert X_{0} \right\Vert ^{q}\right)  T/N 
\hspace{1cm} \forall N \in \mathbb{N}.
\end{equation}
\end{theorem}


We now obtain that the new schemes converge in $L^2 \left( \mathbb{P} \right)$ with rate $1/2$
under the assumptions of the fundamental mean-square convergence theorem proved in \cite{TretyakovZhang2013},
where $b$, $\sigma^1, \ldots, \sigma^m $ grow polynomially at infinity and satisfy a one-sided Lipschitz condition.

\begin{theorem}
\label{th:StrongRateConvergence}
Suppose that for any $\beta \geq 1$ there exists $K_{\beta} \geq 0$ such that
$$
\left\langle x-y,b(x)-b(y)\right\rangle 
+
\beta 
\sum_{k=1}^{m}\left\Vert \sigma^{k}(x)-\sigma^{k}(y)\right\Vert ^{2}
\leq K_{\beta} \left\Vert x-y\right\Vert ^{2}
\hspace{0.5cm}
\forall x,y\in\mathbb{R}^{d} ,
$$
and that
$
\left\Vert b(x)-b(y)\right\Vert ^{2}
\leq 
K \left(1+\left\Vert x\right\Vert ^{q}+\left\Vert y\right\Vert ^{q} \right) \left\Vert x-y\right\Vert ^{2}
$
for all $x,y\in\mathbb{R}^{d}$.
Consider 
Scheme  \ref{scheme:EulerStableG} with 
$ \hat{W}^k_{n+1} = \left( W^k_{T_{n+1}} -  W^k_{T_{n}} \right) / \sqrt{\Delta}$
and
$\Delta = T/N$,
where $T > 0$  and $N \in \mathbb{N}$.
Assume that for any $q\in\mathbb{N}$,
$\mathbb{E}\left\Vert X_{0}\right\Vert ^{q}<+\infty$ 
and 
$
\left(
 \mathbb{E}\left(\left\Vert X_{0}-\bar{X}_{0}\right\Vert ^{q}\right)
\right)^{1/q}
\leq
K_q  \, T / N 
$,
with $K_q \geq 0$ independent of $N$.
Then, for any $p \in \mathbb{N}$ we have
\begin{equation}
\label{eq:SEq3}
 \mathbb{E}\left(\left\Vert X_{T_{n}}-\bar{X}_{n}\right\Vert ^{2p}\right)
\leq
K_p \left(1+\mathbb{E}\left\Vert X_{0}\right\Vert ^{q_p}\right)  \left( T / N \right)^{p}
\hspace{1cm}
\forall n = 0, \ldots ,N,
\end{equation}
where $K_p \geq 0$ and $q_p \in \mathbb{N}$ are independent of $N$.
\end{theorem}

\section{Numerical Experiments}
\label{Sec:Numerical Experiments}

In what follows, 
the potential of Scheme \ref{scheme:EulerStable} is illustrated by means of numerical experiments involving the integration of linear and nonlinear equations with a variety of asymptotic behavior. 

\subsection{Bilinear SDE}
\label{sec:NumExpBilinear}
We start with  the test equation 
\begin{equation}
\label{eq:1.16}
X_{t}
=
X_{0}
+\int_{0}^{t} \begin{pmatrix} b & 0 \\ 0 & b \end{pmatrix} X_{s} \, ds
+\int_{0}^{t} \begin{pmatrix} \sigma & 0 \\ 0 & \sigma \end{pmatrix} X_{s} \, dW^{1}_{s}
+\int_{0}^{t} \begin{pmatrix} 0 & -\epsilon \\ \epsilon & 0 \end{pmatrix}X_{s} \, dW^{2}_{s} ,
\end{equation}
where $X_{t} = \left( X^1_{t} , X^2_{t} \right) ^{\top} \in \mathbb{R}^2$
and $b, \sigma, \epsilon \in \mathbb{R}$ (see, e.g., \cite{Berkolaiko2012,Buckwar2010}).
We solve \eqref{eq:1.16} by Scheme \ref{scheme:EulerStableBilinear},
the backward Euler method \eqref{scheme:BackwardEuler} 
and the balanced scheme 
\begin{equation}
 \begin{aligned}
\label{scheme:Balanced}
 \bar{B}_{n+1}
& =
\bar{B}_n+b\left(\bar{B}_n\right)\Delta+\sum_{k=1}^m\sigma^{k}\left(\bar{B}_n\right)\sqrt{\Delta} \hat{W}^k_{n+1}
-\frac{1}{2} Jb \left(\bar{B}_{n}\right)\Delta \left(\bar{B}_n-\bar{B}_{n+1}\right)
\\
& \quad
+\left(\sum_{k=1}^m \sqrt{ J \sigma^k\left(\bar{B}_n\right)^{\top} J \sigma^k\left(\bar{B}_n\right)}
 \left\vert \hat{W}^k_{n+1}
\right\vert \sqrt{\Delta}\right)\left(\bar{B}_n-\bar{B}_{n+1}\right) 
\end{aligned}
\end{equation}
(see, e.g., \cite{Alcock2006,Schurz2005}).
Here, 
$ \sqrt{\Delta} \hat{W}^{k}_{n+1} = W^k_{T_{n+1}} - W^k_{T_{n}} $.
Moreover,
we apply to \eqref{eq:1.16} the S-ROCK scheme \cite{Abdulle2008,Abdulle2008b,Abdulle2013}
\begin{equation*}
\begin{cases}
  \bar{S}_{n+1}
   =
 K_s+\sum_{k=1}^m\sigma^k\left(K_s\right)\sqrt{\Delta}  \hat{W}_{n+1}^k 
 &
 \\
 K_1=\bar{S}_n+\Delta \frac{y_1}{y_0}b\left(K_0\right)
&
 \\
 K_j  = 
2\Delta y_1 \frac{T_{j-1}\left(y_0\right)}{T_{j}\left(y_0\right)}b\left(K_{j-1}\right)+2y_0\frac{T_{j-1}\left(y_0\right)}{T_{j}\left(y_0\right)}K_{j-1}-\frac{T_{j-2}\left(y_0\right)}{T_{j}\left(y_0\right)}K_{j-2}
&
\forall j=2,\ldots,s
\end{cases} ,
\end{equation*}
where $s=3$,
$K_0=\bar{S}_n$,
$y_0=1+ 2.2 / s^2 $,
$y_1= T_s\left(y_0\right) / T^{\prime}_s\left(y_0\right)$,
and $T_j$ denotes the $j$-th Chebyshev polynomial,
that is,  
$T_0\left( y \right)=1$, $T_1\left( y \right)= y$
and
$T_j\left( y \right)=2 y T_{j-1}\left( y\right)-T_{j-2}\left( y \right)$.

\begin{figure}[bt]
\centering
   \includegraphics[height= 3.5in,width=5.7in]{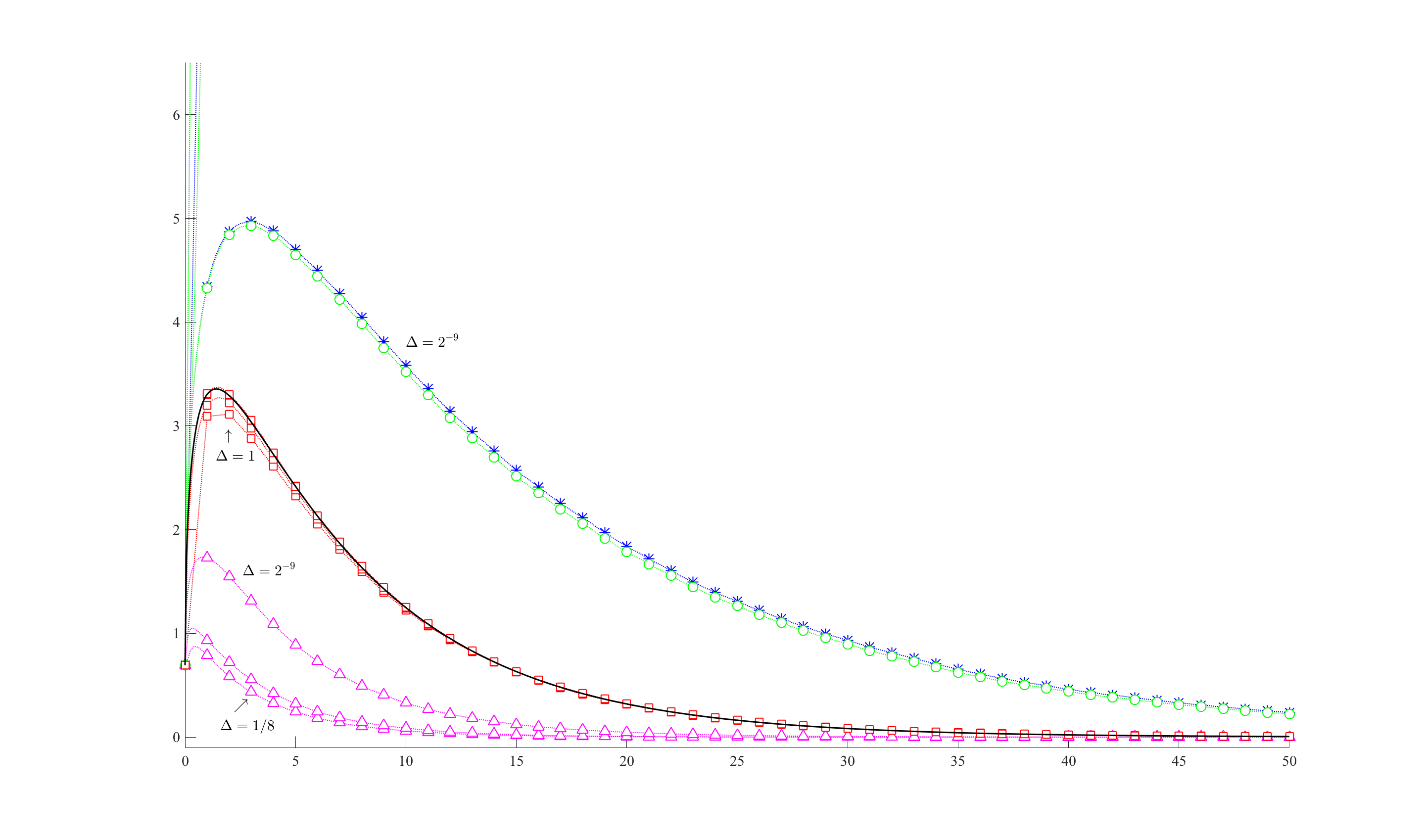}
  \caption{Computation of $\mathbb{E}  \log \left( 1 + \left( X_t^1 \right)^2 \right)$,
where $t\in\left[0,50\right]$ and 
$X_t$ solves \eqref{eq:1.16} with $b = -4$, $\sigma = \epsilon = 8$ and $X_0=\left(1,2\right)^{\top}$. 
Scheme \ref{scheme:EulerStableBilinear}, the backward Euler method, the balanced scheme and  SROCK scheme
are represented by squares, stars, triangles and circles, respectively. The true values are plotted with a solid line.
\label{fig:LS1}
}
\end{figure}

First,
we compute $\mathbb{E} \log \left( 1 + \left(X^1_t\right)^2 \right)$ in case $b = -4$, $\sigma = \epsilon = 8$ and $X_0=\left(1,2\right)^{\top}$.
Since  $- \lambda = b+\left(\epsilon^2-\sigma^2\right)/2 < 0$,
$\mathbb{E} \left(  \left\Vert X_t \right\Vert ^p \right)$ converges exponentially to $0$ as $t \rightarrow + \infty$
provided that $p > 0$ is sufficiently small.
We thus get $\lim_{t \rightarrow + \infty} \mathbb{E} \log \left( 1 + \left(X^1_t\right)^2 \right) = 0$,
because
$ \log \left( 1 + x^2 \right) \leq 2 x^p / p$ for all $x \geq 0$ and $p \in \left] 0 , 1 \right[ $.
Figure \ref{fig:LS1} displays the numerical approximations of $\mathbb{E}  \log \left( 1 + \left( X_t^1 \right)^2 \right)$
obtained from sampling $10^6$ times $\bar{E}_{n}$, $\bar{B}_n$, $\bar{S}_{n}$ and Scheme \ref{scheme:EulerStableBilinear}
with step sizes $1/8$, $1/32$ and $2^{-9}=0.002$.
We also include $\bar{X}_n$ (i.e., Scheme \ref{scheme:EulerStableBilinear}) with $\Delta = 1$.
In good agreement with Theorem \ref{th:MomentStable},
the first coordinate of $ \bar{X}_{n}$ decays to $0$ with the same speed that $X_t$ even when $\Delta=1$.
Scheme \ref{scheme:EulerStableBilinear} also reproduces very well the  transient behavior of the exact solution.
Indeed,
it's difficult to distinguish  between $\bar{X}_n$ and $X_{T_n}$ whenever $\Delta=1/8,2^{-5},2^{-9}$.
In contrast, 
the estimates given by $\bar{E}_n$ and the S-ROCK scheme $\bar{S}_n$
grow disproportionately when  $\Delta=1/8,1/32$.
The balanced scheme $\bar{B}_n$ shows a slow speed  of convergence.

\renewcommand{\arraystretch}{1}
\renewcommand\tabcolsep{3pt}
\begin{table}[tb]
\begin{center}
\footnotesize
\def\arraystretch{0.5}
\begin{tabular}{c|c|c|c|c|c|c|c|c|}
\cline{2-9}
&\multicolumn{8}{c|}{$\Delta$}\\
&1/16&1/32&1/64&1/128&1/256&1/512&1/1024&1/2048\\
\cline{1-9}
\multicolumn{1}{|c|}{\phantom{$x^{2^2}$} \hspace{-17pt} $\epsilon_a \left(\bar{X},\Delta \right)$}
& 1.6811e-3& 0.16279 & 3.0648 e-3& 1.3733e-3 & 1.6362e-3 & 1.6405e-3 & 2.3889e-3 & 1.5285e-3
\\
\multicolumn{1}{|c|}{  $\epsilon_a \left(\bar{E},\Delta \right)$}
& 187.3152 & 257.7663 & 279.6334 & 201.0138 & 49.9386 & 2.3838 & 0.31756 & 0.13123
\\
\multicolumn{1}{|c|}{$\epsilon_a \left(\bar{B},\Delta \right)$}
& 2.6598 & 2.5687 & 2.4313 & 2.2477 & 2.0212 &1.7589  & 1.4777 & 1.2023
\\
\multicolumn{1}{|c|}{$\epsilon_a \left(\bar{S},\Delta \right)$}
& 173.0071 & 250.4513 &275.9211 & 199.142 & 49.173 & 2.3248 & 0.30093 & 0.12372
\\
\hline
\end{tabular}
%
%
%
\end{center}
\caption{ Absolute errors involved in the computation of
$\mathbb{E}  \log \left( 1 + \left( X^1_T \right)^2 \right)$,
where $X_t$ satisfies  \eqref{eq:1.16} with $b=-4$, $\sigma=\epsilon=8$ and $X_0=\left(1,2\right)^{\top}$.
}
\label{TableLS1}
\end{table}


Table \ref{TableLS1} lists the errors
\begin{equation}
\label{eq:4.4}
\epsilon_a \left( \widetilde{Y},\Delta\right)
=
\max_{t \in \{0, h, 2h, \ldots, 10\}} 
\left\vert \mathbb{E}\log \left( 1 + \left(X^1_t\right)^2 \right)
-
\mathbb{E}\log \left( 1 + \left(\widetilde{Y}^1_{t/\Delta}\right)^2 \right) \right\vert
\end{equation}
obtained from the sample means of $10^8$ observations of the schemes 
$\widetilde{Y}_n = \left(\widetilde{Y}_n^1,\widetilde{Y}_n^2\right)^{\top}= \bar{X}_n, \bar{E}_{n}, \bar{B}_n, \bar{S}_{n}$,
where $h = \max\{\Delta, 2^{-4}\}$.
In this subsection,
the \textquotedblleft true\textquotedblright  \  values of $X_t$
have been calculated by sampling $10^8$ times the explicit solution of \eqref{eq:1.16}.
The length of the $99\%$ confidence intervals for $\mathbb{E}\log \left( 1 + \left(X^1_t\right)^2 \right)$ and 
$\mathbb{E}\log \left( 1 + \left(\widetilde{Y}^1_{t/\Delta}\right)^2 \right)$
are at least of order $10^{-3}$ (= e-3) except for $\bar{E}$ and $\bar{S}$ 
that are sometimes of order $10^{-2}$ (= e-2);
they  have been estimated as in \cite{Kloeden1992}.
According to Table \ref{TableLS1},
the errors corresponding to Scheme \ref{scheme:EulerStableBilinear} 
are similar to the length of the $99\%$ confidence intervals associated with the sample,
except for $\Delta = 1/32$.

\renewcommand{\arraystretch}{1}
\renewcommand\tabcolsep{3pt}
\begin{table}[tb]
\begin{center}
\footnotesize
\begin{tabular}{c|c|c|c|c|c|c|c|c|c|c|c|}
\cline{2-12}
 &\multicolumn{11}{c|}{$\Delta$}\\
 & $1$ & $2^{-1}$ & $2^{-2}$ & $2^{-3}$ & $2^{-4}$ & $2^{-5}$ & $2^{-6}$ & $2^{-7}$ & $2^{-8}$ & $2^{-9}$ & $2^{-10}$  \\ 
\cline{1-12}
\multicolumn{1}{|c|}{\phantom{$x^{2^2}$} \hspace{-17pt} 
$100 \, \epsilon_r \left(\bar{X},\Delta \right)$}
& 2.1e-2 & 6.8e-4 & 9.7e-3 & 1.6e-2 & 1.6e-2 & 1.1e-3 & 4.6e-3 & 5.6e-3 & 7.4e-3 & 1.6e-2  & 1.3e-2
 \\ 
\multicolumn{1}{|c|}{$100 \, \epsilon_r \left(\bar{E},\Delta \right)$}
& 61.7 & 144 & - & 943 & 784 & 518 & 175 & 10.4 & 6.07 & 3.21 &  1.58
\\ 
\multicolumn{1}{|c|}{$100 \, \epsilon_r \left(\bar{B},\Delta \right)$}
& 8.42 & 24.6 & 43.8 & 65.0 & 84.3 & 98.8 & 106  & 106 & 98.5 & 85.9 & 70.9 
\\ 
\multicolumn{1}{|c|}{$100 \, \epsilon_r \left(\bar{S},\Delta \right)$}
& 336 & 477 & 617 & 711 & 685 & 472 & 154 & 2.95 &  9.58 & 5.03 & 2.51
\\ 
\hline
\end{tabular}
\end{center}
\caption{Relative weak global errors \eqref{eq:4.5} at $T=10$
for  \eqref{eq:1.16} with $b=\sigma=4$, $\epsilon=3$ and $X_0=\left(2,4\right)^{\top}$.
The means values are estimated from the sample means of $10^8$ observations of 
$\bar{X}_n, \bar{E}_{n}, \bar{B}_n, \bar{S}_{n}$.
}
\label{TableLS3}
\end{table}

\begin{figure}[tb]
\centering
     \includegraphics[height= 3.0in,width=5.7in]{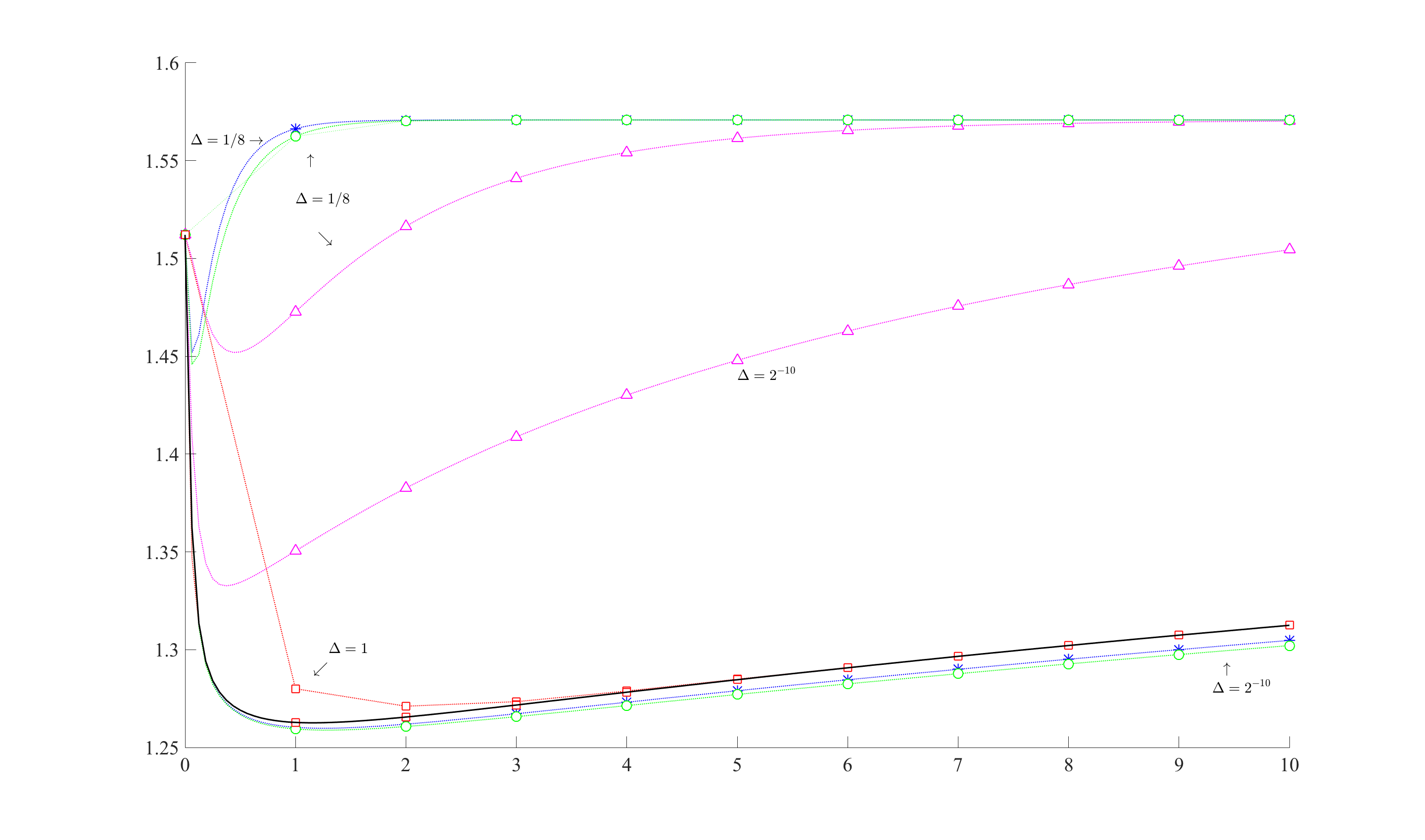}
   \caption{Computation of $\mathbb{E}  \arctan \left( 1 + \left( X_t^2 \right)^2 \right)$,
where $X_t$ solves \eqref{eq:1.16} with $b = \sigma = 4$, $\epsilon = 3$ and $X_0=\left(2,4\right)^{\top}$.
Scheme \ref{scheme:EulerStableBilinear}, the backward Euler method, the balanced scheme and  SROCK scheme
are represented by squares, stars, triangles and circles, respectively.
The true values are plotted with a solid line.
}
\label{fig:LS2}
\end{figure}

Second,
we take $b = \sigma =  4$, $\epsilon =3$ and $X_0=\left(2,4\right)^{\top}$.
In this situation,
\eqref{eq:1.16} satisfies the assumptions of Theorem \ref{th:NoStable}.
Furthermore,
it follows from the proofs of  Theorem 12 of \cite{Appleby2008} and Theorem \ref{th:NoStable} that
the norms of  $X_t$  and $ \bar{X}_{n} $ converge almost surely to $+ \infty$
since $A_n = \left( n+1 \right) \Delta \sigma^2$. 
Figure \ref{fig:LS2} presents the computation of $\mathbb{E}  \arctan \left( 1 + \left( X_t^2 \right)^2 \right)$
by using $\bar{X}_n, \bar{E}_{n}, \bar{B}_n, \bar{S}_{n}$ with step sizes $1/16$ and $2^{-10} = 1/1024$,
as well as Scheme \ref{scheme:EulerStableBilinear} with $\Delta = 1$.
The true values have been plotted with a solid line.
In Figure \ref{fig:LS2},
we can see that Scheme \ref{scheme:EulerStable} replicates very well the growth of $X_t$.
Indeed,
the difference between $\bar{X}_n$ and $X_{T_n}$ is not significant when $\Delta=1/16,2^{-10}$.  
Table \ref{TableLS3} gives the errors 
\begin{equation}
\label{eq:4.5}
\epsilon_r \left( \widetilde{Y}, \Delta \right)
=
\frac{
\left\vert  
\mathbb{E}\log \left( 1 + \left(X^1_T\right)^2 \right)
-
\mathbb{E}\log \left( 1 + \left(\widetilde{Y}^1_{T/\Delta}\right)^2 \right) 
\right\vert
}{
\left\vert \mathbb{E}\log \left( 1 + \left(X^1_T\right)^2 \right)  \right\vert 
} 
\end{equation}
with $T=10$
and 
$\widetilde{Y}_n = \left(\widetilde{Y}_n^1,\widetilde{Y}_n^2\right)^{\top}= \bar{X}_n, \bar{E}_{n}, \bar{B}_n, \bar{S}_{n}$.
Table \ref{TableLS3} shows that the error
$
\left\vert  
\mathbb{E}\log \left( 1 + \left(X^1_T\right)^2 \right)
-
\mathbb{E}\log \left( 1 + \left(\bar{X}^1_{T/\Delta}\right)^2 \right) 
\right\vert
$
is at most $0.01 \%$ of $\left\vert \mathbb{E}\log \left( 1 + \left(X^1_T\right)^2 \right)  \right\vert$
for any $\Delta$. 
Recall that $\bar{X}_n$ stands for  Scheme \ref{scheme:EulerStable}.
The relative weak errors of $\bar{E}_{n}$, $\bar{B}_n$ and $\bar{S}_{n}$ are large.

\renewcommand{\arraystretch}{1}
\renewcommand\tabcolsep{3pt}
\begin{table}[htbp]
\begin{center}
\footnotesize

\begin{tabular}{c|c|c|c|c|c|c|c|c|c|c|}
\cline{2-11}
&\multicolumn{10}{c|}{$\Delta$}\\
&1/16&1/32&1/64&1/128&1/256&1/512&1/1024&1/2048&1/4096&1/8192\\
\cline{1-11}
\multicolumn{1}{|c|}{\phantom{$x^{2^2}$} \hspace{-17pt} 
$\hat{\epsilon} \left(\bar{X},\Delta \right)$}
&2.0003&2.0002&1.9963&1.8715&1.3854&0.72545&0.28044&0.08892&0.02517&6.71e-3\\
\multicolumn{1}{|c|}{\phantom{$x^{2^2}$} \hspace{-17pt} $\hat{\epsilon} \left(\bar{E},\Delta \right)$}&5.94e73&7.73e77 &1.05e69&1.84e56&9.35e33&3.68e16&4.16e7&5.42e3&63.833&6.725\\
\multicolumn{1}{|c|}{$\hat{\epsilon} \left(\bar{B},\Delta \right)$}&3.98e42&3.57e38&3.98e33&1.43e29&3.79e24&5.17e20&4.99e14&5.58e10&2.42e7&5.04e4\\
\multicolumn{1}{|c|}{$\hat{\epsilon} \left(\bar{S},\Delta \right)$}&3.39e72&1.79e77&5.03e68&1.26e56&7.80e33&3.35e16&3.97e7&5.29e3&63.084&6.6864\\
\hline
\end{tabular}

\caption{Relative strong errors appearing in the pathwise solution of  \eqref{eq:1.16} with $b=-4$, $\sigma=\epsilon=8$, $X_0=\left(1,2\right)^{\top}$ and $t \in \left[ 0 , 2 \right]$.
}
\label{TableLS2}

\begin{tabular}{c|c|c|c|c|c|c|c|c|c|c|}
\cline{2-11}
&\multicolumn{10}{c|}{$\Delta$}\\
&1&1/2&1/4&1/8&1/16&1/32&1/64&1/128&1/256&1/512\\
\cline{1-11}
\multicolumn{1}{|c|}{ \phantom{$x^{2^2}$} \hspace{-17pt} 
$\hat{\epsilon}  \left(\bar{X},\Delta \right)$}
&1.9998&2.0001&1.9978&1.2963&0.71472&0.35089&0.14508&5.09e-2&1.572e-2&4.44e-3\\
\multicolumn{1}{|c|}{ $\hat{\epsilon} \left(\bar{E},\Delta \right)$}&8.92e22&3.27e25 &-&3.48e35&1.59e30&1.3e19&7.62e8&1.43e3&16.577&2.6258\\
\multicolumn{1}{|c|}{$\hat{\epsilon} \left(\bar{B},\Delta \right)$}&6.63e17&1.42e16&7.5e15&3.55e13&1.5e12&2.52e10&5.2e8&1.32e7&2.27e5&8.58e3\\
\multicolumn{1}{|c|}{$\hat{\epsilon} \left(\bar{S},\Delta \right)$}&7.17e28&2.67e30&1.2e34&1.31e32&5.63e28&2.69e18&3.52e8&9.772e2&13.7402&2.3963\\
\hline
\end{tabular}
\end{center}
\caption{Relative strong errors appearing in the pathwise solution of   \eqref{eq:1.16} with $b=\sigma=4$, $\epsilon=3$, $X_0=\left(2,4\right)^{\top}$ and and $t \in \left[ 0 , 2 \right]$.
}
\label{TableLS4}
\end{table}

Third, 
Tables \ref{TableLS2} and \ref{TableLS4}  provide  the relative strong errors
$$
\hat{\epsilon}\left(\widetilde{Y},\Delta\right)=\sup_{n=0,\ldots,T/\Delta}
\mathbb{E} \left( \left\Vert X_{T_n} - \widetilde{Y}_n \right\Vert^2 / \left\Vert X_{T_n} \right\Vert^2
\right)
$$
in the above two examples:
(i) $b=-4$, $\sigma=\epsilon=8$, $X_0=\left(1,2\right)^{\top}$;
and
(ii) $b = \sigma = 4$, $\epsilon = 3$,$X_0=\left(2,4\right)^{\top}$.
We take $T=2$ and $ \widetilde{Y}_n = \bar{X}_n, \bar{E}_{n}, \bar{B}_n, \bar{S}_{n}$.
The $99\%$ confidence interval are at least of order $10^{-4}$ (= e-4) for the DND scheme $\bar{X}_n$,
and 
all the sample sizes are equal to $10^8$.
From Tables \ref{TableLS2} and \ref{TableLS4} we can see that  
Scheme \ref{scheme:EulerStableBilinear} reproduces very well the trajectories of $X_t$,
complementing Theorem \ref{th:StrongRateConvergence}.

Finally, 
we discuss the effect of round-off errors on Scheme \ref{scheme:EulerStable} applied to \eqref{eq:1.16}.
Using simple algebraic transformations we get
\begin{equation}
\label{eq:4.23}
 \left\Vert  \bar{Z}_{n+1} \right\Vert ^2 = \left( 1- \epsilon^2 \Delta/2 \right)^2 + \epsilon^2 \Delta \left(\hat{W}^{2}_{n} \right)^2 .
\end{equation}
If $\hat{W}^{2}_{n}$ takes values  $\pm 1$ with probability $1/2$,
then
$\left\Vert  \bar{Z}_{n+1} \right\Vert ^2 = 1 + \epsilon^4 \Delta^2 / 4$,
and so we can calculate  
$  \bar{Z}_{n+1} /  \left\Vert  \bar{Z}_{n+1} \right\Vert $
without problems.
From \eqref{eq:4.23} it follows that 
$\bar{Z}_{n+1} \approx 0$ 
if and only if   $\Delta \approx  2 / \epsilon^2$ and  $\hat{W}^{2}_{n} \approx 0$.
The latter happens with an extremely  low probability in case $\hat{W}^{2}_{n} $ is uniformly distributed on $\left[ - \sqrt{3}, \sqrt{3} \right]$.
If $\hat{W}^{2}_{n} $ is obtained by means of a normal pseudorandom number generator,
then
$\hat{W}^{2}_{n} \approx 0$ with small probability,
and hence the performance of Scheme \ref{scheme:EulerStable} is not affected when  $\Delta =  2 / \epsilon^2 = 2^{-5}$,
as we can see in Figure \ref{fig:LS1}.

\subsection{Scalar SDE}
\label{sec:NumExpScalar}
This subsection examines the behavior of Scheme \ref{scheme:EulerStable}
applied to  the stochastic Ginzburg-Landau equation
\begin{equation}
\label{eq:5.1}
X_{t}
=
X_0
+
\int_{0}^{t} \left(  a \, X_s -  b \, \left( X_s \right) ^3  \right)  ds
 + 
 \int_{0}^{t} \sigma \, X_{s} \, dW^{1}_{s} ,
\end{equation}
which constitutes a classical test  in the theory of stochastic bifurcation 
(see, e.g., \cite{Arnold1998,Baxendale1994}).
Here, $X_{t}$ takes values in $ \mathbb{R}$,
$a \in \mathbb{R}$ and $b, \sigma > 0$.
We compute $\mathbb{E}  \phi \left(X_T\right)$ for the test function $\phi\left(x\right):=\log \left( 1 + x^2 \right)$ in the following situations:
\textbf{(Ex1)} $a=b=1$, $\sigma =2$, $X_0 = 1$, $T=5$; 
\textbf{(Ex2)} $a = 6$, $b = 9$, $\sigma = 3$,  $X_0 = 1$, $T=10$; and
\textbf{(Ex3)} $a = 9$, $b = 1$, $\sigma = 4$,  $X_0 = 10^{-6}$, $T=20$.
%
%
%
In Example \textbf{Ex1}, the motivating problem of \cite{Higham2007},
$
 b\left(x\right)/x-\left(\sigma^1\left(x\right)/x \right)^{2}  /2
 \leq -1
$,
and hence
$\mathbb{E}  \log \left( 1 + \left(X_t \right)^2 \right) $
converges exponentially to $0$ as $t \rightarrow + \infty$.
In the test problems \textbf{Ex2} and \textbf{Ex3}, 
(\ref{eq:5.1}) has three invariant forward Markov measures (see, e.g., p. $480$ in \cite{Arnold1998}).

\renewcommand{\arraystretch}{1}
\renewcommand\tabcolsep{3pt}
\begin{table}[htbp]
\begin{center}
\footnotesize
\begin{tabular}{cc|c|c|c|c|c|c|c|c|}
\cline{3-10}
& &\multicolumn{8}{c|}{$\Delta$}\\
& &1&1/2&1/4&1/8&1/16&1/32&1/64&1/128\\ 
\cline{1-10}
\multicolumn{1}{|c}{\multirow{3}{*}{
\begin{sideways} {\bf Ex1} \end{sideways}}} & 
\multicolumn{1}{c|}{\phantom{$x^{2^2}$} \hspace{-17pt} $\epsilon_r \left(\bar{X},\Delta \right)$}&0.075304&0.18365&0.11585&0.063313&0.033357&0.016634&0.0071119&0.0023089\\
\multicolumn{1}{|c}{}&
\multicolumn{1}{c|}{$\epsilon_r \left(\hat{E},\Delta \right)$}&105.7917&70.0728&14.5311&0.078833&0.4021&0.17222&0.085249&0.041843\\
\multicolumn{1}{|c}{}&
\multicolumn{1}{c|}{$\epsilon_r \left(\bar{B},\Delta \right)$}&14.7149&11.6118&8.62&6.1368&4.2603&2.9208&1.9922&1.3595\\
\hline
\hline
\multicolumn{1}{|c}{\multirow{3}{*}{
\begin{sideways} {\bf Ex2} \end{sideways}
}} &
\multicolumn{1}{c|}{\phantom{$x^{2^2}$} \hspace{-17pt} $\epsilon_r \left(\bar{X},\Delta \right)$}&0.6434&0.39757&0.46459&0.23047&0.089759&0.039578&0.018567&0.0087071\\
\multicolumn{1}{|c}{}&
\multicolumn{1}{c|}{$\epsilon_r \left(\hat{E},\Delta \right)$}&111.4751&118.4597&65.1188&6.0323&0.097374&0.25646&0.12314&0.056297\\
\multicolumn{1}{|c}{}&
\multicolumn{1}{c|}{$\epsilon_r \left(\bar{B},\Delta \right)$}&3.6856&3.4593&3.0512&2.6293&2.169&1.7251&0.12314&1.0222\\
\hline
\hline
\multicolumn{1}{|c}{\multirow{3}{*}{
\begin{sideways} {\bf Ex3} \end{sideways}
}} &
\multicolumn{1}{c|}{\phantom{$x^{2^2}$} \hspace{-17pt} $\epsilon_r \left(\bar{X},\Delta \right)$}&0.31386&0.60358&0.66093&0.51422&0.26152&0.11126&0.049885&0.023443\\
\multicolumn{1}{|c}{}&
\multicolumn{1}{c|}{$\epsilon_r \left(\hat{E},\Delta\right)$}&94.3589&141.2312&141.0024&47.1114&1.0046&0.7024&0.38043&0.16122\\
\multicolumn{1}{|c}{}&
\multicolumn{1}{c|}{$\epsilon_r \left(\bar{B},\Delta \right)$}&8.026&7.8669&7.3796&6.7487&5.9521&5.0546&4.1603&3.3413\\
\hline
\end{tabular}
\end{center}
\caption{Relative errors involved in the computation of
$\mathbb{E}  \log \left( 1 + \left( X_T \right)^2 \right)$ by using the schemes $\bar{X}_n, \hat{E}_n$ and $\bar{B}_n$, each with uniform distributed random variables.
}
\label{Tabla1}
\end{table}

We solve  (\ref{eq:5.1}) using 
Scheme \ref{scheme:EulerStableScalar} with $\hat{W}^{k}_{n}$ distributed uniformly on $\left[ - \sqrt{3}, \sqrt{3} \right]$,
the balanced scheme \eqref{scheme:Balanced},
and the weak version of the tamed Euler scheme  \cite{Hutzenthaler2012}
\begin{equation*}
\hat{E}_{n+1}
=
\hat{E}_n+ \left( b\left(\hat{E}_n\right) \Delta 
\right)/\left( 1+\left\vert b\left(\hat{E}_n\right)\right\vert\Delta \right)
+\sum_{k=1}^m\sigma^{k}\left(\hat{E}_n\right)\sqrt{\Delta} \hat{W}^k_n .
\end{equation*}
From \cite{Higham2007} we have that
the Euler-Maruyama scheme applied to (\ref{eq:5.1}) 
blows up, with positive probability, at a geometric rate.
Table \ref{Tabla1} presents the relative errors 
$$
\epsilon_r \left( \widetilde{Y}, \Delta \right)
=
\left\vert \mathbb{E}  \phi \left( X_T \right) - \mathbb{E} \phi \left( \widetilde{Y}_{n_T} \right) \right\vert
/
\left\vert \mathbb{E}  \phi \left( X_T \right)  \right\vert ,
$$
where $\Delta \, n_T= T$ and $\widetilde{Y}_n$ stands for
$\bar{X}_n$ (Scheme \ref{scheme:EulerStableScalar}), $\hat{E}_n$ and $\bar{B}_n$.
We have computed 
$\mathbb{E} \phi \left( \widetilde{Y}_{n_T} \right)$
by Monte Carlo simulation with $10^8$ realizations of $\widetilde{Y}_n$.
Similarly, the reference values for  $\mathbb{E}  \phi \left( X_{T} \right)$ have been obtained by sampling $10^8$ times 
the backward Euler method \eqref{scheme:BackwardEuler}
with $\mathbb{P} \left( \hat{W}_n^k = \pm 1 \right) = 1/2$ and $\Delta=2^{-11}$.
Following \cite{Kloeden1992},
we estimated that the length of the $99\%$ confidence intervals for the \textquotedblleft true\textquotedblright \  values
are at least of order $10^{-4}$.
In Examples \textbf{Ex1} and  \textbf{Ex2}, 
$\mathbb{E}  \phi \left( X_{t} \right)$ decreases as $t \rightarrow T$ to 
$0.035434702074070$ and $0.120956832017214$, respectively.
On the other hand,
$\mathbb{E}  \phi \left( X_{t} \right)$ grows to $0.260641703098477$ as $t \rightarrow T$ in Example \textbf{Ex3}. 
We select $T$ so that 
$\mathbb{E}  \phi \left( X_{T} \right) \approx \lim_{t \rightarrow \infty } \mathbb{E}  \phi \left( X_{t} \right) $.

Table \ref{Tabla1} shows the very good accuracy of Scheme \ref{scheme:EulerStableScalar}  applied to (\ref{eq:5.1}).
In the three examples,
$\epsilon_r \left(\bar{X},\Delta \right) < \epsilon_r \left(\bar{B},1/218 \right)$ 
and $\epsilon_r \left(\hat{E},\Delta \right)$ takes large values when $\Delta$ is not small.

\subsection{System of non-linear SDEs}
\label{sec:NumExpNLinearSystem}

We  compare the performance of Scheme \ref{scheme:EulerStable}, 
the Euler method \eqref{scheme:BackwardEuler}
and the balanced scheme \eqref{scheme:Balanced} in the non-linear multidimensional setting,
for which $ \sqrt{\Delta} \hat{W}^{k}_{n+1} = W^k_{T_{n+1}} - W^k_{T_{n}} $.
To this end,
we compute 
$\mathbb{E}  \log \left( 1 + \left( X^1_t \right)^2 \right) $
where  
\begin{equation}
\label{eq:4.1} 
d \begin{pmatrix} X^1_{t} \\ X^2_{t} \end{pmatrix} 
=
 a  \sqrt{2 + \cos \left( X^1_{t} \right) }  \begin{pmatrix} X^1_{t} \\ X^2_{t} \end{pmatrix}  \, dW^{1}_{t}
+
 b \sqrt{2 + \sin \left( X^2_{t} \right) }  \begin{pmatrix} - X^2_{t} \\ X^1_{t} \end{pmatrix} \, dW^{2}_{t} ,
\end{equation}
with $X_0 = \left( 4 , 2 \right) ^{\top}$ and  $a,b \in \mathbb{R}$.
We consider two situations:
{\bf (Ex4)} $a = 6$, $b = 3$;
and 
{\bf (Ex5)} $a = 2.5$, $b = 5$.

\begin{figure}[htbp]
\centering
 \includegraphics[height= 2.5in,width=5.7in]{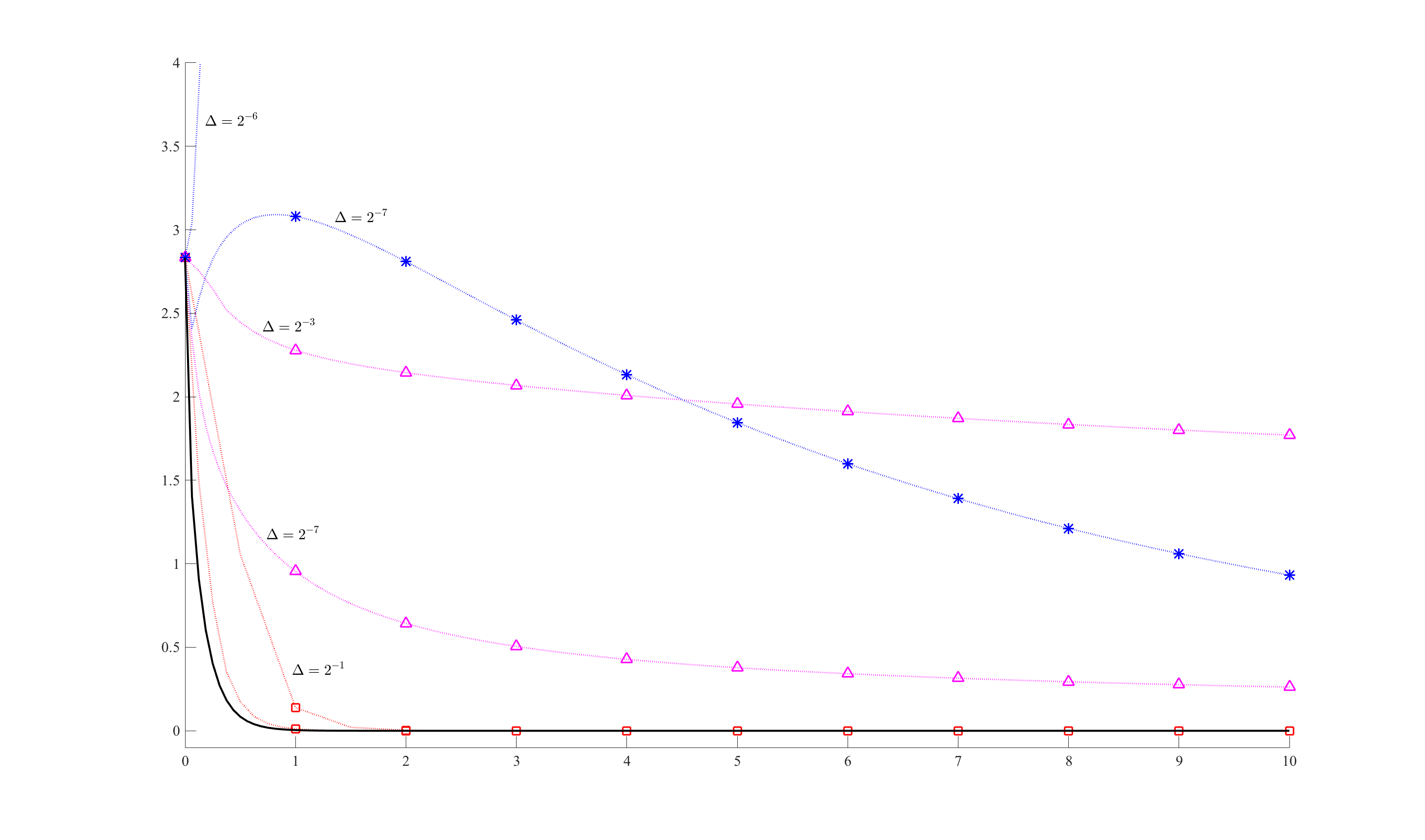}
    \caption{Computation of $\mathbb{E} \log\left(1 + \left(X^1_t\right)^2\right)$
	where $X_t$ solves \eqref{eq:4.1} with $a = 6$ and $b = 3$.
	The reference values for $\mathbb{E} \log\left(1 + \left(X^1_t\right)^2\right)$ are shown with a solid line.
	Stars stand for the backward Euler method, 
	triangles for the balanced scheme \eqref{scheme:Balanced},
	and squares for Scheme \ref{scheme:EulerStable}
	with larger values of $\Delta$, namely $\Delta = 2^{-1}$ and $\Delta = 2^{-3}$.
	}
	\label{fig:ex4}
\end{figure}

\renewcommand{\arraystretch}{1}
\renewcommand\tabcolsep{3pt}
\begin{table}[htbp]
\begin{center}
\footnotesize
\begin{tabular}{c|c|c|c|c|c|c|c|c|}
\cline{2-9}
 &\multicolumn{8}{c|}{$\Delta$}\\
 & 1 & $2^{-1}$ & $2^{-2}$ & $2^{-3}$ & $2^{-4}$ & $2^{-5}$ & $2^{-6}$ & $2^{-7}$\\ 
\cline{1-9}
\multicolumn{1}{|c|}{\phantom{$x^{2^2}$} \hspace{-17pt} $\epsilon_a \left(\bar{X},\Delta \right)$}
& 0.50795 & 0.13444 & 0.029171  & 6.97309e-3 & 1.85285e-3 & 1.0012e-4 & 4.3949e-4 & 6.5198e-4\\
\multicolumn{1}{|c|}{$\epsilon_a \left(\bar{E},\Delta \right)$} & 39.303 & 63.392 & 98.320 & 142.858 & 185.801 & 98.692 & 117.436 & 3.07578\\ 
\multicolumn{1}{|c|}{$\epsilon_a \left(\bar{B},\Delta \right)$} & 2.7791 & 2.67505 & 2.49069 & 2.27442 & 1.97407 & 1.60898 & 1.24567 & 0.95072\\ 
\hline
\end{tabular}
\caption{Weak errors \eqref{eq:4.4} appears in the numerical solution of Example {\bf Ex4} by 
$\bar{X}_n, \bar{E}_n$ and $\bar{B}_n$.}
\label{TablaEx4}
\end{center}
\end{table}

In Example {\bf Ex4},
$$
- \lambda 
=
\frac{1}{2} \sup_{ \left( x, y \right) \neq 0} 
\bigl(
b^2 \left( 2 + \sin \left( y \right) \right) - a^2 \left( 2 + \cos \left( x \right) \right)
\bigr)
\leq \frac{3b^2 - a^2}{2} < 0 ,
$$
and so $0$ is a stable equilibrium point of \eqref{eq:4.1} in both senses, almost sure and small moment.
In concordance with Theorem  \ref{th:MomentStable},
Figure \ref{fig:ex4} shows that the dynamics of $\bar{X}_n$ resembles that of  $X_t$.
Table \ref{TablaEx4} reveals the much higher precision of Scheme \ref{scheme:EulerStable}  for all values of $\Delta$.

\begin{figure}[tb]
\centering
 \includegraphics[height= 2.5in,width=5.7in]{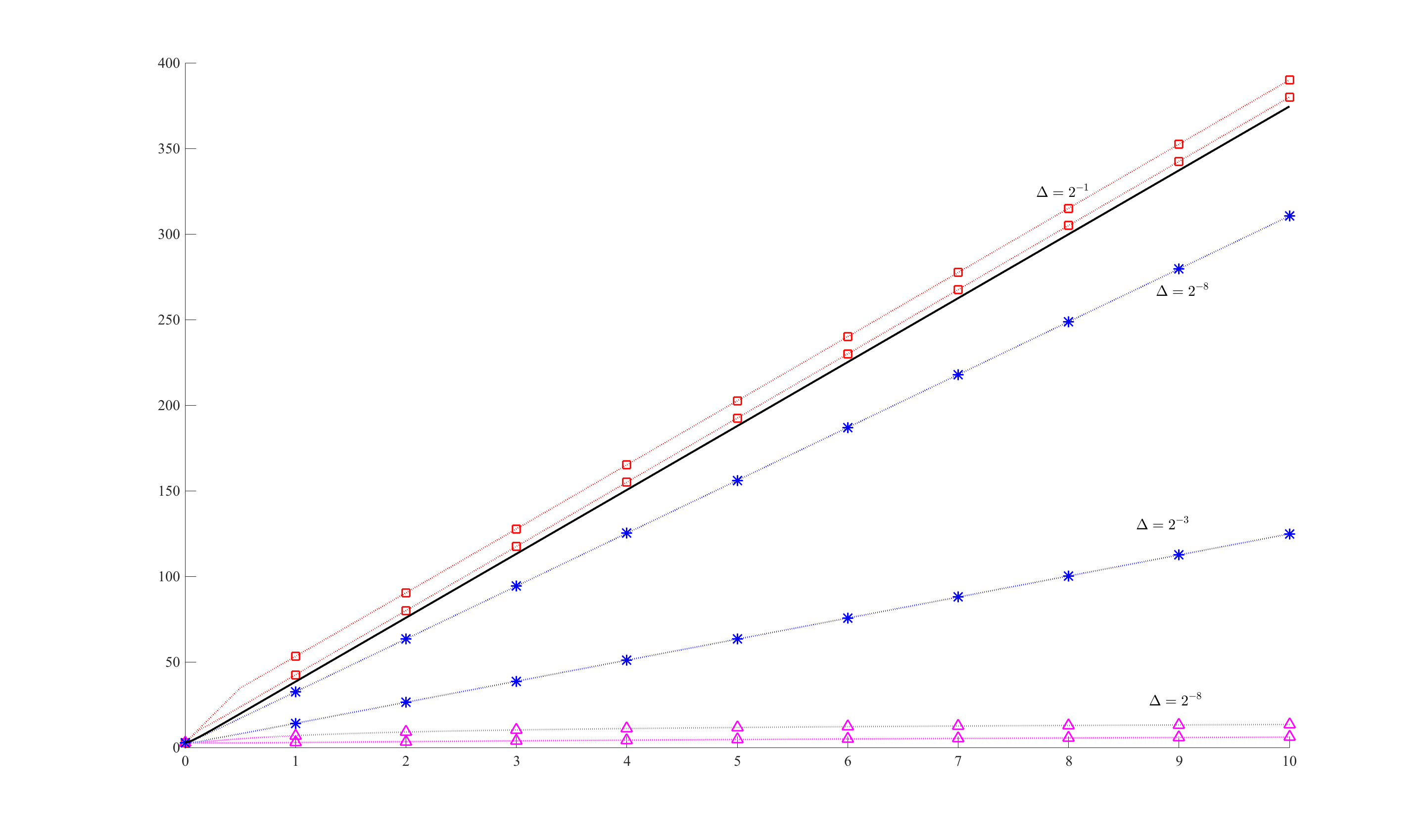}
    \caption{Computation of $\mathbb{E} \log\left(1 + \left(X^1_t\right)^2\right)$
	where $X_t$ solves \eqref{eq:4.1} with $a = 2.5$ and $b = 5$.
	The reference values for $\mathbb{E} \log\left(1 + \left(X^1_t\right)^2\right)$ are shown with a solid line.
	Scheme \ref{scheme:EulerStable}, the backward Euler method and the balanced scheme are represented by squares, stars and triangles, respectively.}
	\label{fig:ex5}
 \end{figure}

We turn to Example {\bf Ex5}. Since
\begin{equation*}
\frac{
\langle x, b \left( x \right)  \rangle
+
\frac{1}{2} \sum_{k=1}^{m}  \left\Vert \sigma^k \left( x \right)  \right\Vert^2 
 }{
 \left\Vert x \right\Vert^2 
 }
- \left( 1 + \frac{1}{3 a^2} \right) 
\frac{ 
\sum_{k=1}^{m}
\left\langle x,  \sigma^k \left( x \right)  \right\rangle^2
}{
\left\Vert x \right\Vert^4
}
\geq 
5 ,
\end{equation*}
the assumptions of Theorem \ref{th:NoStable} are satisfied.
From the proofs of  Theorem 12 of \cite{Appleby2008} and Theorem \ref{th:NoStable}
we deduce that $\left\Vert X_t \right\Vert$ and $\left\Vert \bar{X}_{n}  \right\Vert $ converge almost surely to $+ \infty$,
because $A_n \geq \left( n+1 \right) \Delta a^2$.
Figure \ref{fig:ex5} shows that 
$\mathbb{E} \log\left(1 + \left(\bar{X}^1_n\right)^2\right)$
goes to $+ \infty$ with the same speed of $ t \mapsto \mathbb{E} \log\left(1 + \left(X^1_t\right)^2\right)$,
while the balanced scheme presents a rigid behavior. 
Curiously,
for large $\Delta$, 
the Euler method \eqref{scheme:BackwardEuler} goes fast to $+ \infty$ in Example {\bf Ex4},
whereas it growths too slow in Example {\bf Ex5}.
From Table \ref{TablaEx5} we have that
the weak  error 
$
\left\vert  
\mathbb{E}\log \left( 1 + \left(X^1_T\right)^2 \right)
-
\mathbb{E}\log \left( 1 + \left(\bar{X}^1_{T/\Delta}\right)^2 \right) 
\right\vert
$
of Scheme \ref{scheme:EulerStable}
is around $1 \%$ of 
$\left\vert  
\mathbb{E}\log \left( 1 + \left(X^1_T\right)^2 \right)
\right\vert
$
for $\Delta = 1/8, 1/16$.
For the balanced scheme,   
$
\epsilon_r \left(\bar{B},\Delta \right) \approx 1
$.

\renewcommand{\arraystretch}{1}
\renewcommand\tabcolsep{3pt}
\begin{table}[htbp]
\begin{center}
\footnotesize
\begin{tabular}{c|c|c|c|c|c|c|c|c|c|c|c|}
\cline{2-12}
 &\multicolumn{11}{c|}{$\Delta$}\\
 & $1$ & $2^{-1}$ & $2^{-2}$ & $2^{-3}$ & $2^{-4}$ & $2^{-5}$ & $2^{-6}$ & $2^{-7}$ & $2^{-8}$ & $2^{-9}$ & $2^{-10}$ \\ 
\cline{1-12}
\multicolumn{1}{|c|}{\phantom{$x^{2^2}$} \hspace{-17pt} $100 \, \epsilon_r \left(\bar{X},\Delta \right)$}
& 7.77 & 4.13 & 2.34 & 1.45 & 1.00 & 0.78 & 0.66 & 0.58 &  0.55 & 0.52 & 0.50
 \\ 
\multicolumn{1}{|c|}{$100 \, \epsilon_r \left(\bar{E},\Delta \right)$}
& 90.3 & 84.7 & 76.8 & 66.6 & 55.5 & 45.1 & 35.9 & 26.4 & 17.0 & 9.98 & 5.37
\\ 
\multicolumn{1}{|c|}{$100 \, \epsilon_r \left(\bar{B},\Delta \right)$}
& 99.0 & 98.8 & 98.6 & 98.3 & 97.9 & 97.5 & 97.0 &  96.7 & 96.4 & 96.1 & 95.8
\\ 
\hline
\end{tabular}
\caption{
Relative weak global errors \eqref{eq:4.5} at $T=10$ in {\bf Ex5}.
}
\label{TablaEx5}
\end{center}
\end{table}

\subsection{SDE with no equilibrium point at $0$}

\begin{figure}[tb]
\centering
    \includegraphics[height= 2.5in,width=5.7in]{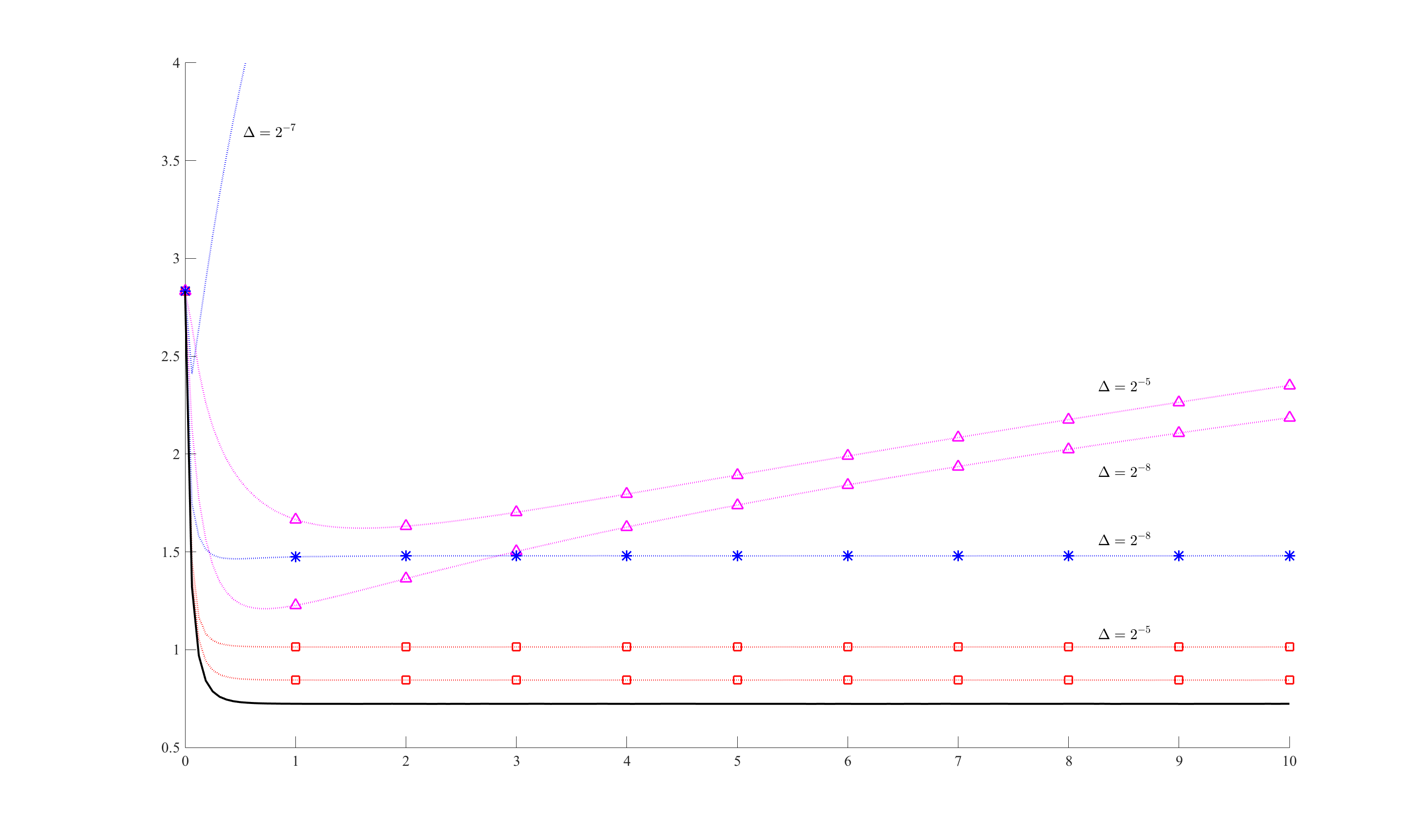}
    \caption{Computation of $\mathbb{E} \log\left(1 + \left(X^1_t\right)^2\right)$
	where $X_t$ solves \eqref{eq:4.3}.
	The reference values are shown with a solid line.
	Scheme \ref{scheme:EulerStableG} ($\Delta = 2^{-5},2^{-6}$), 
	the backward Euler method ($\Delta = 2^{-7}, 2^{-8} $) 
	and the balanced scheme ($\Delta = 2^{-5}, 2^{-8} $) are represented by squares,  stars and triangles, respectively.
	}\label{fig:Ex6}
 \end{figure}

\renewcommand{\arraystretch}{1}
\renewcommand\tabcolsep{3pt}
\begin{table}[tb]
\begin{center}
\footnotesize
\begin{tabular}{c|c|c|c|c|c|c|c|c|}
\cline{2-9}
 &\multicolumn{8}{c|}{$\Delta$}\\
 & $2^{-2}$ & $2^{-3}$ & $2^{-4}$ & $2^{-5}$ & $2^{-6}$ & $2^{-7}$ & $2^{-8}$ & $2^{-10}$ \\ 
\cline{1-9}
 \multicolumn{1}{|c|}{\phantom{$x^{2^2}$} \hspace{-17pt} 
$\epsilon \left(\bar{X},\Delta \right)$}& 0.82459 
& 0.69621 & 0.49943 & 0.29090 & 0.12168 & 0.050324 & 0.015892 & 7.5456e-3\\ 
\multicolumn{1}{|c|}{$\epsilon \left(\bar{E},\Delta\right)$}&  97.6195 & 142.156 & 
185.1003 & 195.013 & 117.509 & 6.73808 & 0.75687 & 0.0106861\\ 
\multicolumn{1}{|c|}{$\epsilon \left(\bar{B},\Delta \right)$}& 1.99544 & 1.907807 & 
1.78435 & 1.627503 & 1.486202 & 1.42107 & 1.46227 & 1.69665 \\ 
\hline
\end{tabular}
\caption{ Weak errors \eqref{eq:4.4} involved in the numerical solution of \eqref{eq:4.3}
by means of $\bar{X}_n, \bar{E}_n$ and $\bar{B}_n$.}
\label{TablaEx6}
\end{center}
\end{table}

In order to evaluate the performance of Scheme \ref{scheme:EulerStableG} 
we consider the test problem 
\begin{equation}
\label{eq:4.3}
\begin{aligned}
 X_{t}
 & =
\begin{pmatrix} 4 \\  2 \end{pmatrix} 
+
\int_{0}^t
\left( 6  \sqrt{2 + \cos \left( X^1_{s} \right) } X_{s} + \begin{pmatrix} 2 \\  0.5 \end{pmatrix} \right)  dW^{1}_{s}
\\ 
& \quad
 +
\int_{0}^t  \left( 3 \sqrt{2 + \sin \left( X^2_{s} \right) }  \begin{pmatrix} - X^2_{s} \\ X^1_{s} \end{pmatrix} 
+ \begin{pmatrix} 1 \\ -1 \end{pmatrix}  \right) dW^{2}_{s} ,
\end{aligned}
\end{equation}
which arises from adding the vectors $\left( 2, 0.5\right)^{\top}$ and  $\left( 1, -1\right)^{\top}$
to the diffusion coefficients of Example {\bf Ex4}.

Figure \ref{fig:Ex6} compares the behavior of Scheme \ref{scheme:EulerStableG}, 
the Euler method \eqref{scheme:BackwardEuler} and the balanced scheme \eqref{scheme:Balanced}.
We estimate the mean values of $\log\left(1 + \left(X^1_t\right)^2\right)$
by sampling $10^8$ times the schemes $\bar{X}_n, \bar{E}_{n}, \bar{B}_n$ for different $\Delta$'s. 
In Figure \ref{fig:Ex6},
we can see that the Euler method \eqref{scheme:BackwardEuler} blows up 
with step size $2^{-7}$,
as well as the incorrect growth of the estimations of 
$\left(\mathbb{E} \log\left(1 + \left(X^1_t\right)^2\right) \right)_{t \in \left[ 2, 10 \right]}$
provided by the balanced scheme \eqref{scheme:Balanced} with $\Delta = 2^{-5}, 2^{-8} $.
We can observe that
Scheme \ref{scheme:EulerStableG} converges to $\mathbb{E} \log\left(1 + \left(X^1_t\right)^2\right)$
reproducing the dynamics of $t \mapsto \mathbb{E} \log\left(1 + \left(X^1_t\right)^2\right)$.
Table \ref{TablaEx6} displays the errors  \eqref{eq:4.4}.
According to Table \ref{TablaEx6} and Figure \ref{fig:Ex6},
the Euler method \eqref{scheme:BackwardEuler} has large errors when $\Delta  \geq 2^{-7}$,
and the balanced scheme \eqref{scheme:Balanced} presents a low speed of convergence.
Scheme \ref{scheme:EulerStableG} is particularly accurate.

\section{Proofs}
\label{Sec:Proofs}

\subsection{Proof of Lemma \ref{lem:ZnoNulo}}
Fix $n \in \mathbb{Z}_+$, $\eta \in \mathbb{R}$ and $x \in \mathbb{R}^d$ satisfying $\left\Vert x \right\Vert = 1$.
To obtain a contradiction, 
suppose that there exists an event $A$ having positive probability such that 
for all $ \omega \in A$,
\begin{equation}
 \label{eq:2p.3}
\begin{aligned}
 0 & = x 
 + \left( \bar{ b } \left( \eta , x \right) - \left\langle x , \bar{ b } \left( \eta , x \right) \right\rangle x
 + \Psi  \left( \eta , x \right)
 \right) \Delta 
 \\ 
 & \quad 
 + \sum_{k=1}^m  \left( \bar{ \sigma}^k \left( \eta  , x \right)  
- \left\langle  x ,  \bar{ \sigma}^k \left( \eta  , x \right) \right\rangle  x \right)
 \sqrt{\Delta}  \hat{W}^{k}_{n} \left( \omega \right) .
\end{aligned}
\end{equation}
Let $\omega_0 \in A$. 
By \eqref{eq:2p.3},  
\begin{equation}
\label{eq:2p.4}
\sum_{k=1}^m  \left( \bar{ \sigma}^k \left( \eta  , x \right)  
- \left\langle  x ,  \bar{ \sigma}^k \left( \eta  , x \right) \right\rangle  x \right)
\left( \hat{W}^{k}_{n} \left( \omega \right) -  \hat{W}^{k}_{n} \left( \omega_0 \right) \right)
= 0
\hspace{1cm}
\forall \omega \in A .
\end{equation}
Since 
$
\left\{ \hspace{-2pt}
\left( \hat{W}^{1}_{n} \left( \omega \right) -  \hat{W}^{1}_{n} \left( \omega_0 \right) ,
\ldots,  \hat{W}^{m}_{n} \left( \omega \right) -  \hat{W}^{m}_{n} \left( \omega_0  \right) \right) ^{\top}
\hspace{-3pt}:
\omega \in A
\hspace{-1pt} \right\}
$
has  positive 
Lebesgue measure,
this subset  contains a basis of $\mathbb{R}^m$.
Then,  \eqref{eq:2p.4}  leads to
$$
 \bar{ \sigma}^k \left( \eta  , x \right)  
- \left\langle  x ,  \bar{ \sigma}^k \left( \eta  , x \right) \right\rangle  x
=
0
\hspace{1cm} \forall k=1,\ldots,m ,
$$
which implies $ \Psi  \left( \eta , x \right) = 0$.
Therefore, 
from (\ref{eq:2p.3}) it follows that 
$$
x + \left(\bar{ b } \left( \eta , x \right) 
- \left\langle x , \bar{ b } \left( \eta , x \right) \right\rangle x \right) \Delta
=
0 ,
$$
and so using $\left\Vert x \right\Vert =1$ we get
$
\left\Vert x \right\Vert ^2
=
\left\langle  x , 
 \left( \left\langle x , \bar{ b } \left( \eta , x \right) \right\rangle x  - \bar{ b } \left( \eta , x \right)  \right) \Delta
\right\rangle
=
0
$,
contrary to $\left\Vert x \right\Vert = 1$.

We have got that almost surely
$$ x 
 + \left( \bar{ b } \left( \eta , x \right) - \left\langle x , \bar{ b } \left( \eta , x \right) \right\rangle x
 + \Psi  \left( \eta , x \right)
 \right) \Delta 
 + \sum_{k=1}^m  \left( \bar{ \sigma}^k \left( \eta  , x \right)  
- \left\langle  x ,  \bar{ \sigma}^k \left( \eta  , x \right) \right\rangle  x \right)
 \sqrt{\Delta}  \hat{W}^{k}_{n}  
$$
is different from $0$. 
Using regular conditional distributions yields $\bar{Z}_{n+1} \neq 0$ a.s.
\endproof

\subsection{Proof of Theorem \ref{th:Stable}}
Since $\bar{\eta}_{0} >0$,
according to  (\ref{eq:2.5}) we have that $\bar{\eta}_{n} >0$ for all $n \in \mathbb{Z}_+$.
Then,
using (\ref{eq:2.5}) yields   
\begin{equation}
\label{eq:2p.1}
\log \left(\bar{\eta}_{n+1}\right)
=
\log \left(\bar{\eta}_{0}\right) 
+  \Delta \sum_{j=0}^n \mu_j 
+ S_n ,
\end{equation}
where
$
S_n = \sum_{j=0}^n  \sum_{k=1}^m 
\left\langle \hat{Z}_j  ,\bar{ \sigma }^k \left(\bar{\eta}_j  , \hat{Z}_j \right) \right\rangle \sqrt{\Delta} \hat{W}^k_{j+1} 
$.
Combining $\left\Vert \hat{Z}_j  \right\Vert  = 1$, (\ref{eq:StabCond}) and (\ref{eq:2p.1}) 
we obtain
\begin{equation}
\label{eq:2p.2}
 \frac{1}{n+1} \log \left(\bar{\eta}_{n+1}\right)
 \leq
 \frac{1}{n+1} \log \left(\bar{\eta}_{0}\right) -  \lambda \Delta + \frac{1}{n+1} S_n  .
\end{equation}

By $\left\Vert \hat{Z}_j  \right\Vert  = 1$ 
and 
$ \left\Vert  \sigma^{k} \left(x\right) \right\Vert   \leq  K  \left\Vert x \right\Vert$ 
for all $x \in \mathbb{R}^d$,
\begin{align*}
 \mathbb{E} \left( 
\sum_{k=1}^m
\left\langle \hat{Z}_j  ,\bar{ \sigma }^k \left(\bar{\eta}_j  , \hat{Z}_j \right) \right\rangle \sqrt{\Delta} \hat{W}^k_{j+1}
\right)^2
& =
\Delta \sum_{k=1}^m
 \mathbb{E} \left(\left\langle \hat{Z}_j  ,\bar{ \sigma }^k \left(\bar{\eta}_j  , \hat{Z}_j \right) \right\rangle ^2 \right)
 \\
 & \leq
\Delta  \sum_{k=1}^m 
 \mathbb{E} \left( \left\Vert \bar{ \sigma }^k \left(\bar{\eta}_j  , \hat{Z}_j \right) \right\Vert ^2 \right)
\leq K \Delta .
\end{align*}
Therefore
$
\sum_{j=0}^{\infty}
\frac{1}{\left( j+1 \right)^2} 
\mathbb{E} \left( 
\sum_{k=1}^m
\left\langle \hat{Z}_j  ,\bar{ \sigma }^k \left(\bar{\eta}_j  , \hat{Z}_j \right) \right\rangle \sqrt{\Delta} \hat{W}^k_{j+1}
\right)^2
< \infty 
$.
The conditional expectation of 
$
\sum_{k=1}^m
\left\langle \hat{Z}_j  ,\bar{ \sigma }^k \left(\bar{\eta}_j  , \hat{Z}_j \right) \right\rangle \sqrt{\Delta} \hat{W}^k_{j+1}
$
given the $\sigma$-algebra generated by
$\bar{\eta}_{0}, \hat{Z}_0  , \hat{W}^1_1, \ldots, \hat{W}^m_1, \ldots, \hat{W}^1_{j}, \ldots, \hat{W}^m_{j}$
is equal to $0$,
and so applying a generalized law of large numbers 
we deduce that $S_n / \left( n+1 \right) \rightarrow 0$ a.s. (see, e.g., p. 243 of \cite{Feller1971}).
Letting $n \rightarrow \infty$ in (\ref{eq:2p.2}) we obtain \eqref{eq:7.31}.
\endproof

\subsection{Proof of Theorem \ref{th:MomentStable}} 
Let $p > 0$.
Since  $\bar{\eta}_{0}, \hat{Z}_0 $ and $\hat{W}^k_j$ 
are independent random variables,
using \eqref{eq:7.2} yields
\begin{align}
\nonumber
& \mathbb{E} \left(
\exp \left( 
p \mu_n \Delta
+ \sum_{k=1}^m
 p  \left\langle \hat{Z}_n  ,\bar{ \sigma }^k \left(\bar{\eta}_n  , \hat{Z}_n \right) \right\rangle   \right) \sqrt{\Delta} \hat{W}^k_{n+1} 
 \diagup  \mathfrak{G}^N_{n}
\right) 
\\
\nonumber
& =  
\exp \left( p \mu_n \Delta \right) 
\prod_{k=1}^m
\mathbb{E} \left( 
\exp \left( 
 p   \left\langle \hat{Z}_n  ,\bar{ \sigma }^k \left(\bar{\eta}_n  , \hat{Z}_n \right) \right\rangle  \sqrt{\Delta} \hat{W}^k_{n+1}  
 \right) 
 \diagup  \mathfrak{G}^N_{n}
\right) 
\\
\label{eq:5.2}
& \leq   
\exp \left(  p \mu_n \Delta  \right)
\prod_{k=1}^m
\exp \left(  \gamma \, p^2   \left\langle \hat{Z}_n  ,\bar{ \sigma }^k \left(\bar{\eta}_n  , \hat{Z}_n \right) \right\rangle ^2 \Delta  
\right) ,
\end{align}
where $\mathfrak{G}_{n}$ stands for the $\sigma$-algebra generated by
$\bar{\eta}_{0}, \hat{Z}_0 , \hat{W}^1_{\ell}, \ldots, \hat{W}^m_{\ell} $ with $\ell =1, \ldots, n$.
By Hypothesis \ref{hyp:Basica},
combining \eqref{eq:2.5}, \eqref{eq:StabCond} with \eqref{eq:5.2} we obtain
\begin{align*}
 \mathbb{E} \left( \left( \bar{\eta}_{n+1} \right)^p  \diagup  \mathfrak{G}_{n} \right) 
& \leq 
\left( \bar{\eta}_{n}  \right)^p \exp \left(  p \mu_n \Delta  \right)
\prod_{k=1}^m
\exp \left(  \gamma \, p^2   \left\langle \hat{Z}_n  ,\bar{ \sigma }^k \left(\bar{\eta}_n  , \hat{Z}_n \right) \right\rangle ^2 \Delta  
\right)
\\
& \leq 
\left( \bar{\eta}_{n}  \right)^p \exp \left(  p \, \Delta \left( - \lambda + m \,  \gamma \, p \, K^2  \right) \right) ,
\end{align*}
where $K$ is the constant appearing in Hypothesis \ref{hyp:Basica}. Thus
\begin{align}
\nonumber
 \mathbb{E} \left( \left( \bar{\eta}_{n+1} \right)^p  \diagup  \mathfrak{G}_{n} \right) 
& \leq 
\left( \bar{\eta}_{n}  \right)^p \exp \left(  p \, \Delta \left( - \lambda + \epsilon  \right) \right)
\exp \left(  p \, \Delta \left( - \epsilon + m \,  \gamma \, p \, K^2  \right) \right)
\\
\label{eq:5.4}
& \leq 
\left( \bar{\eta}_{n}  \right)^p \exp \left(  p \, \Delta \left( - \lambda + \epsilon  \right) \right)
\end{align}
whenever 
$
p < \epsilon / \left( m \,  \gamma  \, K^2 \right)
$.
Iterating \eqref{eq:5.4} we deduce \eqref{eq:5.3}.
\endproof

\subsection{Proof of Theorem \ref{th:NoStable}}
We return to the proof of Theorem \ref{th:Stable}.
By \eqref{eq:2.7}, 
it follows from \eqref{eq:2p.1} that
\begin{equation}
 \label{eq:2p.5}
 \log \left(\bar{\eta}_{n+1}\right)
\geq
 \log \left(\bar{\eta}_{0}\right) +  \theta  A_n  + S_n ,
\end{equation}
where 
$
A_n
=
\Delta \sum_{j=0}^n \sum_{k=1}^m \left\langle \hat{Z}_j , \bar{ \sigma }^k \left( \bar{\eta}_j ,  \hat{Z}_j \right) \right\rangle^2 
$.
Since $S_n$ is a square integrable martingale with quadratic variation process $A_n$,
$S_n$ converges a.s. on $\left\{ A_\infty < \infty \right\}$ to a finite random variable $\eta$
(see, e.g., Section 2.6.1 of \cite{DacunhaCastelle1983}).
Therefore  
$$
\liminf_{n \rightarrow \infty} \log \left(\bar{\eta}_{n}\right)
\geq
\log \left(\bar{\eta}_{0}\right) +  \theta  A_\infty  + \eta
\hspace{1cm}
\text{a.s. on } \left\{ A_\infty < \infty \right\} ,
$$
and so 
$
\liminf_{n \rightarrow \infty} \bar{\eta}_{n}  >  0
$
a.s. on the event  $\left\{ A_\infty < \infty \right\}$.

Applying the strong law of large numbers for martingales 
(see, e.g., Section 2.6.1 of \cite{DacunhaCastelle1983})
we  obtain that
$
S_n / A_n \longrightarrow_{n \rightarrow \infty} 0 
$ 
a.s. on  
$
\left\{ A_\infty = \infty \right\}
$.
Hence \eqref{eq:2p.5} yields
$
\liminf_{n \rightarrow \infty}  \log \left(\bar{\eta}_{n}\right) / A_n  \geq \theta
$
a.s. on 
$ \left\{ A_\infty = \infty \right\} $,
which implies 
$
\liminf_{n \rightarrow \infty} \bar{\eta}_{n} = + \infty
$
$a.s.$ on 
$
\left\{ A_\infty = \infty \right\}
$.
\endproof

\subsection{Proof of Theorem \ref{th:GeneralWeakConvergence}}
We abbreviate $\ell \left( N \right)$ to $\ell$.
For all $ t \in \left[ 0, T \right]$ and $x \in \mathbb{R}^d$,
we define $u\left( t,x\right) =\mathbb{E} \varphi \left( X_{T-t}^{x}\right) $.
From
$$
\mathbb{E} f \left(  \bar{Y}^N_{\ell}  \right) 
=
\mathbb{E}u\left( 0, \bar{Y}^N_0 \right) 
+ \sum_{n =0}^{\ell  -1}  \mathbb{E} \left( u \left( \tau^N_{n+1},  \bar{Y}^N_{n+1} \right)  
- u \left( \tau^N_{n}, \bar{Y}^N_{n} \right)
\right) 
$$
we have
$
\mathbb{E} f \left(  \bar{Y}^N_{ \ell  } \right) 
=
\mathbb{E}u\left( 0, \bar{Y}^N_0 \right) 
+
\sum_{n =0}^{ \ell -1}  \mathbb{E}  \left(  H_{1}^{n}  +  H_{2}^{n} \right) 
$,
where 
$$
H_{1}^{n} = 
u \left( \tau^N_{n+1},  \bar{Y}^N_{n+1}  \right) - u \left( \tau^N_{n+1},\bar{Y}^N_{n} \right) 
+ u \left( \tau^N_{n+1}, \bar{Y}^N_{n} \right)  -  u \left( \tau^N_{n+1}, E^n_{\tau^N_{n+1}} \right) 
$$
and 
$
H_{2}^{n} =  u \left( \tau^N_{n+1}, E^n_{\tau^N_{n+1}}  \right) - u \left( \tau^N_{n}, \bar{Y}^N_n \right) 
$.
Since $\mathbb{E} f \left( X_{T}\right) =  \mathbb{E}u\left( 0, X_0 \right)$, 
$$
\left\vert 
\mathbb{E} f \left( X_{T}\right) - \mathbb{E} f \left( \bar{Y}^N_{ \ell}  \right) \right\vert
 \leq 
\left\vert \mathbb{E}u\left( 0, \bar{Y}^N_{0}  \right) - \mathbb{E}u\left( 0, X_0 \right)   \right\vert 
+
\sum_{n =0}^{\ell - 1} \left( 
\left\vert \mathbb{E} \left( H_{1}^{n}\right) \right\vert + \left\vert \mathbb{E}\left( H_{2}^{n}\right) \right\vert
\right) .
$$
Combining $ u\left( 0, \cdot \right) \in C_{p}^{4} \left(  \mathbb{R}^{d},\mathbb{R}\right)  $ 
with  Condition (i) yields 
\begin{equation}
\label{eq:7.8}
 \left\vert 
\mathbb{E} f \left( X_{T}\right) - \mathbb{E} f \left(  \bar{Y}^N_{\ell } \right) \right\vert
 \leq  K \left(  1+\mathbb{E}\left\Vert X_{0} \right\Vert ^{q}\right)  \frac{T}{N}
+
\sum_{n =0}^{\ell-1} \left(  \left\vert \mathbb{E} \left( H_{1}^{n}\right) \right\vert  
+ \left\vert \mathbb{E}\left( H_{2}^{n}\right) \right\vert \right).
\end{equation}

By $u \left( \tau^N_{n+1}, \cdot \right) \in  C_{p}^{4} \left(  \mathbb{R}^{d},\mathbb{R}\right) $,
applying Taylor's formula we find
\begin{align*}
 H_{1}^{n} 
 & =
 \sum_{ k =1}^{3}\frac{1}{ k !}\sum_{\vec{p}\in \mathcal{P}_{k}}\partial _{x}^{\vec{p}}
 u \left( \tau^N_{n+1}, \bar{Y}^N_{n} \right)
\left( F_{\vec{p}}\left(  \bar{Y}^N_{n+1} -  \bar{Y}^N_{n} \right) 
-
F_{\vec{p}}\left(
E^n_{\tau^N_{n+1}}  -  \bar{Y}^N_{n} \right) 
\right) 
\\
& \quad 
+ R_{n}\left( \bar{Y}^N_{n+1} \right) +R_{n}\left( E^n_{T_{n+1}}  \right) ,
\end{align*}
where
$$
R_{n}\left( x\right) 
=
\frac{1}{4!}\sum_{\vec{p}\in P_{4}}\partial _{x}^{\vec{p}}u \left( \tau^N_{n+1},  \bar{Y}^N_{n} +\Xi_{\vec{p},n}\left( x\right)
\left( x-  \bar{Y}^N_{n} \right) \right) 
F_{\vec{p}}\left( x-  \bar{Y}^N_{n} \right) ,
$$
with $\Xi _{\vec{p},n} \left( x\right)$ diagonal $\mathbb{R}^{d \times d}$-matrix whose components belong to
$\left[ 0,1\right] $.
Using
$b, \sigma^k \in \mathcal{C}_{P}^{ 0 } \left( \mathbb{R}^d,\mathbb{R}\right)  $
and 
$ u  \in  \mathcal{C}_{P}^{ 4 } \left( \left[ 0, T \right] \times \mathbb{R}^d,\mathbb{R}\right)$
we obtain 
$$
\mathbb{E} \left( R_{n}\left( E^n_{\tau^N_{n+1}}  \right)  \diagup  \mathfrak{G}^N_{n}  \right)
\leq
K \left( 1+\left\Vert \bar{Y}^N_n \right\Vert ^{q}\right) \left( \tau^N_{n+1} - \tau^N_n \right)^{2} 
\hspace{0.5cm} \forall n \leq \ell-1 \text{ and } N \in \mathbb{N}.
$$
Moreover, $ u  \in  \mathcal{C}_{P}^{ 4 } \left( \left[ 0, T \right] \times \mathbb{R}^d,\mathbb{R}\right)$
and Condition (iii) lead to 
$$
\mathbb{E} \left( R_{n}\left( \bar{Y}^N_{n+1} \right)  \diagup  \mathfrak{G}^N_{n} \right)
\leq
K \left( 1+\left\Vert \bar{Y}^N_n \right\Vert ^{q}\right) \left(  \tau^N_{n+1} - \tau^N_n \right) ^{2} 
\hspace{0.8cm} \forall n \leq \ell-1 \text{ and } N \in \mathbb{N}.
$$
Therefore
\begin{align*}
& \left\vert \mathbb{E} \left( H_{1}^{n} \diagup  \mathfrak{G}^N_{n} \right) \right\vert 
\leq
 K \left( 1+\left\Vert \bar{Y}^N_n \right\Vert ^{q}\right) \left(  \tau^N_{n+1} - \tau^N_n \right) ^{2}
\\
&   \hspace{1cm} +
 \sum_{ k =1}^{3} \sum_{\vec{p}\in \mathcal{P}_{k}}
 \frac{\left\vert  \partial _{x}^{\vec{p}} u \left( \tau^N_{n+1}, \bar{Y}^N_{n} \right) \right\vert }{ k !}
\left\vert 
\mathbb{E} 
\left( F_{\vec{p}}\left(  \bar{Y}^N_{n+1} -  \bar{Y}^N_{n} \right) 
-
F_{\vec{p}}\left(
E^n_{\tau^N_{n+1}}  -  \bar{Y}^N_{n} \right) 
\diagup  \mathfrak{G}^N_{n}  \right) 
\right\vert ,
\end{align*}  
and so Condition (iv) gives
\begin{equation}
\label{eq:7.7}
 \left\vert \mathbb{E} \left( H_{1}^{n} \diagup  \mathfrak{G}^N_{n} \right) \right\vert 
 \leq
K \left( 1+\left\Vert \bar{Y}^N_n \right\Vert ^{q}\right)  \left(  \tau^N_{n+1} - \tau^N_n \right) T/N  
\hspace{0.8cm} \forall n \leq \ell-1 \text{ and } N \in \mathbb{N}.
\end{equation}

Using  It\^o's formula we obtain that $u$ 
is the unique $\mathcal{C}^{1,2 } \left( \left[ 0, T \right] \times \mathbb{R}^d,\mathbb{R}\right)$
solution of  \eqref{eq:Kolmogorov}
(see, e.g., proof of Theorem 5.7.6 of \cite{KaratzasShreve1998} 
or proof of Theorem 7.14 of  \cite{GrahamTalay2013}).
Hence,
using It\^o's formula  yields
\begin{align*}
 \mathbb{E} \left( H_{2}^{n}  \diagup  \mathfrak{F}_{\tau^N_n} \right)
& =
\int_{\tau^N_n}^{\tau^N_{n+1}} \mathbb{E}  \left(
\frac{\partial}{dt} u \left( t,  E^n_{t} \right) + \mathcal{L}_{\bar{Y}^N_{n}}\left( u \right) \left(t, E^n_t \right)
 \diagup  \mathfrak{F}_{\tau^N_n} \right) dt 
\\
& = 
\int_{\tau^N_n}^{\tau^N_{n+1}} \mathbb{E}  \left(
-  \mathcal{L}\left( u \right) \left(t, E^n_t \right) + \mathcal{L}_{\bar{Y}^N_{n}}\left( u \right) \left(t, E^n_t \right)
 \diagup  \mathfrak{F}_{\tau^N_n} \right) dt ,
\end{align*}
where 
$
\mathcal{L}_y=\sum_{k=1}^{d}b_{k} \left( y \right) \partial _{x}^{k}+
\frac{1}{2} \sum_{k, \ell =1}^{d}
\left( \sum_{j =1}^{m} \sigma_k^j \left( y \right)  \sigma_{\ell}^j \left( x \right) 
\right)
\partial_{x}^{k, \ell}
$.
Since
$u$ belongs to $\mathcal{C}^{1,4 } \left(  \left[ 0, T \right] \times \mathbb{R}^d,\mathbb{R}\right)$,
combining It\^o's formula with \eqref{eq:KolmogorovG}
we deduce that
\begin{align*}
 \mathbb{E} \left( H_{2}^{n}  \diagup  \mathfrak{F}_{\tau^N_n} \right)
& =
\int_{\tau^N_n}^{\tau^N_{n+1}} \int_{\tau^N_n}^{t}
 \mathbb{E}  \left(
 \mathcal{L}^2\left( u \right) \left(s, E^n_s \right) - 
 \mathcal{L}_{\bar{Y}^N_{n}}\left(   \mathcal{L}\left( u \right) \right) \left(s, E^n_s \right)
 \diagup  \mathfrak{F}_{\tau^N_n} \right) ds dt
\\
& \quad 
+ \int_{\tau^N_n}^{\tau^N_{n+1}} \int_{\tau^N_n}^{t}
 \mathbb{E}  \left(
- \mathcal{L}_{\bar{Y}^N_{n}}\left(   \mathcal{L}\left( u \right) \right) \left(s, E^n_s \right)
+ \mathcal{L}_{\bar{Y}^N_{n}}^2\left( u \right) \left(s, E^n_s \right)
 \diagup  \mathfrak{F}_{\tau^N_n} \right) ds dt ,
\end{align*}
and so employing 
$ u  \in  \mathcal{C}_{P}^{ 4 } \left( \left[ 0, T \right] \times \mathbb{R}^d,\mathbb{R}\right)$
and 
$b, \sigma^k \in \mathcal{C}_{P}^{ 0 } \left( \mathbb{R}^d,\mathbb{R}\right)  $
we get 
\begin{equation}
\label{eq:7.15}
 \left\vert 
\mathbb{E} \left( H_{2}^{n}  \diagup  \mathfrak{F}_{\tau^N_n} \right)
\right\vert 
\leq
K \left( 1+\left\Vert \bar{Y}^N_n \right\Vert ^{q}\right)  \left (  \tau^N_{n+1} - \tau^N_n \right) ^{2}  
\hspace{1cm} \forall n \leq \ell-1 \text{ and } N \in \mathbb{N}.
\end{equation}

From Condition (ii),  \eqref{eq:7.7} and \eqref{eq:7.15}
it follows that
$$
 \left\vert \mathbb{E} \left( H_{1}^{n}\right) \right\vert  
+ \left\vert \mathbb{E}\left( H_{2}^{n}\right) \right\vert
\leq
K \left( 1+\mathbb{E} \left( \left\Vert \bar{Y}^N_0 \right\Vert ^{q}\right) \right) \left(  \tau^N_{n+1} - \tau^N_n \right) T/N 
\hspace{1cm}
\forall N \in \mathbb{N}.
$$
Hence,
Condition (i) and \eqref{eq:7.8} lead to  \eqref{eq:7.16}.
\endproof

\subsection{Proof of Theorem \ref{th:RateConvergence}} 
\label{subsec:Proofth:RateConvergence}

\

\begin{notation}
We abbreviate 
$ \bar{ b }_n := \bar{ b } \left( \bar{\eta}_n , \hat{U}_n \right)$
and
$  \bar{ \sigma }_n^k :=  \bar{ \sigma }^k \left(\bar{\eta}_n  , \hat{U}_n \right)$.
Let $\mathfrak{G}_{n}$ denote the $\sigma$-algebra generated by
$ \bar{X}_{0}, \hat{W}^k_\ell$ and  $W^k_{T_{\ell}}$,
where $k=1,\ldots,m$ and $\ell =1, \ldots, n$.
We use the same symbol 
$\left( \mathcal{O}_n \right)_{n \geq 1}$
for different stochastic processes such that 
\begin{equation}
\label{eq:7.32}
 \mathbb{E} \left( \left\Vert \mathcal{O}_{n+1} \right\Vert ^p \diagup \mathfrak{G}_{n} \right)
\leq 
 K_p \left(T \right) \left( 1 +  \left\Vert \bar{X}_n \right\Vert^{p q }  \right) 
\end{equation}
for all $n=0,\ldots, N-1$ and $p  \in \mathbb{N}$,
where $K_p$ stands for a non-negative increasing function that is independent of $N$.
\end{notation}

First, Lemma  \ref{lem:ExpansionSolucionN} below
provides an asymptotic expansion of $\bar{\rho}_{n+1} \bar{U}_{n+1}$
as $\Delta \rightarrow 0+$.
To this end,
using \eqref{eq:7.2} and \eqref{eq:7.1} we next obtain
an upper bound for the conditional moments of $\bar{\rho}_{n+1}$.

\begin{lemma}
\label{lem:AcotacionInicial} 
Assume \eqref{eq:7.2} and Hypothesis \ref{hyp:Acotacion}.
Suppose that $\bar{X}_0$ has finite moments of all orders.
Then, for any $q > 0 $ there exits  $K_q \geq 0$ such that
for any $\Delta > 0$,
\begin{equation*}
 \mathbb{E} \left(
\exp \left( 
q \mu_n \Delta
+ \sum_{k=1}^m
\left( 
 q  \left\langle \hat{U}_n  ,\bar{ \sigma }^k_n \right\rangle \sqrt{\Delta} \hat{W}^k_{n+1}  
 \right)  \right) 
 \diagup  \mathfrak{G}_{n}
\right) 
\leq
\exp \left( \Delta K_q \right) 
\qquad
\forall n \geq 0 .
\end{equation*}
\end{lemma}

\begin{proof}
Let $\alpha \geq 0$ be as in Scheme \ref{scheme:EulerStableG}.
For 
 any $\beta \geq 1$
there exits $K_\beta \geq 0$ satisfying 
\begin{equation}
 \label{eq:7.9}
 \langle x , b \left( x \right) \rangle
 +
\beta \sum_{k=1}^{m}  \left\Vert \sigma^k \left( x \right)  \right\Vert^2
\leq 
K_\beta \left( \alpha + \left\Vert x \right\Vert^2 \right)
\qquad \quad
\forall x \in \mathbb{R}^d .
\end{equation}
Indeed, if $\alpha > 0$,
then  \eqref{eq:7.9} follows directly from \eqref{eq:7.1}.
In case $\alpha = 0$,
using the locally Lipschitz property of $b$ and $\sigma^k$ we deduce that
for all $\left\Vert  x  \right\Vert \leq 1$ we have
\begin{equation}
\label{eq:7.41}
\left\Vert b \left( x \right)  \right\Vert
=
\left\Vert b \left( x \right) -  b \left( 0 \right)  \right\Vert
\leq K  \left\Vert  x  \right\Vert
\text{ and } 
\left\Vert \sigma^k \left( x \right)  \right\Vert
=
\left\Vert \sigma^k \left( x \right)  -  \sigma^k \left( 0 \right)  \right\Vert
\leq K  \left\Vert  x  \right\Vert ,
\end{equation}
and so  \eqref{eq:7.1} yields
$
 \langle x , b \left( x \right) \rangle
 +
\beta \sum_{k=1}^{m}  \left\Vert \sigma^k \left( x \right)  \right\Vert^2
\leq 
K \left\Vert x \right\Vert^2 
$
for all $x \in \mathbb{R}^d$.
This gives  \eqref{eq:7.9}.

Proceeding as in the proof of  \eqref{eq:5.2},
from \eqref{eq:7.2} we obtain
$$
\mathbb{E} \left(
\exp \left( 
q \mu_n \Delta
+ \sum_{k=1}^m
 q  \left\langle \hat{U}_n  ,\bar{ \sigma }^k_n \right\rangle \sqrt{\Delta} \hat{W}^k_{n+1}  \right) 
 \diagup  \mathfrak{G}_{n}
\right) 
\leq   
e^{   q \mu_n \Delta }
\prod_{k=1}^m
e^{  \gamma \, q^2  \left\langle \hat{U}_n  ,\bar{ \sigma }^k_n \right\rangle^2 \Delta } .
$$
Then,
applying $\left\Vert  \hat{U}_n  \right\Vert \leq 1$ we get
\begin{align*}
& \mathbb{E} \left(
\exp \left( 
q \mu_n \Delta
+ \sum_{k=1}^m
\left( 
 q  \left\langle \hat{U}_n  ,\bar{ \sigma }^k_n \right\rangle \sqrt{\Delta} \hat{W}^k_{n+1}  
 \right)  \right) 
 \diagup  \mathfrak{G}_{n}
\right) 
\\
&
\leq
\exp \left( \Delta
 \left( q \left\langle \hat{U}_n , \bar{ b }_n \right\rangle 
+\left( \frac{q}{2} + \left\vert  \gamma q^2 - q \right\vert \right)   \sum_{k=1}^m \left\Vert \bar{\sigma}^k_n \right\Vert^2
\right) \right) .
\end{align*}
According to \eqref{eq:7.9} we have
$$
q \left\langle \bar{X}_n , b \left( \bar{X}_n \right) \right\rangle
+
\left( \frac{q}{2} + \left\vert  \gamma q^2 - q \right\vert \right)   \sum_{k=1}^m \left\Vert \sigma^k  \left( \bar{X}_n \right) \right\Vert^2
\leq 
K_q \, \left( \bar{\eta}_{n} \right)^2
\hspace{1cm} \forall n \geq 0 ,
$$
with  $K_q \geq 0$.
Therefore, \eqref{eq:3.5} leads to
$$
\exp \left( \Delta
 \left( q \left\langle \hat{U}_n , \bar{ b }_n \right\rangle 
+\left( \frac{q}{2} + \left\vert  \gamma q^2 - q \right\vert \right)   \sum_{k=1}^m \left\Vert \bar{\sigma}^k_n \right\Vert^2
\right) \right)
\leq 
\exp \left( \Delta K_q \right)
\hspace{0.5cm} \forall n \geq 0 .
$$
\end{proof}

To prove Lemma  \ref{lem:ExpansionSolucionN},
we multiply \eqref{eq:7.13} by 
the local asymptotic expansion of
$
\exp  \left( 
 \mu_n \Delta
+
\sum_{k=1}^m \left\langle \hat{U}_n  ,\bar{ \sigma }_n^k  \right\rangle \sqrt{\Delta} \hat{W}^k_{n+1} 
\right)
$.
This, together Lemma \ref{lem:AcotacionInicial}, yields  \eqref{eq:7.5}.

\begin{lemma}
\label{lem:ExpansionSolucionN}
Assume Hypothesis \ref{hyp:Acotacion}, together with \eqref{eq:7.2} and 
$b, \sigma^k \in \mathcal{C}_{P}^{ 0 } \left( \mathbb{R}^d,\mathbb{R}\right)  $.
Let
$
\mathbb{E} \left( \left\Vert \bar{X}_0 \right\Vert ^q \right) < \infty
$
for all $q \in \mathbb{N}$,
and let $ \left( \bar{U}_{n}, \bar{V}_{n}  \right) \neq 0$
for $n=1,\ldots, N$.
Then 
\begin{align}
\label{eq:7.5}
\bar{\rho}_{n+1} \bar{U}_{n+1}
& =
\bar{\eta}_{n} \hat{U}_n
+ \Delta \bar{\eta}_{n} \bar{ b }_n
+  \sqrt{\Delta} \bar{\eta}_{n} \sum_{k=1}^m   \bar{ \sigma }_n^k    \hat{W}^k_{n+1}
\\ \nonumber
& \quad
+ \Delta \bar{\eta}_{n}  \sum_{k=1}^m \left\langle \hat{U}_n  ,\bar{ \sigma }_n^k  \right\rangle 
\left( \bar{ \sigma }_n^k - \frac{1}{2} \left\langle \hat{U}_n  ,\bar{ \sigma }_n^k  \right\rangle  \hat{U}_n 
\right)
\left( \left( \hat{W}^k_{n+1} \right)^2 - 1 \right)
\\ \nonumber
& \quad
+ \Delta \bar{\eta}_{n}  \sum_{k \neq j} \left\langle \hat{U}_n  ,\bar{ \sigma }_n^k  \right\rangle 
\left( \bar{ \sigma }_n^j - \frac{1}{2} \left\langle \hat{U}_n  ,\bar{ \sigma }_n^j  \right\rangle  \hat{U}_n 
\right) \hat{W}^k_{n+1} \hat{W}^j_{n+1}
\\ \nonumber
& \quad
+ \Delta^{3/2} \bar{\eta}_{n} \widetilde{\Gamma}_{n+1} + \Delta^2 \mathcal{O}_{n+1} ,
\end{align}
where $n=0,\ldots, N-1$,  $\mathcal{O}_{n+1}$ satisfies \eqref{eq:7.32} and  
\begin{align*}
\widetilde{\Gamma}_{n+1} 
 & =
 \left( \mu_n  + \frac{1}{2} 
 \left(
\sum_{k=1}^m \left\langle \hat{U}_n  ,\bar{ \sigma }_n^k  \right\rangle \hat{W}^k_{n+1} 
\right)^2
 \right)
\sum_{j=1}^m  \left( \bar{ \sigma }_n^j
- \left\langle  \hat{U}_n ,  \bar{ \sigma }_n^j  \right\rangle  \hat{U}_n \right) \hat{W}^j_{n+1} 
\\
& \quad
+ \sum_{k=1}^m  \left\langle \hat{U}_n  ,\bar{ \sigma }_n^k  \right\rangle 
\left(
\mu_n \hat{U}_n + \bar{ b }_n  - \left\langle  \hat{U}_n , \bar{ b }_n   \right\rangle  \hat{U}_n
+ \Psi  \left( \bar{\eta}_n , \hat{U}_n \right) 
\right)  \hat{W}^k_{n+1}
\\
& \quad
+
\left(
\sum_{k=1}^m \left\langle \hat{U}_n  ,\bar{ \sigma }_n^k  \right\rangle \hat{W}^k_{n+1} 
\right)^3 \hat{U}_n / 6 .
\end{align*}

\end{lemma}

\begin{proof}
Using the Taylor expansion of $x \mapsto \exp \left( x \right)$ we deduce that
\begin{equation}
\label{eq:7.3}
 \begin{aligned}
\exp  \left( 
 \mu_n \Delta
+
\sum_{k=1}^m \left\langle \hat{U}_n  ,\bar{ \sigma }_n^k  \right\rangle \sqrt{\Delta} \hat{W}^k_{n+1} 
\right)
& =
1
+
\mu_n \Delta
+
 \sum_{k=1}^m \left\langle \hat{U}_n  ,\bar{ \sigma }_n^k  \right\rangle \sqrt{\Delta} \hat{W}^k_{n+1}
\\
& \hspace{-2cm}
+ \frac{\Delta}{2} \left( \sum_{k=1}^m \left\langle \hat{U}_n  ,\bar{ \sigma }_n^k  \right\rangle  \hat{W}^k_{n+1} \right)^2
 + \Delta^{\frac{3}{2}} \bar{\Gamma}_{n+1} + \Delta^2 R_{n+1} ,
\end{aligned}
\end{equation}
where
$ 
\bar{\Gamma}_{n+1}
=
\mu_n \, \sum_{k=1}^m \left\langle \hat{U}_n  ,\bar{ \sigma }_n^k  \right\rangle \hat{W}^k_{n+1} 
+
\left(
\sum_{k=1}^m \left\langle \hat{U}_n  ,\bar{ \sigma }_n^k  \right\rangle \hat{W}^k_{n+1} 
\right)^3 / 6
$
and
\begin{align*}
 R_{n+1} & =  \left( \mu_n \right)^2 
 + \frac{\Delta }{6} \left( \mu_n \right)^3 
 + \frac{ \sqrt{\Delta} }{2} \left( \mu_n \right)^2 
 \sum_{k=1}^m \left\langle \hat{U}_n  ,\bar{ \sigma }_n^k  \right\rangle \hat{W}^k_{n+1} 
 \\
 & \quad
 + \frac{1 }{2}  \mu_n 
  \left(\sum_{k=1}^m \left\langle \hat{U}_n  ,\bar{ \sigma }_n^k  \right\rangle \hat{W}^k_{n+1} \right)^2
+ \frac{\exp \left( \xi_n \right)}{24}
 \left( 
 \mu_n \sqrt{\Delta}
+
\sum_{k=1}^m \left\langle \hat{U}_n  ,\bar{ \sigma }_n^k  \right\rangle  \hat{W}^k_{n+1} 
\right)^4
\end{align*}
with 
$\xi_n$  between $0$ and 
$
 \mu_n \Delta
+
\sum_{k=1}^m \left\langle \hat{U}_n  ,\bar{ \sigma }_n^k  \right\rangle \sqrt{\Delta} \hat{W}^k_{n+1} 
$.
Multiplying \eqref{eq:7.13} by \eqref{eq:7.3}
we  obtain,  after a careful computation, that \eqref{eq:7.5} holds 
with $\mathcal{O}_{n+1}$ replaced by
\begin{align*}
\widetilde{R}_{n+1}
& := 
R_{n+1} \bar{\eta}_n \hat{U}_n 
+ \Delta R_{n+1} \bar{\eta}_n  \left( 
 \bar{ b }_n  
  - \left\langle  \hat{U}_n , \bar{ b }_n  \right\rangle  \hat{U}_n
  + \Psi  \left( \bar{\eta}_n , \hat{U}_n \right)
 \right) 
\\
&  \quad
+ \Delta^{1/2} R_{n+1} \bar{\eta}_n  \sum_{k=1}^m  \left( \bar{ \sigma}^k_n  
- \left\langle  \hat{U}_n ,  \bar{ \sigma}^k_n \right\rangle  \hat{U}_n \right) \hat{W}^k_{n+1}
\\
&  \quad 
+ 
 \bar{\eta}_n \bar{\Gamma}_{n+1}  \sum_{k=1}^m  \left( \bar{ \sigma}^k_n  
- \left\langle  \hat{U}_n ,  \bar{ \sigma}^k_n \right\rangle  \hat{U}_n \right) \hat{W}^k_{n+1}  
+  
  \mu_n  \bar{\eta}_n \bar{ b }_n  
  - \mu_n \left\langle  \hat{U}_n , \bar{ b }_n  \right\rangle    \bar{\eta}_n  \hat{U}_n
 \\
 & \quad
  +  \mu_n  \bar{\eta}_n \Psi  \left( \bar{\eta}_n , \hat{U}_n 
 \right) 
 + \Delta^{\frac{1}{2}}  \bar{\eta}_n \bar{\Gamma}_{n+1} \left( 
 \bar{ b }_n  
  - \left\langle  \hat{U}_n , \bar{ b }_n  \right\rangle  \hat{U}_n
  + \Psi  \left( \bar{\eta}_n , \hat{U}_n \right)
 \right) 
  \\
& \quad
+  \frac{1}{2} \bar{\eta}_n 
\left( \sum_{k=1}^m \left\langle \hat{U}_n  ,\bar{ \sigma }_n^k  \right\rangle  \hat{W}^k_{n+1} \right)^2
\left( 
 \bar{ b }_n  
  - \left\langle  \hat{U}_n , \bar{ b }_n  \right\rangle  \hat{U}_n
  + \Psi  \left( \bar{\eta}_n , \hat{U}_n \right)
 \right)  .
\end{align*}

Let $p \geq 1$. 
Combining $b, \sigma^k \in \mathcal{C}_{P}^{0}\left( \mathbb{R}^{d},\mathbb{R}\right) $
with \eqref{eq:7.41} we get
\begin{equation}
\label{eq:7.42}
\left\Vert \bar{ b }_n \right\Vert 
\leq 
 K  \left( 1 +  \left\Vert \bar{X}_n \right\Vert^{ q }  \right) 
 \text{ and }
 \left\Vert \bar{ \sigma }_n^k \right\Vert 
\leq 
 K  \left( 1 +  \left\Vert \bar{X}_n \right\Vert^{ q }  \right) 
 \hspace{1cm} \forall n \geq 0.
\end{equation}
Applying H\"{o}lder's inequality gives
\begin{equation}
\label{eq:7.43}
\left( \sum_{j=1}^M a_j \right)^q 
\leq
M^{q-1}  \sum_{j=1}^M \left( a_j \right)^q
\end{equation}
whenever $a_j \geq 0$.
Using $\left\Vert \hat{U}_n \right\Vert  \leq 1$, \eqref{eq:7.42} and \eqref{eq:7.43} 
we obtain
$$
\left( 
 \mu_n \sqrt{\Delta}
+
\sum_{k=1}^m \left\langle \hat{U}_n  ,\bar{ \sigma }_n^k  \right\rangle  \hat{W}^k_{n+1} 
\right)^4
\leq
K \left( T \right) \left( 1 +  \left\Vert \bar{X}_n \right\Vert^q \right) 
\left( 1 + \sum_{k=1}^m \left(\hat{W}^k_{n+1} \right)^4 \right) 
$$
for all $n \geq 0$.
Therefore, \eqref{eq:7.2} and \eqref{eq:7.43} lead to
$$
\mathbb{E} \left( 
\left\Vert 
\bar{\eta}_n \hat{U}_n
\left( 
 \mu_n \sqrt{\Delta}
+
\sum_{k=1}^m \left\langle \hat{U}_n  ,\bar{ \sigma }_n^k  \right\rangle  \hat{W}^k_{n+1} 
\right)^4
 \right\Vert^{2p}
 \diagup  \mathfrak{G}_{n} \right)
 \leq
 K_p \left(T \right) \left( 1 +  \left\Vert \bar{X}_n \right\Vert^{p q }  \right) 
$$
for all $n \geq 0$.
On the other hand, 
from Lemma \ref{lem:AcotacionInicial}  it follows that
$$
\mathbb{E} \left( e^{ 2 p  \, \xi_n } \diagup  \mathfrak{G}_{n} \right)
\leq 
1 + \mathbb{E} \left(  
e^{  2 p  \left(
\mu_n \Delta + \sum_{k=1}^m \left\langle \hat{U}_n  ,\bar{ \sigma }_n^k  \right\rangle \sqrt{\Delta} \hat{W}^k_{n+1} 
\right) }
\diagup  \mathfrak{G}_{n} \right)
\leq 
1 + \exp \left( K_p \Delta \right) ,
$$
with $K_{ p } > 0$.
Applying the conditional Holder's inequality we deduce that for all $n \geq 0$,
$$
\mathbb{E} \left( 
\left\Vert 
e^{ \xi_n }
\left( 
 \mu_n \sqrt{\Delta}
+
\sum_{k=1}^m \left\langle \hat{U}_n  ,\bar{ \sigma }_n^k  \right\rangle  \hat{W}^k_{n+1} 
\right)^4
\bar{\eta}_n \hat{U}_n
 \right\Vert^{p}
 \diagup  \mathfrak{G}_{n} \right)
 \leq
 K_p \left(T \right) \left( 1 +  \left\Vert \bar{X}_n \right\Vert^{p q }  \right) ,
$$
where $q$ is independent of $N$ and $p$.
In much the same way,
we handle the other terms of  $\widetilde{R}_{n+1}$ to obtain
\begin{equation}
\label{eq:7.45}
 \mathbb{E} \left( \left\Vert \widetilde{R}_{n+1} \right\Vert ^p \diagup \mathfrak{G}_{n} \right)
\leq 
 K_p \left(T \right) \left( 1 +  \left\Vert \bar{X}_n \right\Vert^{p q }  \right) 
 \hspace{1.5cm} \forall n \geq 0 .
\end{equation}
\end{proof}

Second,
Lemma \ref{lem:LocalAsymptotic} studies 
the effect  of normalizing $\left( \bar{U}_{n+1} , \bar{V}_{n+1} \right)$.
The following asymptotic expansion of $x / \left\Vert x \right\Vert $
plays a key role in this analysis.

\begin{lemma}
\label{lem:DesigualdaProyecccion}
For any  $x \in \mathbb{R}^d \setminus \left\{ 0 \right\}$,
\begin{equation*}
\left\Vert 
\frac{x}{ \left\Vert x \right\Vert } - x - \frac{1}{2} \left( 1-  \left\Vert x \right\Vert^2 \right) x
\right\Vert
\leq
\left( 1-  \left\Vert x \right\Vert^2 \right) ^2 .
\end{equation*}
\end{lemma}

\begin{proof}
Let $x  \neq 0$.
Since 
$
1 /  \left\Vert x \right\Vert - 1
=
\left( 1 -  \left\Vert x \right\Vert^2 \right) / \Bigl( \left\Vert x \right\Vert \left( 1 +  \left\Vert x \right\Vert \right)   \Bigr) 
$,
\begin{align*}
 \frac{1}{ \left\Vert x \right\Vert} - 1
- 
\frac{1}{2} \left( 1 - \left\Vert x \right\Vert^2 \right) 
& =
\left( 1 - \left\Vert x \right\Vert^2 \right) 
\left( \frac{1}{ \left\Vert x \right\Vert \left( 1 +  \left\Vert x \right\Vert \right) }
 - \frac{1}{2} \right) 
\\
& =
\left( 1 - \left\Vert x \right\Vert^2 \right) 
\frac{ 1 - \left\Vert x \right\Vert^2 + 1 -  \left\Vert x \right\Vert }
{2  \left\Vert x \right\Vert \left( 1 +  \left\Vert x \right\Vert \right) }
\\
& =
\frac{ \left( 1 - \left\Vert x \right\Vert^2 \right) ^2
}{ 
2 \left\Vert x \right\Vert
}
\left( 
\frac{ 1 }{  1 +  \left\Vert x \right\Vert  }
+
\frac{ 1 }{  \left( 1 +  \left\Vert x \right\Vert \right)^2 }
\right) .
\end{align*}
Then
$
\left\vert
\frac{1}{ \left\Vert x \right\Vert} - 1
- 
\frac{1}{2} \left( 1 - \left\Vert x \right\Vert^2 \right) 
\right\vert
\leq
\left( 1 - \left\Vert x \right\Vert^2 \right) ^2 / \left\Vert x \right\Vert  
$,
and the lemma follows.
\end{proof}

Next,
using Lemmata  \ref{lem:AcotacionInicial}, \ref{lem:ExpansionSolucionN} 
and \ref{lem:DesigualdaProyecccion} 
we characterize the local asymptotic behavior of $\bar{X}_{n+1}$
as $\Delta \rightarrow 0+$.

\begin{lemma}
\label{lem:LocalAsymptotic}
Under the assumptions of Lemma \ref{lem:ExpansionSolucionN},
\begin{equation}
 \label{eq:7.14}
 \begin{aligned}
\bar{X}_{n+1}
& =
\bar{X}_{n} + \Delta b \left( \bar{X}_{n} \right)
+  \sum_{k=1}^m   \sigma^k  \left( \bar{X}_{n} \right)  \sqrt{\Delta}  \hat{W}^k_{n+1}
\\ 
&
+ \Delta  \sum_{k=1}^m \left( 
\left\langle \hat{U}_n  ,\bar{ \sigma }_n^k  \right\rangle 
 \bar{\eta}_{n} \bar{ \sigma }_n^k - \frac{1}{2} \left\Vert \bar{ \sigma }_n^k  \right\Vert^2  \bar{\eta}_{n}\hat{U}_n 
\right)
\left( \left( \hat{W}^k_{n+1} \right)^2 - 1 \right)
\\  
& 
+ \Delta   \sum_{k \neq j} \left(
 \left\langle \hat{U}_n  ,\bar{ \sigma }_n^k  \right\rangle 
\bar{\eta}_{n} \bar{ \sigma }_n^j - \frac{1}{2} \left\langle \bar{ \sigma }_n^k  ,\bar{ \sigma }_n^j  \right\rangle  \bar{\eta}_{n}\hat{U}_n 
\right) \hat{W}^k_{n+1} \hat{W}^j_{n+1}
\\  
&
+ \Delta^{3/2}  \bar{\eta}_{n}  \Gamma_{n+1} + \Delta^2 \mathcal{O}_{n+1} ,
\end{aligned}
\end{equation}
where $n=0,\ldots, N-1$, 
$\mathcal{O}_{n+1}$ is as in \eqref{eq:7.32}, and 
\begin{align*}
 \Gamma_{n+1} 
 & = \widetilde{\Gamma}_{n+1} 
 - \frac{1}{2}   \sum_{k, i =1}^m   \left( 
\left\langle \hat{U}_n    ,   \bar{ \sigma }_n^k \right\rangle^2 - \left\Vert  \bar{ \sigma }_n^k \right\Vert^2
\right) \bar{ \sigma }_n^i \left( 1-  \left( \hat{W}^k_{n+1} \right) ^2 \right)  \hat{W}^i_{n+1}
\\
& \quad 
- \frac{1}{2}     \sum_{i=1}^m \sum_{k \neq j}^m  \left(
\left\langle
\bar{ \sigma }_n^k   , \bar{ \sigma }_n^j   
\right\rangle 
-
\left\langle  \hat{U}_n ,  \bar{ \sigma }_n^k  \right\rangle   \left\langle  \hat{U}_n ,  \bar{ \sigma }_n^j  \right\rangle
\right) \bar{ \sigma }_n^i \, \hat{W}^k_{n+1} \hat{W}^j_{n+1}     \hat{W}^i_{n+1}
\\
& \quad
- 
  \sum_{k=1}^{m}
\left\langle     
\bar{ b }_n ,  
 \bar{ \sigma }_n^k   
- \left\langle  \hat{U}_n ,  \bar{ \sigma }_n^k  \right\rangle  \hat{U}_n 
\right\rangle \hat{U}_n \hat{W}^k_{n+1}
\\
& \quad
-  \sum_{k,j=1}^{m} \left(
\left\langle  \hat{U}_n ,  \bar{ \sigma }_n^j  \right\rangle 
\left\langle   \bar{ \sigma }_n^j  ,  \bar{ \sigma }_n^k  \right\rangle
-
\left\langle  \hat{U}_n ,  \bar{ \sigma }_n^j  \right\rangle ^2 
\left\langle  \hat{U}_n ,  \bar{ \sigma }_n^k  \right\rangle 
\right)  \hat{U}_n  \hat{W}^k_{n+1} .
\end{align*}
\end{lemma}

\begin{proof}
Since 
$
\left\Vert \hat{U}_n \right\Vert^2 + \left( \alpha / \bar{\eta}_n \right)^2 = 1 
$,
expanding and collecting like terms we obtain
\begin{align}
\label{eq:7.44}
& \langle \bar{U}_{n+1} , \bar{U}_{n+1} \rangle + \left( \bar{V}_{n+1} \right)^2
\\ \nonumber
& =
1
+
\Delta \sum_{k=1}^m \left( 
\left\langle \hat{U}_n    ,   \bar{ \sigma }_n^k \right\rangle^2 - \left\Vert  \bar{ \sigma }_n^k \right\Vert^2
\right) \left( 1-  \left( \hat{W}^k_{n+1} \right) ^2 \right) 
\\ \nonumber
& \quad 
+ \Delta \sum_{k \neq j}^m \left(
\left\langle
\bar{ \sigma }_n^k   , \bar{ \sigma }_n^j   
\right\rangle 
-
\left\langle  \hat{U}_n ,  \bar{ \sigma }_n^k  \right\rangle   \left\langle  \hat{U}_n ,  \bar{ \sigma }_n^j  \right\rangle
\right) \hat{W}^k_{n+1} \hat{W}^j_{n+1}
\\ \nonumber
& \quad
+ 
2 \Delta^{3/2} \sum_{k=1}^{m}
\left\langle     
\bar{ b }_n ,  
 \bar{ \sigma }_n^k   
- \left\langle  \hat{U}_n ,  \bar{ \sigma }_n^k  \right\rangle  \hat{U}_n 
\right\rangle  \hat{W}^k_{n+1}
\\ \nonumber
& \quad
+ 
2 \Delta^{3/2} \sum_{k,j=1}^{m} \left(
\left\langle  \hat{U}_n ,  \bar{ \sigma }_n^j  \right\rangle 
\left\langle   \bar{ \sigma }_n^j  ,  \bar{ \sigma }_n^k  \right\rangle
-
\left\langle  \hat{U}_n ,  \bar{ \sigma }_n^j  \right\rangle ^2 
\left\langle  \hat{U}_n ,  \bar{ \sigma }_n^k  \right\rangle 
\right)  \hat{W}^k_{n+1}
+ \Delta^2 R_{n+1} ,
\end{align}
where
\begin{align*}
 R_{n+1} 
 & =
 \left(  \frac{ \alpha }{ \bar{\eta}_n } \right)^2 
 \left(
  \sum_{k=1}^m   \left( \frac{3}{2} \left\langle  \hat{U}_n ,\bar{  \sigma }_n^k  \right\rangle^2
- \frac{1}{2} \left\Vert \bar{ \sigma }_n^k  \right\Vert ^2 
 \right)  
 -   \left\langle  \hat{U}_n , \bar{ b }_n \right\rangle
 \right)^2
 \\
 & \quad +
 \left\Vert  
   \bar{ b }_n 
  - \left\langle  \hat{U}_n , \bar{ b }_n  \right\rangle  \hat{U}_n
  + \Psi  \left( \bar{\eta}_n , \hat{U}_n \right)
 \right\Vert^2 .
\end{align*}
From \eqref{eq:7.42}, 
$
\left\Vert \hat{U}_n \right\Vert \leq 1
$
and
$
\left( \alpha / \bar{\eta}_n \right)^2 \leq 1 
$
we have  
$$
\left\Vert R_{n+1}  \right\Vert 
\leq 
 K \left( 1 +  \left\Vert \bar{X}_n \right\Vert^{q }  \right) 
 \hspace{2cm} \forall n \geq 0 .
$$

As $ \left( \bar{U}_{n+1}, \bar{V}_{n+1}  \right) \neq 0$,
combining Lemma \ref{lem:DesigualdaProyecccion} with  \eqref{eq:7.44} we deduce that
$$
\frac{\bar{U}_{n+1} }{  \left\Vert  \left( \bar{U}_{n+1}, \bar{V}_{n+1}  \right) \right\Vert }
- 
\bar{U}_{n+1}
- 
\frac{1}{2} \left( 1 	- \left\Vert  \left( \bar{U}_{n+1}, \bar{V}_{n+1}  \right) \right\Vert ^2 \right) \bar{U}_{n+1} 
=
 \Delta^2 \mathcal{O}_{n+1} .
$$
To this end, we use an analysis similar to that in the proof of \eqref{eq:7.45}.
Now, Lemma \ref{lem:AcotacionInicial} yields 
$$
\bar{X}_{n+1} 
- 
\bar{\rho}_{n+1} \bar{U}_{n+1}
- 
\frac{1}{2} \left( 1 	- \left\Vert  \left( \bar{U}_{n+1}, \bar{V}_{n+1}  \right) \right\Vert ^2 \right) \bar{\rho}_{n+1} \bar{U}_{n+1} 
=
 \Delta^2 \mathcal{O}_{n+1} .
$$
Applying Lemma \ref{lem:ExpansionSolucionN} and \eqref{eq:7.44}
together with arguments similar to that in the proof of \eqref{eq:7.45},
we get \eqref{eq:7.14}.
\end{proof}

Third, we show that
the moments of $\bar{X}_n$ are uniformly bounded with respect to $\Delta$,
in the time interval $\left[ 0, T \right]$.
This is a challenge problem for numerical schemes solving SDEs with locally Lipschitz coefficients.
In this direction,
using  Lemma \ref{lem:AcotacionInicial} we now estimate crudely the growth of  $\left\Vert \bar{X}_n \right\Vert$.

\begin{lemma}
\label{lem:AcotacionMomentosCruda} 
Let Hypothesis \ref{hyp:Acotacion} hold, along with \eqref{eq:7.2}.
Suppose that  
$
\mathbb{E} \left( \left\Vert \bar{X}_0 \right\Vert ^q \right) < \infty
$
for all $q \in \mathbb{N}$.
Consider $p \in \mathbb{N}$. 
Then,
for any $n \geq 1$ we have
\begin{equation*}
\mathbb{E} \left(  \left\Vert  \bar{X}_{n} \right\Vert^{2 \ell} \right)
 \leq
 e^{ n  \Delta K_p } \mathbb{E} \left(  \left\Vert  \bar{X}_{0} \right\Vert^{2  \ell} \right)
+
K_p \alpha^{2} n^{\ell} e^{ n  \Delta K_p }  
\left( 1 + \mathbb{E} \left(  \left\Vert  \bar{X}_{0} \right\Vert^{2  \ell - 2 } \right) \right) ,
\end{equation*}
where $\ell = 1, \ldots , p$ and $K_p > 0$.
\end{lemma}

\begin{proof}
Combining \eqref{eq:3.3} with Lemma \ref{lem:AcotacionInicial} yields
\begin{equation*}
 \mathbb{E} \left(  \left\Vert  \bar{X}_{n+1} \right\Vert^q \diagup  \mathfrak{G}_{n} \right)
\leq
\mathbb{E} \left(  \left(  \bar{\rho}_{n+1} \right) ^q \diagup  \mathfrak{G}_{n} \right) 
\leq
 \left\Vert \left( \bar{X}_n , \alpha \right) \right\Vert ^q    \exp \left( \Delta K_q \right).
\end{equation*}
Therefore,
applying the binomial theorem we obtain that for any $q \in \mathbb{N}$,
\begin{equation}
\label{eq:7.11}
 \mathbb{E} \left(  \left\Vert  \bar{X}_{n+1} \right\Vert^{2q} \right)
\leq
e^{ \Delta K_q } \mathbb{E} \left(  \left\Vert  \bar{X}_{n} \right\Vert^{2q} \right)
+
e^{ \Delta K_q } \sum_{j=0}^{q-1}
\begin{pmatrix}
 q \\ j
\end{pmatrix}
 \alpha^{2 \left( q - j \right)}
\mathbb{E} \left(  \left\Vert  \bar{X}_{n} \right\Vert^{2j} \right).
\end{equation}
Without loss of generality we can assume $K_{q+1} \geq K_q \geq 0$.

Iterating \eqref{eq:7.11} we deduce that for all $n \geq 1$,
\begin{align*}
\mathbb{E} \left(  \left\Vert  \bar{X}_{n} \right\Vert^{2q} \right)
& \leq
\exp \left( n  \Delta K_q \right)  \mathbb{E} \left(  \left\Vert  \bar{X}_{0} \right\Vert^{2q} \right)
\\
& \quad 
+ \sum_{j=0}^{q-1}
\begin{pmatrix}
 q \\ j
\end{pmatrix}
\alpha^{2 \left( q - j \right)}
\sum_{k=0}^{n-1}  \exp \left( \left( n - k \right) \Delta K_q \right)
 \mathbb{E} \left(  \left\Vert  \bar{X}_{k} \right\Vert^{2j} \right) .
\end{align*}
Using induction,
together with algebraic manipulations,
we get the assertion of this lemma.
\end{proof}

Next,
applying a stopping time technique (see, e.g., proof of Lemma 3.1 of \cite{TretyakovZhang2013})
we show that the moments of $ \bar{X}_{n}$ are uniformly bounded.

\begin{lemma}
\label{lem:AcotacionMomentos}
Let the assumptions of Lemma \ref{lem:ExpansionSolucionN} be fulfilled.
Consider $p \geq 1 $. 
Then
$
\mathbb{E} \left(  \left\Vert  \bar{X}_{n} \right\Vert^{p} \right)
 \leq
 K \left( T \right)  \left( 1 + \mathbb{E} \left(  \left\Vert  \bar{X}_{0} \right\Vert^{q } \right) \right) 
$
for all $n=0,\ldots, N$ and $N \in \mathbb{N}$.
\end{lemma}

\begin{proof}
Let $p \in \mathbb{N}$.
For any $R>0$ we set
$
\Omega_{R,n} =  \bigcap_{k=0}^n \left( \bar{X}_k \right)^{-1} \left( \left[ -R , R \right] \right)
$.
Using the multinomial theorem gives
\begin{align}
\label{eq:7.18}
 & \mathbb{E} \left(  \mathbf{1}_{\Omega_{R,n}} \left\Vert  \bar{X}_{n+1} \right\Vert^{2p} \right)
  =
 \mathbb{E} \left(  \mathbf{1}_{\Omega_{R,n}} \left\Vert  \left( \bar{X}_{n+1} - \bar{X}_{n} \right) + \bar{X}_{n} \right\Vert^{2p} \right)
\\
\nonumber
& \leq 
 \mathbb{E} \left(  \mathbf{1}_{\Omega_{R,n}} \left\Vert  \bar{X}_{n} \right\Vert^{2p} \right)
 +
 I_{R,n}
 +
 K \sum_{k=3}^{2p} \mathbb{E} \left(  \mathbf{1}_{\Omega_{R,n}} \left\Vert   \bar{X}_{n+1} - \bar{X}_{n} \right\Vert^{k} 
 \left\Vert  \bar{X}_{n} \right\Vert^{2p-k} \right) ,
\end{align}
with
$
 I_{R,n}
 =
 \mathbb{E} \left(  \mathbf{1}_{\Omega_{R,n}} \left\Vert  \bar{X}_{n} \right\Vert^{2p-2}
 \left( p \left( 2p -1 \right)  \left\Vert   \bar{X}_{n+1} - \bar{X}_{n} \right\Vert^{2}
 + 2 p \langle   \bar{X}_{n+1} - \bar{X}_{n} ,    \bar{X}_{n}   \rangle  \right) \right) 
$.

Combining Lemma \ref{lem:LocalAsymptotic} with 
$b, \sigma^k \in \mathcal{C}_{P}^{ 0 } \left( \mathbb{R}^d,\mathbb{R}\right)  $
we obtain
$
\bar{X}_{n+1} - \bar{X}_{n}  = \sqrt{\Delta} \,  \mathcal{O}_{n}
$.
Hence 
\begin{align}
\label{eq:7.17}
& \mathbb{E} \left(  \mathbf{1}_{\Omega_{R,n}} \left\Vert   \bar{X}_{n+1} - \bar{X}_{n} \right\Vert^{k} 
 \left\Vert  \bar{X}_{n} \right\Vert^{2p-k} \right)
 \\ \nonumber
 &  =
 \mathbb{E} \left(  \mathbf{1}_{\Omega_{R,n}} \left\Vert  \bar{X}_{n} \right\Vert^{2p-k}
 \mathbb{E} \left( \left\Vert   \bar{X}_{n+1} - \bar{X}_{n} \right\Vert^{k}  \diagup  \mathfrak{G}_{n}  \right) \right)
 \\ \nonumber
 & \leq 
K  \left(T \right)  \Delta^{k/2} \mathbb{E} \left(  \mathbf{1}_{\Omega_{R,n}} \left\Vert  \bar{X}_{n} \right\Vert^{2p-k}
\left( 1 +  \left\Vert \bar{X}_n \right\Vert^{k q }  \right) \right) 
\hspace{0.8cm} \forall n=0,\ldots, N-1
\end{align}
for any $k=3,\ldots, 2p$.
As $\hat{W}^1_{n+1}, \ldots, \hat{W}^m_{n+1}$ are symmetric random variables
we have
$ \mathbb{E} \left(  \Gamma_{n+1} \diagup  \mathfrak{G}_{n}  \right) = 0$,
and so applying Lemma \ref{lem:LocalAsymptotic} we get
$$
 \mathbb{E} \left( \langle   \bar{X}_{n+1} - \bar{X}_{n} , \bar{X}_{n} \rangle  \diagup  \mathfrak{G}_{n}  \right) 
 \leq 
 \langle  b \left( \bar{X}_{n} \right), \bar{X}_{n} \rangle \Delta
 +
 K \left(T \right) \Delta^{2} \left( 1 +  \left\Vert \bar{X}_n \right\Vert^{ q }  \right) 
$$
for all $n=0,\ldots, N-1$. Similarly 
$$
 \mathbb{E} \left( \langle   \bar{X}_{n+1} - \bar{X}_{n} ,  
 \bar{X}_{n+1} - \bar{X}_{n} \rangle  \diagup  \mathfrak{G}_{n}  \right) 
 \leq 
 \sum_{k=1}^m \left\Vert  \sigma^k \left(  \bar{X}_n \right) \right\Vert ^2 \Delta 
 +
 K \left(T \right) \Delta^{2} \left( 1  +  \left\Vert \bar{X}_n \right\Vert^{ q }  \right)  .
$$
Then, using \eqref{eq:7.1} yields 
\begin{equation}
\label{eq:7.19}
I_{R,n}
\leq  K \left(T \right) \mathbb{E} \left(    \mathbf{1}_{\Omega_{R,n}} \left\Vert  \bar{X}_{n} \right\Vert^{2p-2}
 \left( 
 \Delta \left( 1 +  \left\Vert \bar{X}_n \right\Vert^{ 2 }  \right)
 +
 \Delta^{2} \left( 1 +  \left\Vert \bar{X}_n \right\Vert^{ q }  \right)
 \right) \right) .
\end{equation}
We assume without loss of generality that  $q \geq 2$ in \eqref{eq:7.17} and \eqref{eq:7.19}.
From \eqref{eq:7.18}, \eqref{eq:7.17} and \eqref{eq:7.19} we deduce that
for all $n=0,\ldots, N-1$,
\begin{align}
\label{eq:7.20}
& \mathbb{E} \left(  \mathbf{1}_{\Omega_{R,n}} \left\Vert  \bar{X}_{n+1} \right\Vert^{2p} \right)
   \leq 
 \mathbb{E} \left(  \mathbf{1}_{\Omega_{R,n}} \left\Vert  \bar{X}_{n} \right\Vert^{2p} \right)
 \left( 1 + K \left(T \right)   \Delta \right)
  \\  \nonumber
& 
\hspace{1.5cm} 
+ K \left(T \right)  \Delta \, \mathbb{E} \left(  \mathbf{1}_{\Omega_{R,n}} \left\Vert  \bar{X}_{n} \right\Vert^{2p-2}  \right)
+ K \left(T \right) \Delta^{2} \mathbb{E} \left(   \mathbf{1}_{\Omega_{R,n}} \left\Vert  \bar{X}_{n} \right\Vert^{2p-2+q} \right) 
 \\  \nonumber
 & 
 +  \hspace{-1pt}
 K \left(T \right) \hspace{-2pt} \sum_{k=3}^{2p}  \hspace{-2pt}
  \Delta^{k/2} \mathbb{E} \left(  \mathbf{1}_{\Omega_{R,n}} \left\Vert  \bar{X}_{n} \right\Vert^{2p-k} \right)
 \hspace{-1pt} +  \hspace{-1pt}
K \left(T \right)  \Delta \sum_{k=3}^{2p}  \hspace{-2pt}
  \Delta^{k/2-1} \mathbb{E} \left(  \mathbf{1}_{\Omega_{R,n}} \left\Vert  \bar{X}_{n} \right\Vert^{2p + k \left( q - 1 \right)} \right)  \hspace{-1.5pt}.
\end{align}

In the sequel,  we take
$
R = \Delta^{-1/\left( 6 q - 4\right)}
$.
Then, 
for any for any $ k \geq 3$ we have
\begin{align*}
 \mathbf{1}_{\Omega_{R,n}} \left\Vert  \bar{X}_{n} \right\Vert^{2p + k \left( q - 1 \right)} \Delta^{k/2-1}
 & \leq
 \mathbf{1}_{\Omega_{R,n}} \left\Vert  \bar{X}_{n} \right\Vert^{2p} R^{k \left( q - 1 \right)} \Delta^{k/2-1}
 \leq 
 K \left( T \right) \mathbf{1}_{\Omega_{R,n}} \left\Vert  \bar{X}_{n} \right\Vert^{2p} ,
\end{align*}
because
$
- \frac{ k \left( q-1 \right) }{ 6q-4  } + \frac{k}{2} - 1 > 0
$.
Similarly, using $q \geq 2$ we find
$$
 \mathbf{1}_{\Omega_{R,n}} \left\Vert  \bar{X}_{n} \right\Vert^{2p-2+q} \Delta
 \leq
 \mathbf{1}_{\Omega_{R,n}} \left\Vert  \bar{X}_{n} \right\Vert^{2p} R^{q-2} \Delta
 \leq 
 K \left( T \right) \mathbf{1}_{\Omega_{R,n}} \left\Vert  \bar{X}_{n} \right\Vert^{2p}  .
$$
Therefore, employing \eqref{eq:7.20} we obtain that for all $n=0,\ldots, N-1$,
\begin{align*}
 \mathbb{E} \left(  \mathbf{1}_{\Omega_{R,n}} \left\Vert  \bar{X}_{n+1} \right\Vert^{2p} \right)
 &  \leq 
 \mathbb{E} \left(  \mathbf{1}_{\Omega_{R,n}} \left\Vert  \bar{X}_{n} \right\Vert^{2p} \right)
 \left( 1 + K \left(T \right)   \Delta \right)
  \\ 
 & \quad
 +
 K \left(T \right)  \sum_{k=1}^{p} 
  \Delta^{k} \mathbb{E} \left(  \mathbf{1}_{\Omega_{R,n}} \left\Vert  \bar{X}_{n} \right\Vert^{2p-2k} \right) .
\end{align*}
By Young's inequality,
$
 \left\Vert  \bar{X}_{n} \right\Vert^{2p-2k} \Delta^{k-1}
 \leq
 \frac{p-k}{p} \left\Vert  \bar{X}_{n} \right\Vert^{2p} + \frac{k}{p} \Delta^{p \left( k-1 \right) / k}
$
whenever $k=1, \ldots, p-1$.
In consequence, for any $n=0,\ldots, N-1$,
\begin{equation}
\label{eq:7.21}
 \mathbb{E} \left(  \mathbf{1}_{\Omega_{R,n}} \left\Vert  \bar{X}_{n+1} \right\Vert^{2p} \right)
 \leq 
 \mathbb{E} \left(  \mathbf{1}_{\Omega_{R,n}} \left\Vert  \bar{X}_{n} \right\Vert^{2p} \right)
 \left( 1 + K \left(T \right)   \Delta \right)
 + K \left(T \right)   \Delta .
\end{equation}

Combining \eqref{eq:7.21} with  
$
\mathbb{E} \left(  \mathbf{1}_{\Omega_{R,n+1}} \left\Vert  \bar{X}_{n+1} \right\Vert^{2p} \right)
 \leq  
 \mathbb{E} \left(  \mathbf{1}_{\Omega_{R,n}} \left\Vert  \bar{X}_{n+1} \right\Vert^{2p} \right)
$
yields
$$
\mathbb{E} \left(  \mathbf{1}_{\Omega_{R,n+1}} \left\Vert  \bar{X}_{n+1} \right\Vert^{2p} \right)
\leq 
 \mathbb{E} \left(  \mathbf{1}_{\Omega_{R,n}} \left\Vert  \bar{X}_{n} \right\Vert^{2p} \right)
 \left( 1 + K \left(T \right)   \Delta \right)
 + K \left(T \right)   \Delta 
$$
for all $n=0,\ldots, N-1$.
Now, using Gronwall's lemma we deduce that 
\begin{equation}
\label{eq:7.23}
 \mathbb{E} \left(  \mathbf{1}_{\Omega_{R,n}} \left\Vert  \bar{X}_{n} \right\Vert^{2p} \right)
 \leq
 K \left(T \right) \left( 1 +  \mathbb{E} \left(  \left\Vert  \bar{X}_{0} \right\Vert^{2p} \right)   \right)
 \hspace{1cm}
 \forall n=0,\ldots, N ,
\end{equation}
and so \eqref{eq:7.21} leads to
\begin{equation}
\label{eq:7.22}
 \mathbb{E} \left(  \mathbf{1}_{\Omega_{R,n}} \left\Vert  \bar{X}_{n+1} \right\Vert^{2p} \right)
 \leq
 K \left(T \right) \left( 1 +  \mathbb{E} \left(  \left\Vert  \bar{X}_{0} \right\Vert^{2p} \right)   \right)
 \hspace{0.8cm}
 \forall n=0,\ldots, N -1.
\end{equation}

Since
$
\mathbf{1}_{\left(\Omega_{R,n} \right)^c}
=
\mathbf{1}_{ \left\{ \left\Vert  \bar{X}_{0} \right\Vert > R \right\} }
+ \sum_{k=1}^n \mathbf{1}_{\Omega_{R,k-1}} \mathbf{1}_{ \left\{ \left\Vert  \bar{X}_{k} \right\Vert > R \right\} }
$,
the Cauchy-Schwarz inequality shows that
\begin{align*}
 \mathbb{E} \left( \mathbf{1}_{\left(\Omega_{R,n} \right)^c} \left\Vert  \bar{X}_{n} \right\Vert^{2p} \right)
 & \leq
 \sqrt{ \mathbb{E} \left( \left\Vert  \bar{X}_{n} \right\Vert^{4p} \right)} 
 \sqrt{ \mathbb{E} \left( \mathbf{1}_{ \left\{ \left\Vert  \bar{X}_{0} \right\Vert > R \right\} } \right) }
 \\
 & \quad
 + \sum_{k=1}^n 
  \sqrt{ \mathbb{E} \left( \left\Vert  \bar{X}_{n} \right\Vert^{4p} \right)}  
   \sqrt{ \mathbb{E} \left(
  \mathbf{1}_{\Omega_{R,k-1}} \mathbf{1}_{ \left\{ \left\Vert  \bar{X}_{k} \right\Vert > R \right\} }
  \right) } .
\end{align*}
By 
$
R = \Delta^{-1/\left( 6 q - 4\right)}
$,
applying Markov's inequality we obtain
\begin{align*}
\mathbb{E} \left( \mathbf{1}_{\left(\Omega_{R,n} \right)^c} \left\Vert  \bar{X}_{n} \right\Vert^{2p} \right)
& \leq
 \sqrt{ \mathbb{E} \left( \left\Vert  \bar{X}_{n} \right\Vert^{4p} \right)} 
 \sqrt{ \mathbb{E} \left( \left\Vert  \bar{X}_{0} \right\Vert ^{2 \left(p+1 \right) \left( 6 q -4\right)} \right) }
 \Delta^{p+1}
 \\
 & \quad
 + \sqrt{ \mathbb{E} \left( \left\Vert  \bar{X}_{n} \right\Vert^{4p} \right)}   \sum_{k=1}^n 
   \sqrt{ \mathbb{E} \left(
  \mathbf{1}_{\Omega_{R,k-1}} \left\Vert  \bar{X}_{k} \right\Vert ^{2 \left(p+1 \right) \left( 6 q -4\right)}
  \right) }   \Delta^{p+1} .
\end{align*}
From \eqref{eq:7.22} it follows that
$$
 \sum_{k=1}^n 
   \sqrt{ \mathbb{E} \left(
  \mathbf{1}_{\Omega_{R,k-1}} \left\Vert  \bar{X}_{k} \right\Vert ^{2 \left(p+1 \right) \left( 6 q -4\right)}
  \right) }  
\leq
n K \left(T \right) \sqrt{ 1 +  \mathbb{E} \left(  \left\Vert  \bar{X}_{0} \right\Vert^{2 \left(p+1 \right) \left( 6 q -4\right)} \right) } .
$$
Lemma \ref{lem:AcotacionMomentosCruda} implies 
\begin{equation*}
\mathbb{E} \left(  \left\Vert  \bar{X}_{n} \right\Vert^{4 p} \right)
 \leq
 K \left(T \right) \left( \mathbb{E} \left(  \left\Vert  \bar{X}_{0} \right\Vert^{4  p} \right)
+
n^{2p}  \left( 1 + \mathbb{E} \left(  \left\Vert  \bar{X}_{0} \right\Vert^{4 p  } \right) \right) 
\right) .
\end{equation*}
Therefore, for all $n=0,\ldots, N$ we have
\begin{align*}
 \mathbb{E} \left( \mathbf{1}_{\left(\Omega_{R,n} \right)^c} \left\Vert  \bar{X}_{n} \right\Vert^{2p} \right)
  \leq
   \left( n +1 \right) \left( T / N \right)^{p+1} K \left(T \right)   \sqrt{ 1 +  \mathbb{E} \left(  \left\Vert  \bar{X}_{0} \right\Vert^{2 \left(p+1 \right) \left( 6 q -4\right)} \right) }  \hspace{10pt} &
 \\
. \sqrt{ n^{2p} + \left( 1 + n^{2p} \right) \mathbb{E} \left( \left\Vert  \bar{X}_{0} \right\Vert^{4p} \right)} & ,
\end{align*}
hence
$
\mathbb{E} \left( \mathbf{1}_{\left(\Omega_{R,n} \right)^c} \left\Vert  \bar{X}_{n} \right\Vert^{2p} \right)
\leq
 K \left( T \right)  \left( 1 + \mathbb{E} \left(  \left\Vert  \bar{X}_{0} \right\Vert^{q } \right) \right) 
$,
and so  \eqref{eq:7.23} leads to
\begin{equation}
\label{eq:7.24}
 \mathbb{E} \left(  \left\Vert  \bar{X}_{n} \right\Vert^{2p} \right)
 \leq
 K \left( T \right)  \left( 1 + \mathbb{E} \left(  \left\Vert  \bar{X}_{0} \right\Vert^{q } \right) \right) 
 \hspace{1cm}
 \forall n=0,\ldots, N .
\end{equation}
Finally, 
using Jensen's inequality we extend \eqref{eq:7.24} to any $p \geq 1$.
\end{proof}

We are now in a position to show Theorem \ref{th:RateConvergence}.

\begin{proof}[of Theorem \ref{th:RateConvergence}] 
Let $\bar{X}_{n}$, $\bar{\eta}_{n}$ and $\bar{V}_{n}$ be as 
in the first paragraph of this subsection.
In case $\alpha \neq 0 $,
$\left( \bar{U}_{n+1}, \bar{V}_{n+1}  \right)$ is
the vector $\bar{Z}_{n+1}$ arising from one iteration of  Scheme \ref{scheme:EulerStable} applied to
\begin{align*}
 \begin{pmatrix}
 \widetilde{X_t} 
 \\  \widetilde{V}_t
\end{pmatrix}  
& = 
\begin{pmatrix}
 \bar{X}_n \\ \alpha
\end{pmatrix} 
+
\int_{T_n}^t 
\begin{pmatrix}
 b\left(  \widetilde{X}_s \right) - b\left( 0 \right) +   \frac{b \left( 0 \right) }{\alpha} \widetilde{V}_s 
 \\ 0
\end{pmatrix} 
ds
\\
& \quad
+ \sum_{k=1}^{m}\int_{T_n}^{t}  
\begin{pmatrix}
 \sigma^{k} \left(  \widetilde{X}_s \right) - \sigma^{k}\left( 0 \right) +   \frac{\sigma^{k} \left( 0 \right) }{\alpha} \widetilde{V}_s 
 \\ 0
\end{pmatrix} 
dW^{k}_{s} 
\end{align*}
(see Remark \ref{rem:ObtencionEulerStableG}).
Since the law of $\hat{W}^{k}_{n}$ 
is absolutely continuous with respect to the Lebesgue measure,
using  Lemma \ref{lem:ZnoNulo} yields
$\left( \bar{U}_{n}, \bar{V}_{n}  \right) \neq 0$ a.s.
for all $n\geq 1$,
in both situations $\alpha = 0$ and $\alpha \neq 0$.
Therefore,
the assumptions of Lemma \ref{lem:ExpansionSolucionN} are fulfilled in our framework.
Lemma \ref{lem:AcotacionMomentos} yields Condition (ii).
Applying Lemma \ref{lem:LocalAsymptotic} we deduce Conditions (iii) and (iv).
Then, Theorem \ref{th:GeneralWeakConvergence} leads to
$
\left\vert \mathbb{E} \varphi \left(X_{T}\right) -  \mathbb{E} \varphi \left(  \bar{X}_{N}  \right)   \right\vert 
\leq K \left(  1+\mathbb{E}\left\Vert X_{0} \right\Vert ^{q}\right)  T/N 
$
for all $N \in \mathbb{N}$.
\end{proof}

\subsection{Proof of Theorem \ref{th:StrongRateConvergence}}
Proceeding as in the derivation of the mean-square order of convergence 
of the Euler-Maruyama scheme in the global Lipschitz case,
we obtain 
\begin{equation}
\label{SEq1}
\left\Vert \mathbb{E}
\left( X_{T_{n+1}} - \bar{E}_{n+1} \left( \bar{X}_{n} \right)
 \diagup
\mathcal{\mathfrak{F}}_{T_n} \right)\right\Vert 
\leq K \left(1+\left\Vert \bar{X}_{n}\right\Vert ^{q}\right) \Delta^{3/2}
\end{equation}
\begin{equation}
\label{SEq2}
\text{and } \hspace{1cm}
\mathbb{E}\left(\left\Vert 
X_{T_{n+1}} - \bar{E}_{n+1} \left( \bar{X}_{n} \right) \right\Vert ^{2p}
 \diagup
 \mathcal{\mathfrak{F}}_{T_n} \right)
\leq 
K\left(1+\left\Vert \bar{X}_{n}\right\Vert ^{q }\right)\Delta^{2p} ,
\end{equation} 
where 
$ n = 0, \ldots, N-1$,
$
\bar{E}_{n+1} \left( x \right) 
= 
x + \Delta b( x )
+ \sum_{k=1}^{m}\sigma^{k} (x) \left( W^k_{T_{n+1}} -  W^k_{T_{n}} \right)
$,
and $p \in \mathbb{N}$
(see, e.g., Appendix C of \cite{TretyakovZhang2013}).
By Lemma \ref{lem:LocalAsymptotic},
 \begin{align*}
\bar{X}_{n+1} -  \bar{E}_{n+1} \left( \bar{X}_{n} \right)
& =
\Delta  \sum_{k=1}^m \left( 
\left\langle \hat{U}_n  ,\bar{ \sigma }_n^k  \right\rangle 
 \bar{\eta}_{n} \bar{ \sigma }_n^k - \frac{1}{2} \left\Vert \bar{ \sigma }_n^k  \right\Vert^2  \bar{X}_{n} 
\right)
\left( \left( \hat{W}^k_{n+1} \right)^2 - 1 \right)
\\  
& \quad
+ \Delta   \sum_{k \neq j} \left(
 \left\langle \hat{U}_n  ,\bar{ \sigma }_n^k  \right\rangle 
\bar{\eta}_{n} \bar{ \sigma }_n^j - \frac{1}{2} \left\langle \bar{ \sigma }_n^k  ,\bar{ \sigma }_n^j  \right\rangle  \bar{X}_{n} 
\right) \hat{W}^k_{n+1} \hat{W}^j_{n+1}
\\  
& \quad
+ \Delta^{3/2}  \bar{\eta}_{n}  \Gamma_{n+1} + \Delta^2 \mathcal{O}_{n+1} .
\end{align*}
Hence 
$
\left\Vert \mathbb{E} \left( \bar{X}_{n+1} - \bar{E}_{n+1} \left( \bar{X}_{n} \right)
 \diagup
\mathcal{\mathfrak{F}}_{T_n} \right)\right\Vert 
\leq K \left(1+\left\Vert \bar{X}_{n}\right\Vert ^{q}\right) \Delta^{2}
$
and 
\begin{equation*}
\mathbb{E}\left(\left\Vert 
 \bar{X}_{n+1} - \bar{E}_{n+1} \left( \bar{X}_{n} \right) \right\Vert ^{2p}
 \diagup
 \mathcal{\mathfrak{F}}_{T_n} \right)
\leq 
K\left(1+\left\Vert \bar{X}_{n}\right\Vert ^{q }\right)\Delta^{2p} 
\end{equation*} 
provided that $ n = 0, \ldots, N-1$.
Therefore,
\eqref{SEq1} and \eqref{SEq2} yield 
$$
\left\Vert \mathbb{E}
\left( X_{T_{n+1}} - \bar{X}_{n+1} 
 \diagup
\mathcal{\mathfrak{F}}_{T_n} \right)\right\Vert 
\leq K \left(1+\left\Vert \bar{X}_{n}\right\Vert ^{q}\right) \Delta^{3/2}
$$
and 
$\mathbb{E}\left(\left\Vert 
X_{T_{n+1}} - \bar{X}_{n+1} \right\Vert ^{2p}
 \diagup
 \mathcal{\mathfrak{F}}_{T_n} \right)
\leq 
K\left(1+\left\Vert \bar{X}_{n}\right\Vert ^{ q}\right)\Delta^{2p} 
$
for all $ n \leq N-1$.
Now,
using Lemma \ref{lem:AcotacionMomentos} and Theorem 2.1 of \cite{TretyakovZhang2013} 
we get \eqref{eq:SEq3}.
\endproof

\section{Conclusions}
\label{sec:Conclusions}

In this paper a new approach for designing numerical schemes for SDEs is introduced on the basis of the integration of the system of coupled SDEs that describes the evolution of the norm of the required solution and its projection on the unit sphere.  
In this manner, a numerical scheme for the general class of SDEs with multiplicative noise is proposed and three optimized variants for specific subclasses of SDEs are also introduced. 
Under general conditions, it is proved that the new scheme preserves the almost sure stability of the solutions for any step-size, as well as the property of being distant from $0$. 
Under local Lipschitz conditions, the order $0.5$ of strong convergence of the scheme is obtained 
as well as the order $1$ of weak convergence. 
This last particular result is obtained from the previously derived fundamental weak convergence 
Theorem \ref{th:GeneralWeakConvergence} 
that now complements the  fundamental mean-square convergence theorem for SDEs with locally Lipschitz coefficients recently stated in \cite{TretyakovZhang2013}. 
The proposed scheme is compared in simulations with various state of the art numerical schemes with comparable stability properties. 
The results clearly demonstrate the advantages of the proposed scheme in the integration of a variety of equations with different dynamics and stability properties.

\section*{Acknowledgments}

The authors thanks to the editor Prof. Clayton Webster and referees for their valuable comments and suggestions on the manuscript.


\end{document}